\documentclass[11pt,a4paper]{article}
\pdfoutput=1
\usepackage{latexsym}
\usepackage{amssymb}
\usepackage{amsfonts}
\usepackage{amsmath}
\usepackage{amsthm}
\usepackage{ulem}	%allows esoteric underlinings
\usepackage{sectsty}	%sectionstyle
\usepackage[all,knot]{xy}
\usepackage{longtable}
\usepackage{lscape}
\usepackage[section]{placeins}
\usepackage{graphics}
\usepackage{psfrag, psfig}
\usepackage[dvips]{graphicx}
\usepackage{caption}
\usepackage{subcaption}
\usepackage{rotating} %allows figures and tables to be used in landscape mode
\usepackage{epsfig}
\title{Braid Equivalence in the H\'enon Family I}
\author{A.\ de Carvalho \and T.\ Hall \and P.\ Hazard}
\date{\today}

\theoremstyle{plain}%default
\newtheorem{thm}{Theorem}[section]

\newtheorem{cor}[thm]{Corollary}

\theoremstyle{remark}
\newtheorem{rmk}[thm]{Remark}

\newtheorem{obs}[thm]{Observation}

\theoremstyle{definition}
\newtheorem{defn}[thm]{Definition}
\newtheorem{eg}[thm]{Example}
\newtheorem{constr}[thm]{Construction}

\numberwithin{equation}{section}

\DeclareMathOperator{\Homeo}{\ensuremath{\mathrm{Homeo}}}

\newcommand{\bt}{\ensuremath{\mathrm{bt}}}

\renewcommand{\o}{\ensuremath{\mathrm{o}}}
\newcommand{\init}{\ensuremath{\mathrm{{init}}}}
\renewcommand{\term}{\ensuremath{\mathrm{term}}}
\newcommand{\Disk}{\ensuremath{\Delta}}
\newcommand{\RR}{\ensuremath{\mathbb{R}}}
\newcommand{\NN}{\ensuremath{\mathbb{N}}}
\newcommand{\ZZ}{\ensuremath{\mathbb{Z}}}
\newcommand{\<}{\ensuremath{\langle}}
\renewcommand{\>}{\ensuremath{\rangle}}
\newcommand{\del}{\ensuremath{\partial}}
\newcommand{\tr}{\ensuremath{\mathrm{tr}}}

\begin{document}

\newpage
\pagestyle{plain}

\setcounter{page}{1}

\maketitle
%\tableofcontents
\begin{abstract}
We give two general constructions of braid equivalences which exist between certain deformations of the 2-branched Horsehoe map.
We then give numerical evidence suggesting that these constructions of braid equivalences are always realised in the H{\'e}non family.
\end{abstract}

\section{Introduction}

During the last three decades of the twentieth century, much effort
was devoted to the study of families of low-dimensional dynamical
systems depending on parameters. There is today a very thorough theory
explaining the dynamics of families of one-dimensional (real and complex) 
endomorphisms.  The dynamics of the real quadratic family
\mbox{$f_{a}(x)=a-x^{2}$}, for example, is nearly completely
understood~\cite{Lyubich00}. In the 1970s H\'enon introduced the
family\footnote{H{\'e}non considered a different parametrisation 
of this family,
\begin{equation*}
H_{a,b}(x,y)=(-y+1-ax^2,bx)
\end{equation*}} 
which now bears his name, a two-dimensional analog
of the quadratic family: $F_{a,b}(x,y)=\left(f_{a}(x)-by,x\right)$.
This is a family of plane diffeomorphisms depending on two parameters
which, for $b=0$, degenerates to the quadratic family.  In contrast to
the quadratic family, and despite the existence of several beautiful
results about it, our understanding of the H\'enon family is still
rather rudimentary.  While there are many similarities between the two
families which make it possible to use knowledge of the former to help
in understanding the latter, there are also many fundamental
differences which demand that different techniques be developed.

One of the most basic aspects of the dynamics of a parametrized
family is the way in which the periodic orbit structure changes as the
parameters vary. This article is concerned with periodic orbits in
H\'enon family in the parameter regions close to degeneration, and
exploits both similarities and differences between the quadratic and
H\'enon families.

\medskip

Periodic orbits of endomorphisms of the real line are specified by
their associated cyclic permutation: the way that their points,
ordered on the line, are permuted by the endomorphism. For
homeomorphisms of the plane such as H\'enon maps, the analogous
specification, introduced by Boyland~\cite{Boyland84, Boyland94}, is
the braid type. If $F\colon{\mathbb R}^2\to{\mathbb R}^2$ is an 
orientation-preserving homeomorphism and~$P$ is a periodic orbit of~$F$, 
then the braid type $\bt(P,F)$ is the isotopy class of $F$ relative to~$P$, 
up to topological conjugacy. In other words, the braid type of~$P$ is
determined by fixing the action of~$F$ on~$P$ but allowing it to be
deformed by isotopy in the complement of~$P$, and also allowing a
global change of coordinates.

The periodic orbit structure of maps in the quadratic family -- or,
indeed, of any unimodal map~$f$ -- is easily understood using
techniques of kneading theory. The critical point~$c$ is used to
divide the line into left and right halves, and the kneading sequence
of~$f$ is the itinerary of the critical value~$f(c)$: the sequence of
lefts and rights along the orbit $(f^n(c))_{n\ge 1}$. There is then a
simple recipe for generating the set of all itineraries of
points~$x\in{\mathbb R}$ from this kneading sequence. Permutations
associated to periodic orbits of unimodal maps -- which are called
{\it unimodal permutations} -- are determined by the itineraries of the
points on the orbit. The set of permutations of periodic orbits of a
unimodal map~$f$ is therefore determined by the kneading sequence of
the map, and can be enumerated by a straightforward algorithm.

The situation for the H\'enon family is quite different. We have very
little idea, to this day, of the way in which braid types of periodic
orbits are built up in the family, going from none to a full
horseshoe's worth, as the parameter~$a$ increases. In fact, by the
result of Kan, Ko\c{c}ak and Yorke~\cite{KKY92}, periodic orbits of the
H\'enon family are both created and destroyed near every homoclinic
tangency, and it is not even known whether or not all periodic orbits
which appear in the H\'enon family have the same braid type as
periodic orbits of the horseshoe.

\medskip

\begin{figure}[h]
\begin{subfigure}[b]{0.5\textwidth}
\centering
\includegraphics
%[bb=0 0 320 320]
[scale=0.1]
{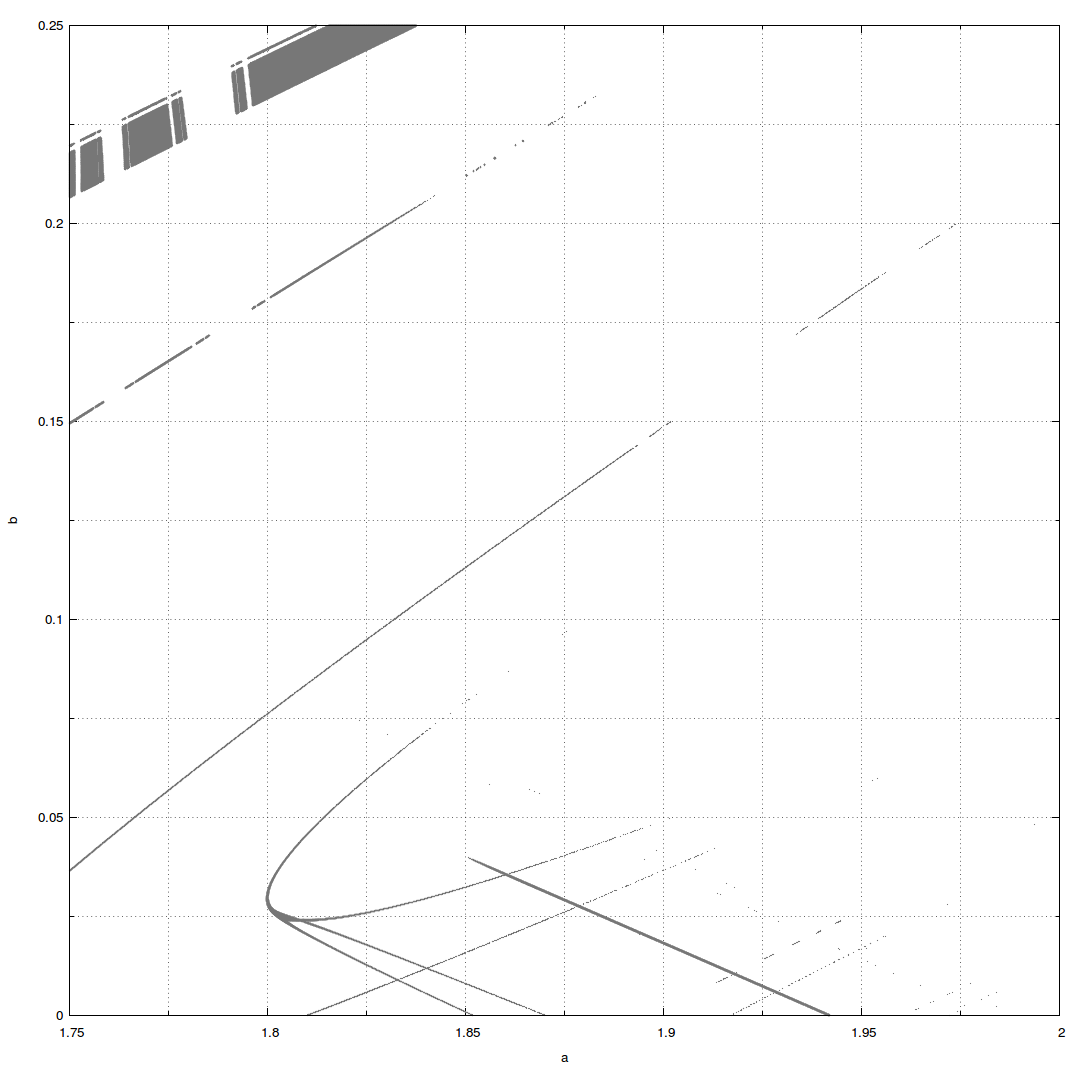}
\caption{Period $p=8$}
\label{fig:per8}
\end{subfigure}
~
\begin{subfigure}[b]{0.5\textwidth}
\centering
\includegraphics
%[bb=0 0 320 320]
[scale=0.1]
{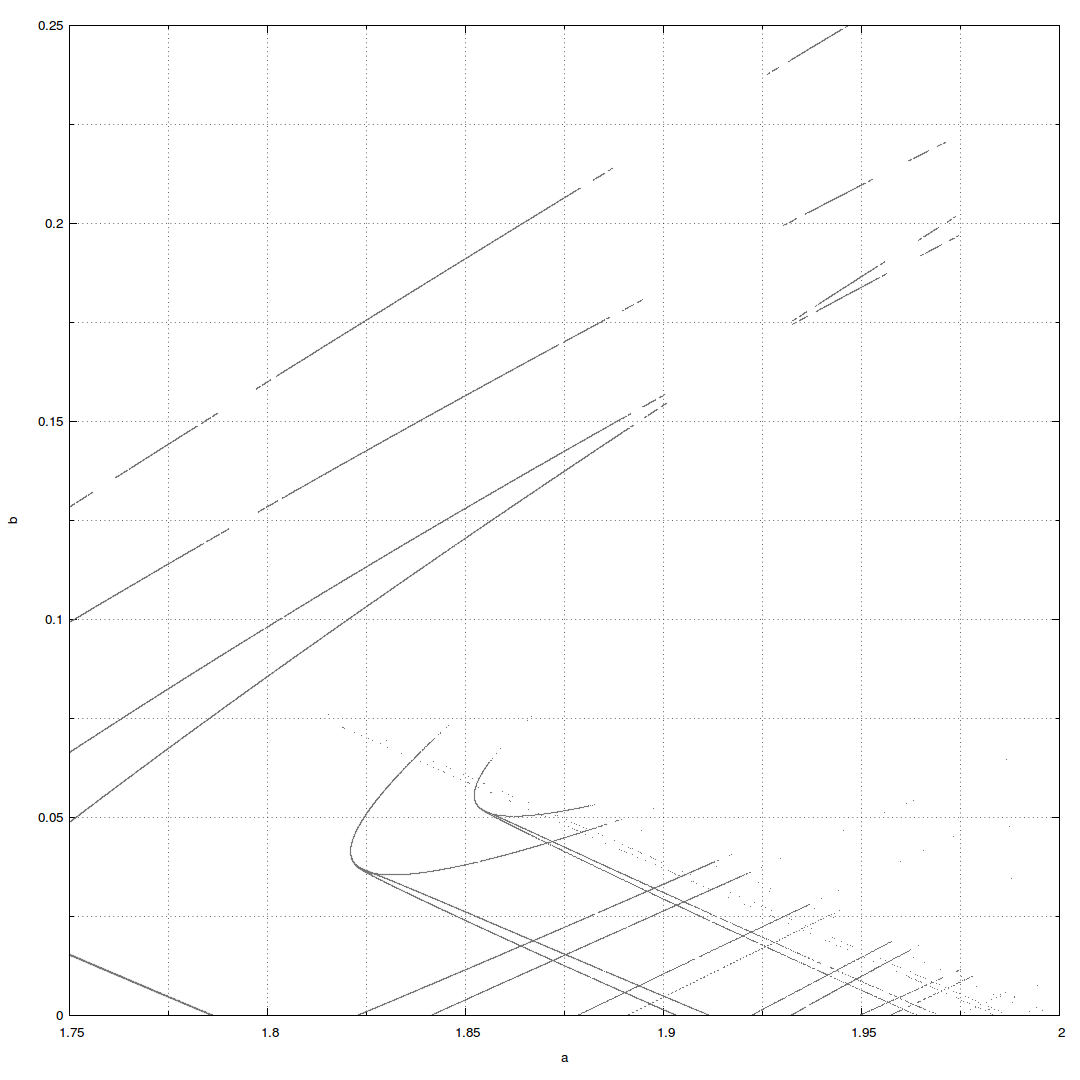}
\caption{Period $p=9$}
\label{fig:per9}
\end{subfigure}

\begin{subfigure}[b]{0.5\textwidth}
\centering
\includegraphics
%[bb=0 0 320 320]
[scale=0.1]
{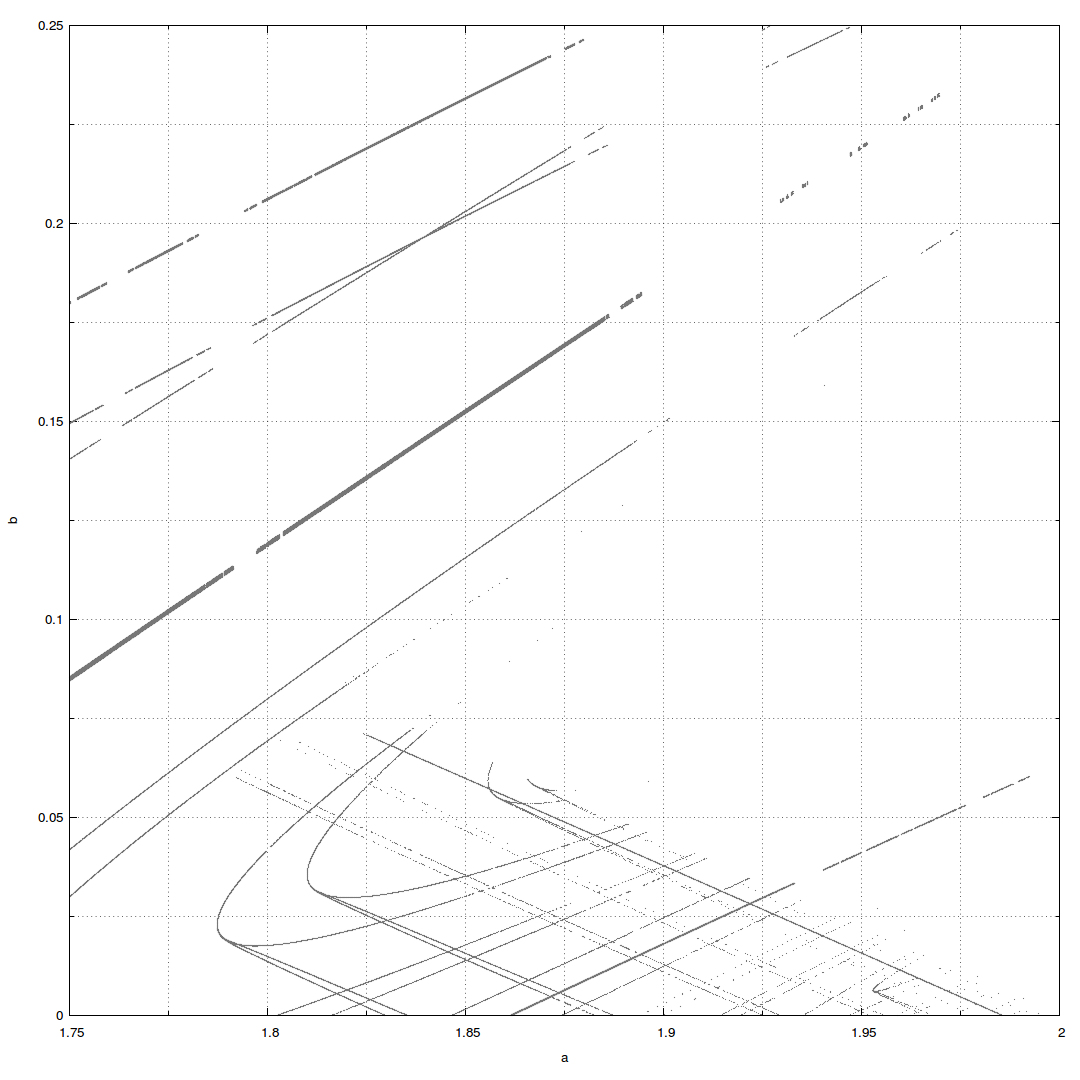}
\caption{Period $p=10$}
\label{fig:per10}
\end{subfigure}
~
\begin{subfigure}[b]{0.5\textwidth}
\centering
\includegraphics
%[bb=0 0 320 320]
[scale=0.1]
{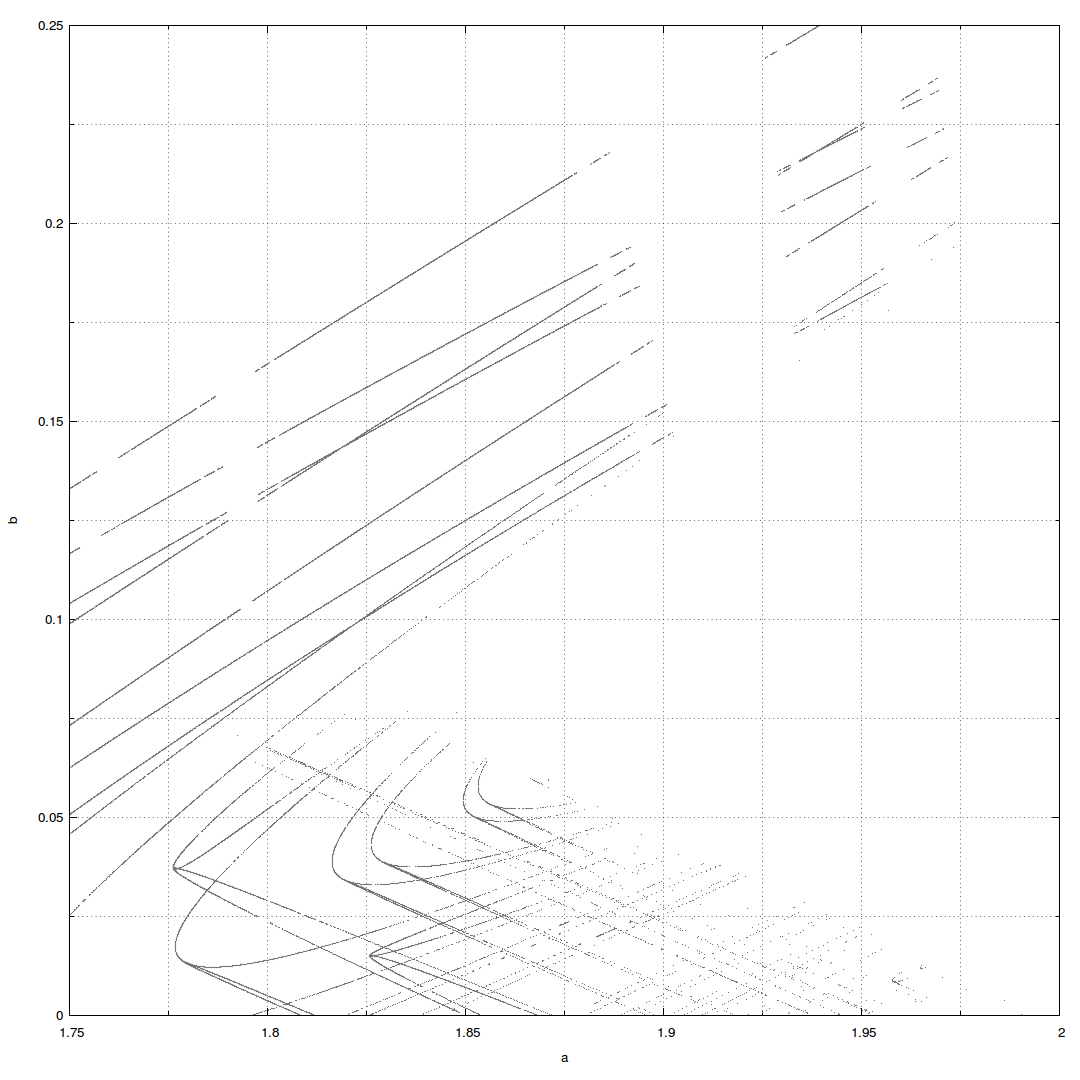}
\caption{Period $p=11$}
\label{fig:per11}
\end{subfigure}
\caption{Scatterplots of parameters $(a,b)\in [1.75,2.00]\times
  [0.0,0.25]$ for which the H\'enon map $F_{a,b}$ possesses an
  attracting periodic orbit of period $p=8,9,10,11$.  One of the
  ``hook-like'' structures for period 8 can be seen intersecting the
  $a$-axis in two intervals at approximately $a=1.8517$ and $a=1.87$
  respectively.  }
\label{fig:per8-11}
\end{figure}

The diagrams in Figure~\ref{fig:per8-11} show regions in the parameter
plane for the H\'enon family where attracting periodic orbits of
periods 8, 9, 10, and 11 were found. 
(Similar loci, for various periods, were first considered by El Hamouly 
and Mira~\cite{ElHamoulyMira82}.) These plots have a very rich
structure, and understanding how the dynamics varies in the H\'enon
family includes explaining this structure.  In this paper we are
particularly interested in the hook-like structures, also called
{\it swallow configurations} by Milnor~\cite{Milnor92} in the
one-dimensional cubic case.  These structures are open sets consisting
of a main body and four limbs, two of which intersect the $a$-axis
in two distinct (small) intervals.
Each of these hooks indicates that there is one attracting periodic
orbit in the H\'enon family which can be deformed into two different
attracting periodic orbits of the quadratic family. That is, we expect 
each of the hooks to be associated to attracting periodic orbits of the
corresponding H\'enon maps whose braid type is constant in the
region~$b>0$ and which degenerate into periodic orbits of the
quadratic family with two different permutations as $b\downarrow 0$
along each of the two ends of the hook. Viewing this in the opposite
direction, the hooks indicate that there are certain pairs of unimodal
permutations which coalesce into a single H\'enon periodic orbit.

In this paper we identify mechanisms which lead to this coalescence on
the level of unimodal permutations.  More precisely, we describe two
mechanisms which associate to unimodal permutations pairs of
equivalent braids, and provide numerical evidence showing
that the hooks in Figure~\ref{fig:per8-11} can all be explained in
terms of these mechanisms. 

There are at least two distinct reasons for which it is useful to be
able to relate unimodal permutations to braid equivalences in the
H\'enon family. First, as mentioned above, our understanding of the
periodic orbit structure of H\'enon maps is very limited, and being
able to connect H\'enon braid types to unimodal permutations provides
a means of deducing information about the former from the latter
(about which we know everything); one of the consequences of our
results is that we predict the existence of infinitely many hooks in
the H\'enon family, each of which associates a braid type with a pair
of intervals in the parameter space of the quadratic family. Second,
the more general problem of deciding whether or not two braid types
are equal is a difficult one, and mechanisms for constructing equivalent
braids can be useful. For example, there are conjectural
constraints~\cite{dCH3} on the possible orders in which horseshoe braid
types can be built up in families which pass from trivial dynamics to
a full horseshoe; and if this conjecture holds, then each pair of
equivalent braids generates an infinite family of pairs of equivalent
braids. 

\bigskip

We next review some definitions and terminology which will make it
possible to give rough descriptions of the two mechanisms mentioned
above. A geometric braid on $n$ strands (see Figure~\ref{fig:braid1-})
is a diagram with~$n$ arcs (strands) connecting two ordered sets
of~$n$ points lined up vertically, so that only double intersections
are allowed and at each of them it is specified which strand goes
above and which goes below. The points at the top are called initial
endpoints and those at the bottom are called terminal endpoints. To
each geometric braid is associated a braid type, and braid types
determine geometric braids up to conjugacy (these facts will be
discussed further in the text). Geometric braids also induce
permutations on~$n$ elements in the obvious way: associate to the
initial endpoint the terminal endpoint along the same strand. Since
this association forgets all information about crossings of strands,
it is far from being one-to-one. It is possible, however, to associate
a unique braid to a unimodal permutation by requiring that each pair
of strands crosses at most once and that, when two strands cross, the
one which started to the left goes above the one which started to the
right. In this way we can talk about the unimodal braid associated
with a unimodal permutation.

 Given a unimodal permutation $\upsilon$, let $f$ be a unimodal map
 realising $\upsilon$ as its critical orbit.  The {\it dynamical}
 preimage of a point in the critical orbit is the preimage under $f$
 which is also contained in the critical orbit (corresponding to the
 unique preimage at the level of the permutation).  The other preimage
 of a point of the critical orbit is called the {\it non-dynamical}
 preimage.

The first mechanism is as follows.  First we `break' $f$ at the
dynamical preimage of the critical point: that is, we perturb $f$ in a
neighbourhood of this point.  Assuming that the break is small, we get a new
point which is a closest return to the critical point.  Make this new
point follow the sequence of forward iterates of the critical point:
we may do this for any finite time by making the break sufficiently
small. If we arrive near the non-dynamical preimage of the critical
point, we can `reconnect' this iterate of the new point to the critical
point. Again, this means that we perturb $f$ in a neighborhood of this
iterate so that its image is the critical point. Depending on which
side of the critical point the new point lies, we get a pair of
distinct critical orbit types.  (In terms of braids this construction
is a generalisation of the cabling construction for a braid, followed
by a half-twist between strands.)

The second mechanism takes a pair of equivalent unimodal permutations
$\upsilon_-$ and $\upsilon_+$, such as those generated by the first
mechanism, and produces another equivalent pair $\upsilon_-^1$ and
$\upsilon_+^1$ of unimodal permutations.  Let $f_-$ and $f_+$ denote
unimodal maps whose critical orbits realise $\upsilon_-$ and
$\upsilon_+$ respectively.  As before, we break $f_-$ and $f_+$ at the
dynamical preimage, but in both cases the break is large so that the
new point is a second closest return to the critical point.  This is
done to ensure that the image of the new point is close to the image
of the first closest return.  We make this image of the new point
follow the forward iterates of the first closest return until it
arrives close to the dynamical preimage of the critical point. It is
then reconnected to the critical point as before. By repeating this
process, a chain of pairs of equivalent braids can be generated.

\bigskip
In Section 2 we present some background on braids and unimodal
maps; in Sections~3 and~4 the first and second constructions of braid
equivalences are described; and in Section~5 we discuss applications
to the H\'enon family and numerical results.

\paragraph{Acknowledgements.}
A.~de Carvalho was partially supported by {\sc FAPESP} grant 2011/16265-8. 
T.~Hall was partially supported by {\sc FAPESP} grant 2011/17581-0. 
P.~Hazard was supported by {\sc FAPESP} grant 2008/10659-1 and Leverhulme Trust grant RPG-279. 
We would like to thank the institutions {\sc IME-USP}, {\sc IMPA} and the {\sc IMS} at Stony Brook for their hospitality during the time in which
parts of this work was done. Thanks also to A.~Hammerlindl for his help with numerical computations.
Finally, we would like to thank C.~Tresser for useful discussions on braids and the H\'enon family.

\section{Notation and Terminology}

%p=period
%P=periodic orbit
\subsection{Braids and Braid Equivalence}\label{subsect:braids}
Let $\Disk\subset \RR^2$ denote the closed unit disk.
Let $\Homeo_+(\Disk)$ denote the group of orientation-preserving homeomorphisms of $\Disk$. 
Given a subset $A\subset \Disk$ let $\Homeo_+(\Disk,A)$ denote the group of orientation-preserving homeomorphisms $f$ of $\Disk$ satisfying $f(A)=A$.
%\verb=in some def's bd D is fixed=
Where necessary we endow these groups with the uniform topology.
%(We consider this currently for simplicity. Later, when considering H\'enon maps, we will need to consider diffeomorphisms onto their images or diffeomorphisms of $\RR^2$, but the changes needed in those cases are minor.)

\subsubsection{Braid Equivalence}
%\verb=DEFINE BRAIDS=
%\verb=DEFINE BRAID EQUIVALENCE=
%\verb=SHOW THAT PERIODIC ORBS INDUCE BRAIDS=
%\verb=SHOW THAT EQUIV BRAIDS GIVE BRAID EQUIVALENT ORBITS=

%this def used in Hall94
Let $p\in \NN$.
Let $P_{\mathrm{can}}\subset \Delta$ denote the unique set of $p$ points contained in the horizontal axis, whose complement in this axis has connected components all of equal length.
%i.e. $P_\mathrm{can}$ is given by
%\begin{equation}
%P_{\mathrm{can}}=\left\{\left(\frac{2k}{n+1}-1,0\right)\in D : 1\leq k\leq n\right\}
%\end{equation}
Let $F\in\Homeo_+ (\Disk,P_{\mathrm{can}})$.
Denote by $[F]_{P_{\mathrm{can}}}$ the isotopy class of $F \ \mathrm{rel} \ P_{\mathrm{can}}$.
When it is clear from the context, we will also use the notation $[F]$.

The group of such isotopy classes under composition is called the {\it mapping class group} of $(\Disk,P_{\mathrm{can}})$ and is denoted by $\mathrm{MCG}_p$.

Let $F\in\Homeo_+(\Disk)$ possess a periodic orbit $P$ of smallest period $p$. 
We assume $P$ is in the interior of $\Delta$: if not, we extend $F$ arbitrarily over a collaring of $\Delta$.
Let $H\colon (\Disk,P_\mathrm{can})\to (\Disk,P)$ be any homeomorphism.
Then the {\it braid type} of $(P,F)$ is the conjugacy class $\<[H^{-1}\circ F\circ H]\>$ of $[H^{-1}\circ F\circ H]$ in $\mathrm{MCG}_p$. 
Denote the braid type of $(P,F)$ by $\mathrm{bt}(P,F)$.
Let $\mathrm{BT}_p$ denote the set of all braid types of a fixed period $p$.
\begin{rmk}
The braid type is independent of the collaring and the choice of homeomorphism $H$ (see~\cite{Hall94} for more details).
\end{rmk}
Let $F_0, F_1\in\Homeo_+(\Disk)$ possess periodic orbits $P_0$ and $P_1$ respectively. 
We say that the pair $(P_0,F_0)$ and $(P_1,F_1)$ are {\it braid equivalent} if $\mathrm{bt}(P_0,F_0)=\mathrm{bt}(P_1,F_1)$. 
Denote this equivalence by $(P_0,F_0)\sim_{BE}(P_1,F_1)$.
%(Clearly $\sim_{BE}$ is an equivalence relation.)
%this def used in dCH3
Equivalently, $(P_0,F_0)\sim_{BE}(P_1,F_1)$ if there exists a homeomorphism $H\colon(\Delta,P_0)\to (\Delta,P_1)$ such that $F_{0}\simeq H^{-1}\circ F_{1}\circ H \ \mathrm{rel} \ P_0$ in $\Delta$.

\subsubsection{Braids}
We now relate the notion of braid equivalence to that of a braid.
Let $B_p$ denote the braid group on $p$ strands (see~\cite{BirmanBook}).
Denote the composition of braids $\alpha, \beta\in B_p$ by $\alpha\cdot\beta$. 
If $\alpha$ and $\beta$ are {\it conjugate} we write $\alpha\sim\beta$.
If $\alpha$ and $\beta$ are {\it reverse-conjugate}, i.e. $\alpha\cdot\gamma=\gamma^{-1}\cdot\beta$ for some $\gamma\in B_p$, we write $\alpha\sim_r\beta$.
\begin{rmk}
%Let us make another disclaimer: 
We will consider braids, geometric braids and braid diagrams.
Typically we will not make the distinction. 
However, where necessary we will use the notation $a\simeq b$ to denote that the geometric braids $a$ and $b$ are isotopic.
We will also denote their product by $a\cdot b$, whenever the set of terminal endpoints of $b$ coincides with the set of initial endpoints of $a$. 
\end{rmk}

Let $Z(B_p)$ denote the centre of $B_p$.
It is known that $B_p/Z(B_p)$ is naturally isomorphic to $\mathrm{MCG}_p$ (see~\cite{VLHansen1}).
Hence $\mathrm{BT}_p$ is in one-to-one correspondence with the set of conjugacy classes of $B_p/Z(B_p)$ whose underlying permutation is a cycle (and hence close up to give a knot, rather than just a link).
Therefore, given $(P_F,F)$ we denote by $\beta(P_F,F)$ the conjugacy class in $B_p/Z(B_p)$ corresponding to $\mathrm{bt}(P_F,F)$.

\begin{rmk}\label{rmk:BE-equiv2}
%In what follows we will consider braid diagrams, rather than geometric braids.
%Therefore we wish to rephrase Remark~\ref{rmk:BE-equiv} in terms of braid diagrams.
%Let $a,b$ be geometric braids and let $\alpha, \beta$ be associated braid diagrams.
%Then $a\simeq b$ if and only if $\alpha$ can be transformed into $\beta$ by a series of isotopies, 2nd and 3rd Reidemeister moves and their inverses (see~\cite{KasselsTuraev1}). 

%Applying this, 
The following is well-known. 
Let $(P_0,F_0)$ and $(P_1,F_1)$ have associated braids $\beta_0$ and $\beta_1$. 
Then $(P_0,F_0)\sim_{BE} (P_1,F_1)$ if and only if there exists a braid $\sigma$ and $m\in \ZZ$ such that $\sigma^{-1}\cdot \beta_0\cdot\sigma$ can be deformed into $\beta_1\cdot\tau^m$ by a sequence of isotopies, 2nd and 3rd Reidemeister moves and their inverses, where $\tau$ is a braid associated to a generator of $Z(B_p)$.
\end{rmk}

Finally, we observe that if $F\in \overline{\Homeo_+(\Disk)}$ possesses a periodic orbit $P$ of smallest period $p$ which corresponds to an isolated fixed point of non-zero index for the iterate $F^p$, then any small perturbation $F'\in \Homeo_+(\Disk)$ will also possess a periodic orbit $P'$ which is a continuation of $P$ and of the same period. 
Therefore we can define $\mathrm{bt}(P,F)=\mathrm{bt}(P',F')$.
Since having an isolated fixed point of non-zero index is an open property, this is well-defined.

\subsubsection{Braid Equivalence in Families}
\begin{defn}[Braid equivalence in parametrised families]
Let $k\in \NN$.
Let $B\subset \RR^k$ be a contractible bounded open set.
Let $F\in C(B,\overline{\Homeo_+(\Disk)})$ be a $k$-parameter family of continuous self-maps.
We will use the notation $F_b$ for the map $F(b)$. 
Let $B^0=\{b\in B : F_b\in \Homeo_+(\Disk)\}$.
%Assume that $\lim{b\to b^0}F_b$ is well-defined for all $b^0\in B\setminus B^0$.
Let $b_0, b_1\in B$ satisfy the property that $F_{b_0}$ and $F_{b_1}$ have periodic orbits $P_0$ and $P_1$ respectively, both of smallest period $p$. 
Then $(P_0,F_{b_0})$ and $(P_1,F_{b_1})$ are {\it braid equivalent in the family $F$} if 
\begin{enumerate}
\item 
there exists a collection of $p$ pairwise distinct points 
\begin{equation}
P(t)=\{P^1(t),P^2(t),\ldots,P^p(t)\}\subset \Disk
\end{equation}
where $P^i(t)$ varies continuously with $t\in [0,1]$ for each $i=1,2,\ldots,p$, such that $P_0=P(0)$ and $P_1=P(1)$;
\item 
there exists a path $\gamma\colon [0,1]\to B$ such that $\gamma(0)=b_0, \gamma(1)=b_1$ and $P(t)$ is a periodic orbit for $F_{\gamma(t)}$, for each $t\in[0,1]$.
\end{enumerate} 
\end{defn}
(Observe, necessarily $P(t)$ must have smallest period $p$ for $F_{\gamma(t)}$ as the $P^i(t)$ are distinct.)
In other words, the one-parameter sub-family $F_{\gamma(t)}$ realises a strong Nielsen equivalence between $(P_0,F_{b_0})$ and $(P_1,F_{b_1})$.
\begin{rmk}
Braid equivalence in a parametrised family implies braid equivalence.
%The implication `braid equivalence in parametrised families implies braid equivalence' is as follows: 
%take $H_t=F_{t}\circ F_{0}^{-1}$ so that $H_t(P(0))=P(t)$. 
%Then $I_t=H_t^{-1}\circ F_t\circ H_t$ is an isotopy between $I_0=F_0$ and $I_1=H_1^{-1}\circ F_1\circ H_1$ preserving $P(0)$. 
%Hence $I_t$ is an isotopy between $F_0$ and a conjugate of $F_1$ which preserves $P(0)$. 
%Therefore $(P_0,F_0)$ and $(P_1,F_1)$ are braid equivalent.
\end{rmk}

%Assume there exists $S\subset\del B$ such that $s\in S$ implies $F_s=\lim_{b\to s}F_b$ exists and is well-defined, with respect to the uniform topology on $C(\Disk,\Disk)$ (or equivalently the compact-open topology).
%\begin{defn}[Braid equivalence in degenerate families]
%Let $s_0,s_1\in S$ satisfy the property that $F_{s_0}$ and $F_{s_1}$ possess periodic orbits $P_0$ and $P_1$ respectively, both of smallest period $p$.
%Then $(P_0,F_{s_0})$ and $(P_1,F_{s_1})$ are {\it braid equivalent in $F$} if 
%\begin{enumerate}
%\item 
%there exists a collection of $p$ pairwise distinct points 
%\begin{equation}
%P(t)=\{P^1(t),P^2(t),\ldots,P^p(t)\}\subset \Disk
%\end{equation}
%where $P^i(t)$ varies continuously with $t\in [0,1]$ for each $i=1,2,\ldots,p$, such that $P_0=P(0)$ and $P_1=P(1)$;
%\item 
%there exists a path $\gamma\colon[0,1]\to B\cup S$ such that $\gamma(0)=s_0$, $\gamma(1)=s_1$, and $P(t)$ is a periodic orbit of $F_{\gamma(t)}$, for all $t\in[0,1]$;
%\end{enumerate}
%\end{defn}
%\begin{rmk}
%We could encapsulate both of these definitions into a single, more general, definition by replacing homeomorphisms with continuous self-maps (and hence the notion of isotopy with homotopy). However, in what follows we will be interested mostly in the behaviour of H\'enon and H\'enon-like diffeomorphisms as they degenerate to unimodal maps.
%\end{rmk}
We will be specifically interested in the case when $b_0, b_1\in B\setminus B^0$ and $\gamma(b_0,b_1)\subset \Homeo_+(\Disk)$. 
(The motivating example will be that $F$ denotes the family of H\'enon maps and $F_{b_0}, F_{b_1}$ will correspond to quadratic maps.)

\subsubsection{Braids and Braid Diagrams}\label{sect:braids}
Recall that braids can be represented by braid diagrams.
Braid diagrams will be normalised in the following way: 
They lie in the unit square $[-1,1]\times [-1,1]$. 
Each strand has an {\it initial endpoint} lying in $[-1,1]^\init=[-1,1]\times \{1\}$.
Each strand has a {\it terminal endpoint} lying in $[-1,1]^\term=[-1,1]\times \{-1\}$.
The vertical line through any initial endpoint contains a terminal endpoint (and vice versa).
All crossings are transverse.
Only double points are allowed.
The strands are directed downwards.
(The last condition ensures no strand can `backtrack'.) 
Consequently, if $\alpha$ and $\beta$ are braids with associated diagrams, then the braid diagram corresponding to the product $\alpha\cdot\beta$ is the braid formed by placing the diagram corresponding to $\beta$ directly above the diagram corresponding to $\alpha$ and rescaling.
\begin{rmk}
We adopt this convention as it coincides with the convention of composing maps or isotopies from the right.
For example, if the isotopy $F^\alpha_t$ realises the braid $\alpha$ and the isotopy $F^\beta_t$ realises the braid $\beta$ then $F^\alpha_t\cdot F^\beta_t$ realises the braid $\alpha\cdot\beta$.
\end{rmk}
Let $\beta$ be a braid diagram.
Denote by $S_\beta$ the set of strands.
Denote by $E^{\init}_\beta$ the set of initial endpoints.
We will denote the points in $E^{\init}_\beta$, ordered from left to right, by $0^{\init},1^{\init},\ldots,(p-1)^{\init}$.
Similarly, denote by $E^{\term}_\beta$ the set of terminal endpoints.
We will denote the points in $E^{\term}_\beta$, ordered from left to right, by $0^{\term},1^{\term},\ldots,(p-1)^{\term}$.
We denote the strand emanating from $i^{\init}$ by $s_\beta(i)$.
Let $\pi\colon [-1,1]^2\to [-1,1]$ denote the projection to the first coordinate.
Then by assumption $\pi(i^\init)=\pi(i^\term)$ for all $i=0,1,\ldots,p-1$.

Given a strand $s\in S_\beta$ let $s^\init$ and $s^\term$ denote its initial and terminal endpoints respectively.
When it is clear which braid $\beta$ is being considered we will drop $\beta$ from our notation,  so that $s_\beta(i)$ becomes $s(i)$, $s_\beta(i)^\init$ becomes $s(i)^\init$, and so on.

%A braid $\beta$ induces a permutation $\sigma$ by assigning $s_\beta(i)^{\term}=\sigma(i)^{\term}$.

\subsection{Unimodal Dynamics}\label{subsect:unimodal}
\subsubsection{Unimodal Permutations}
\begin{rmk}
In this section we consider only intervals in $\{0,1,\ldots,p-1\}$, 
e.g. for $i,j\in \{0,1\ldots,p-1\}$ satisfying $i<j$ we let 
\begin{equation}
[i,j]=\{k\in \{0,1,\ldots,p-1\} : i\leq k\leq j\}
\end{equation}
Define $(i,j)$, $[i,j)$ and $(i,j]$ similarly. 
In later sections we may also assume they are embedded in $\RR$ but it will be clear from the context what is meant.
\end{rmk}
\begin{defn}
Let $p\in \NN$. 
Endow the set $\{0,1,\ldots,p-1\}$ with its natural ordering.
A permutation $\upsilon$ of the set $\{0,1,\ldots,p-1\}$ is {\it unimodal} if there exists $m\in (0,p-1)$, such that 
\begin{enumerate}
\item $\upsilon$ is order-preserving on the interval $[0,m]$
\item $\upsilon$ is order-reversing on the interval $[m,p-1]$
\end{enumerate}
We call $m$ the {\it folding point} of $\upsilon$.
%We call the interval $L=[0,m)$ the {\it left interval} and 
%we call the interval $R=(m,p-1]$ the {\it right interval} of $\upsilon$.
We denote the set of unimodal permutations on $\{0,1,\ldots,p-1\}$ by $U_p$.
Let $U=\bigcup_{p\in\NN}U_p$.
\end{defn}
If a unimodal permutation $\upsilon$ is cyclic then we can introduce the following, which we call {\it cyclic notation}:
Define $\o\colon \{0,1,\ldots,p-1\}\to\{0,1,\ldots,p-1\}$ by $\o(\upsilon^i(m))=i$ for $i=0,1,\ldots,p-1$.
Then we may uniquely represent $\upsilon$ by $(\o(0),\o(1),\o(2),\ldots,\o(p-1))$.
Observe that $\o$ is a bijection.
Hence we can recover $\upsilon$ by setting $\upsilon(\o^{-1}(i))=\o^{-1}(i+1\! \mod \ p)$.
\begin{eg}\label{eg:main}
The cyclic unimodal permutation $\upsilon$ of $\{0,1,2,3,4\}$ given by
\begin{equation}\label{perm:unimodal1}
\upsilon(0)=1, \ 
\upsilon(1)=3, \ 
\upsilon(2)=4, \ 
\upsilon(3)=2, \ 
\upsilon(4)=0
\end{equation}
has folding point $m=2$ and therefore $\o(2)=0$.
Applying $\upsilon$ iteratively gives 
\begin{equation}
\o(4)=1, \
\o(0)=2, \
\o(1)=3, \
\o(3)=4
\end{equation}
In cyclic notation
\footnote{Note that we will use underlinings to distinguish representations, so $\uline{3}$ denotes the third point in the orbit of $\uline{1}=p-1$, but $3$ denotes the third point from the left.} 
we write $(\o(0),\o(1),\o(2),\o(3),\o(4))=(\uline{2},\uline{3},\uline{0},\uline{4},\uline{1})$. 
(In other words, cyclic notation is just shorthand for the collection of inequalities $\upsilon^2(m)<\upsilon^3(m)<m<\upsilon^4(m)<\upsilon(m)$.)
\end{eg}

%\begin{defn}
%Let $\upsilon\in U_p$ be cyclic. 
%Fix $q\in \{0,\ldots,p-1\}$. 
%Then the {\it sequence of closest returns} to $q$ is the sequence of times $t_0<t_1<\ldots$ such that 
%\begin{itemize} 
%\item if $\perm^{t_i}(q)<q$ then $\perm^t(q)<q, t<t_i$ implies $\perm^{t}(q)<\perm^{t_i}(q)$;
%\item if $\perm^{t_i}(q)>q$ then $\perm^t(q)>q, t<t_i$ implies $\perm^{t}(q)>\perm^{t_i}(q)$;
%\end{itemize}
%The times $t_i$ above are called the {\it closest return times} for $q$.
%We distinguish these as a {\it left closest return} in the first case and a {\it right closest return} in the second case. 
%
%\end{defn}
%For example, the unimodal permutation $\upsilon$ given by~\eqref{perm:unimodal1} has return times to the folding point $m=2$ given by $t_1=1, t_2=2, t_3=3, t_4=4$. Of these, $t_2$ and $t_3$ are left return times and $t_1$ and $t_3$ are right return times.

\begin{defn}
Let $\upsilon\in U_p$ be cyclic.
Let $q\in\{0,1,\ldots,p-1\}$. 
We call $q$ a {\it right closest return time} to the folding point if $\upsilon^q(m)<m$ and the interval $(\upsilon^q(m),m)$ does not contain $\upsilon^r(m)$ for any integer $r$.
Similarly, call $q$ a {\it left closest return time} to the folding point if $m>\upsilon^q(m)$ and the interval $(m,\upsilon^q(m))$ does not contain $\upsilon^r(m)$ for any integer $r$.
Finally, we call $q\in\{0,1,\ldots,p-1\}$ a {\it closest return time} to the folding point if the interval $(\upsilon^{q+1}(m),\upsilon(m))$ does not contain any point $\upsilon^{r}(m)$ for any integer $r$.
\end{defn}
Note that a closest return time will either be a left or right closest return time.

\begin{defn}
Let $\upsilon\in U_p$.
Let $i, j, k\in \{0,1,\ldots,p-1\}$, $i<j$.
We say the closed interval $[i,j]$ {\it maps over $k$} if $k\in \upsilon[i,j]$ and {\it maps strictly over $k$} if $k\in \upsilon(i,j)$.
\end{defn}
\begin{rmk}
For each $k\in \{0,1,\ldots,p-1\}$, $k\neq 0$ or $p-1$, it is clear that there exists an interval mapping strictly over $k$.
It is also clear that each interval mapping strictly over the folding point $m$ contains at least one subinterval of shortest length which also maps strictly over $m$. 
(The interval $[\upsilon^{-1}(k)-1,\upsilon^{-1}(k)+1]$ maps strictly over $k$ and no strict subinterval also satisfies this property.)
In fact, there are at most two intervals mapping over $k$ of shortest length.
\end{rmk}

\begin{defn}
Let $\upsilon\in U_p$.
For each $k\in \{0,1,\ldots,p-1\}$ we call $\upsilon^{-1}(k)$ the {\it dynamical preimage} of $k$.
If $k\neq 0, p-1$, the {\it interval containing the dynamical preimage of $k$} is the shortest closed interval mapping strictly over $k$ which contains the dynamical preimage $\upsilon^{-1}(k)$ of $k$.
If $k=0$ or $p-1$, the {\it interval containing the dynamical preimage of $k$} is the shortest closed interval mapping over $k$ which contains $\upsilon^{-1}(k)$. 

The other shortest closed interval strictly mapping over $k$, when it exists, is called the {\it interval containing the non-dynamical preimage of $k$}.
\end{defn}

\begin{eg}\label{eg:missing-preimage}
The unimodal cyclic permutation $\upsilon$ of $\{0,1,2,3\}$ given by
\begin{equation}
\upsilon(0)=2, \ 
\upsilon(1)=3, \ 
\upsilon(2)=1, \ 
\upsilon(3)=0
\end{equation}
(which in cyclic notation is given by $(\uline{2},\uline{0},\uline{3},\uline{1})$) does not have an interval containing the non-dynamical preimage of the folding point $m=1$.
To see this, observe that $1$ has dynamical preimage $2$ lying in the right interval. 
The only non-empty subinterval of the left interval is the left interval itself, $[0,1]$. 
This maps to the interval $[2,3]$. 
As $[2,3]$ doesn't contain the folding point, the unimodal permutation $\upsilon$ does not have an interval containing the non-dynamical preimage of the folding point.
\end{eg}

\begin{defn}
Let $\upsilon\in U_p$ be cyclic.
Given $i, j\in \{0,1,\ldots,p-1\}$, denote by $\kappa=\kappa(i,j)>0$ the smallest positive integer such that $j=\upsilon^\kappa(i)$.
If $\upsilon^l(i)\neq m$ for all $0\leq l<\kappa$, let
\begin{equation}
\rho(i,j)=\mathrm{card} \{0\leq l<\kappa : \upsilon^l(i)>m\}
\end{equation}
\end{defn}

\begin{defn}
Let $\upsilon\in U_p$ be cyclic.
We say that $\upsilon$ is {\it reconnectable at the dynamical preimage} if the following property holds: 
let $D=[d^-,d^+]$ denote the interval containing the dynamical preimage $d$ of $\uline{0}$. 
Then either $\rho(\uline{1},d^-)$ is odd or $\rho(\uline{1},d^+)$ is even, or both.

We say that $\upsilon$ is {\it reconnectable at the non-dynamical preimage} if the following properties hold:
\begin{enumerate}
\item  {\it [Preimage condition]} 
The interval containing the non-dynamical preimage of $m$ exists. 
Denote it by $E=[e^-,e^+]$.
\item  {\it [Parity condition]} 
Either $\rho(\uline{1},e^-)$ is odd or $\rho(\uline{1},e^+)$ is even, or both.
%(i.e. the parity from $1$ to $i$ is odd or the parity from $1$ to $i+1$ is even, or both).
\end{enumerate}
\end{defn}
\begin{rmk}
It will become clear in what follows that this notion also makes sense for other points $k\in \{0,1,\ldots,p-1\}$, $k\neq m$.
However, for simplicity we will only consider reconnections at the folding point.
\end{rmk}
Recall that we aim to construct braid-equivalent pairs of unimodal combinatorial types starting from a given unimodal combinatorial type.
In what follows it will become clear that the construction we propose works precisely when the initial combinatorial type is reconnectable, either at the dynamical preimage or non-dynamical preimage. 

%Roughly speaking, the motivation for this definition is the following. Given a unimodal map with periodic folding point we first `disconnect' the folding orbit via a small perturbation at the dynamical preimage of the folding point. The orbit of the new point we then make to shadow the orbit of the folding point until we can `glue' an appropriate iterate of the new point to the critical point. Necessarily the iterate of the new point must lie in either the interval containing the dynamical preimage of the folding point, or the interval containing the non-dynamical preimage of the folding point. Therefore, we need to ensure intervals containing dynamical or non-dynamical preimages of the folding point exist (the preimage condition) and that when shadowing the folding orbit, we eventually land inside this interval (the parity condition).   

%\begin{rmk}
%We can similarly define the preimages of an arbitrary point in $\{0,1,\ldots,p-1\}$ with respect to a unimodal permutation. 
%In fact it is clear that multimodal permutations can be defined similarly and the dynamical and non-dynamical preimages can be defined in this case as well. However, for simplicity we only consider unimodal permutations and critical preimages as these are the only objects considered below.
%\end{rmk}

\begin{eg}\label{eg:non-reconnectable}
Given a cyclic unimodal permutation $\upsilon$, even if the interval containing the non-dynamical preimage of $m$ exists it may not be reconnectable. 
For example consider, in cyclic notation, $(\uline{2},\uline{7},\uline{3},\uline{8},\uline{0},\uline{5},\uline{4},\uline{9},\uline{6},\uline{1})$. 
The dynamical preimage $\uline{9}$ of the folding point $\uline{0}$ lies to the right of $\uline{0}$.
The interval containing the non-dynamical preimage of $\uline{0}$ is $[\uline{7},\uline{3}]$.
To calculate $\rho(\uline{1},\uline{3})$ we count the number of elements of $\{\uline{1},\uline{2}\}$ lying to the right of $\uline{0}$, of which there is $1$.
Similarly there are $4$ elements of $\{\uline{1},\uline{2},\uline{3},\uline{4},\uline{5},\uline{6}\}$ lying to the right of $\uline{0}$.
Hence $\rho(\uline{1},\uline{3})=1$ and $\rho(\uline{1},\uline{7})=4$. 
Since $\uline{7}<\uline{3}$ the parity condition isn't satisfied.
\end{eg}

\subsubsection{Unimodal Braids}\label{sect:unimodalbraids}
\begin{figure}
\begin{subfigure}[b]{0.5\linewidth}
\centering
\includegraphics[page=1,scale=0.1, bb=0 0 320 320]
{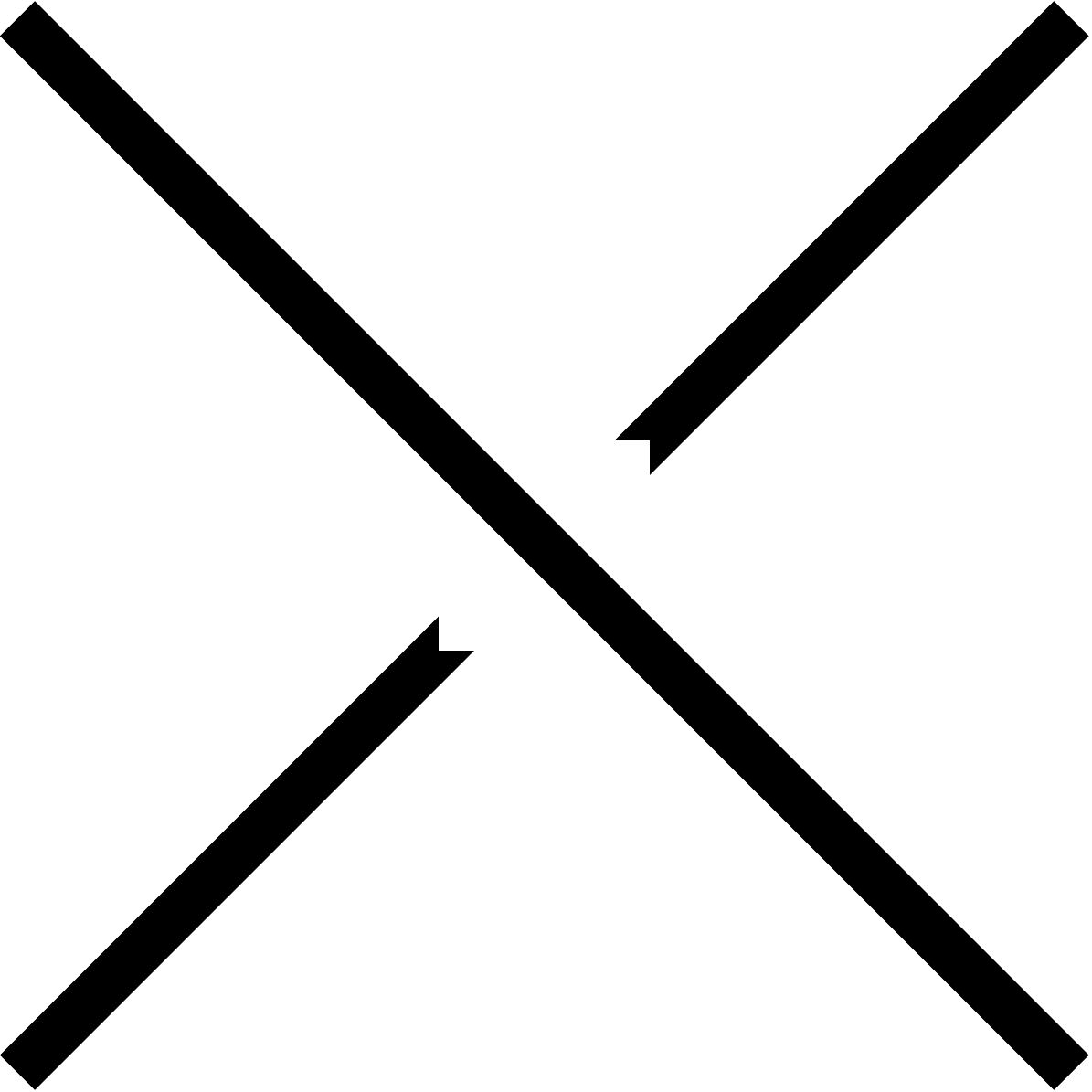}
%{crossing-positive.eps}
\vskip 5 pt
\caption{Positive crossing}
\label{fig:positivecrossing}
\end{subfigure}
\begin{subfigure}[b]{0.5\linewidth}
\centering
\includegraphics[page=2,scale=0.1]
{BEiH1-figs_epsi-2015-03-30-crop.pdf}
%{crossing-negative.eps}
\vskip 5 pt
\caption{Negative Crossing}
\label{fig:negativecrossing}
\end{subfigure}
\caption{The two types of crossings used to construct braid diagrams.}
\label{fig:crossings}
\end{figure}

We call a braid {\it positive} if each crossing is positive (see Figure~\ref{fig:crossings}).
A braid is {\it direct} if each strand crosses any other strand at most once\footnote{
In the literature such braids are called permutation braids -- however, it will be useful to have an adjective to describe this property, as in~\cite{Hall94}.}.
%We shall also use these adjectives for any collection of strands.
It is known that for any permutation $\upsilon\in U_p$ there is a unique positive, direct braid $\beta$ which induces $\upsilon$. 
%(In fact, to any permutation there corresponds a unique positive, direct braid which induces it.)
We call such braids {\it unimodal}.
Denote the set of unimodal braids on $p$ strands by $UB_p$.
Let $UB=\bigcup_{p\in\NN}UB_p$.

A unimodal braid $\beta$ possesses a canonical braid diagram satisfying:
\begin{enumerate}
\item\label{property:leftward} 
All strands initiated at $m$ or to the left of $m$ travel rightwardly.
\item\label{property:rightward} 
All strands initiated at a point to the right of $m$ travel rightwardly, touch the {\it folded line} $\pi^{-1}(p-1)$, 
then travel leftwardly.
\end{enumerate}
We call strands satisfying Properties~\ref{property:leftward} and~\ref{property:rightward} {\it unimodal strands}.
Henceforth we will identify a unimodal braid $\beta$ with the corresponding canonical unimodal braid diagram. 
%(In particular, we will use $\beta$ to denote both the braid and the braid diagram.)
For example, the canonical unimodal braid diagram for the permutation $(\uline{2},\uline{3},\uline{0}, \uline{4},\uline{1})$, 
is shown in Figure~\ref{fig:braid1-} below.

\vspace{10pt}

\begin{figure}[h]
\centering
\psfrag{0in}{$\uline{0}^{\init}$}
\psfrag{1in}{$\uline{1}^{\init}$}
\psfrag{2in}{$\uline{2}^{\init}$}
\psfrag{3in}{$\uline{3}^{\init}$}
\psfrag{4in}{$\uline{4}^{\init}$}
\psfrag{0te}{$\uline{0}^{\term}$}
\psfrag{1te}{$\uline{1}^{\term}$}
\psfrag{2te}{$\uline{2}^{\term}$}
\psfrag{3te}{$\uline{3}^{\term}$}
\psfrag{4te}{$\uline{4}^{\term}$}
\includegraphics[page=3,scale=0.5]{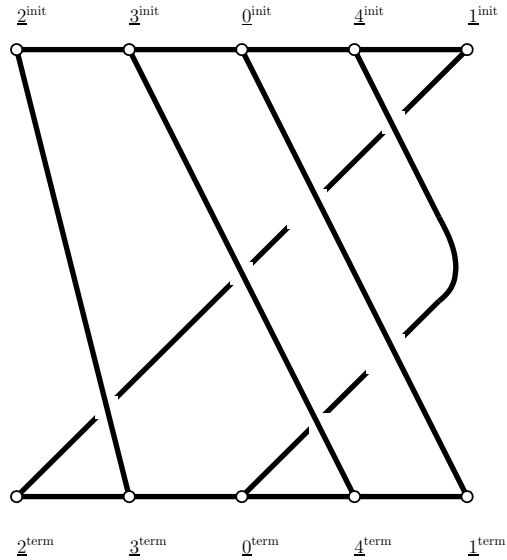}
%{2-3-0-4-1-braid.eps}
\caption{The braid $\beta$ for the cyclic unimodal permutation $(\uline{2},\uline{3},\uline{0},\uline{4},\uline{1})$}
\label{fig:braid1-}
\end{figure}
Let $S$ be a collection of unimodal strands, not necessarily forming a braid, which are positive, direct and which contains a folding strand $s(\uline{0})$. 
Let $s$ be a strand in $S$.
The strand $r$ {\it unimodally follows} $s$ in $S$ if
\begin{enumerate}
\item $r$ is unimodal, $r^\init$ neighbours $s^\init$, and $r^\term$ neighbours $s^\term$
\item the collection $S\cup\{r\}$ is positive and direct
\item $r$ has crossings of the following types
\begin{itemize}
\item if $t$ is a strand crossing $s$: then it also crosses $r$ and the crossing is of the same type. All such crossings occur in the same order for $r$ as they do for $s$. 
\item if $s$ hits the folded line: then $r$ hits the folded line and they cross once, when $s$ is moving rightwardly but $r$ is still moving leftwardly, or vice versa. 
Otherwise there is no crossing. 
\end{itemize}
\end{enumerate} 
%Observe that once a side of $s^\init$ is chosen there is, up to isotopy, a unique unimodal strand $r$ following $s$.
%IS THE FOLLOWING NEC.?
Let $s$ and $s'$ be strands in $S$ such that $s^\init< s'^\init$ are neighbours lying on the same side of $\uline{0}$.
A strand $r$ {\it unimodally follows} the pair $s$, $s'$ in $S$ if 
\begin{enumerate}
\item $r$ is unimodal, $r^\init$ lies between $s^\init$ and $s'^\init$, and $r^\term$ lies between $s^\term$ and $s'^\term$
\item the collection $S\cup\{r\}$ is positive and direct. 
%--alternate version
%\item $r$ lies between $s$ and $s'$ whenever $s$ and $s'$ both run either leftwardly or rightwardly;
\item $r$ has crossings of three types
\begin{itemize}
\item if $t$ is a strand crossing both $s$ and $s'$: then $t$ makes the same type of crossing with $r$ between its crossings with $s'$ and $s$;
%--alternate version
%\item if $t$ is a strand also between $s$ and $s'$ then $r$ and $t$ make a positive crossing at the folded line;
\item if $t$ is a strand crossing $s'$ but not $s$: if $r^\term$ lies between $s'^\term$ and $t^\term$, then $r$ makes the same type of crossing with $t$ after $s'$ crosses $t$. Otherwise there is no crossing.
Similarly, if $t$ crosses $s$ but not $s'$.
\item if $s$ and $s'$ hit the folded line, then $r$ makes a single crossing both with $s$ and with $s'$, either side of the single crossing between $s$ and $s'$. 
Otherwise there are no crossings.
\end{itemize}
\end{enumerate}
%Observe that if $r^\init$ lies between $s^\init$ and $s'^\init$ and $r^\term$ lies between $s^\term$ and $s'^\term$ then up to isotopy there is a unique unimodal strand $r$ with initial and terminal endpoints $r^\init$ and $r^\term$ respectively.

Before proceeding, let us make the following trivial observation concerning the action of half-twists on a cabled pair of strands.
\begin{obs}[Fundamental Observation]\label{obs:cabling1}
Let $\beta$ be a braid.
Let $s$ be a strand of $\beta$.
Form a new strand $s'$ which follows the strand $s$ (i.e. makes the same crossings with all other strands and in the same order).
Allow an arbitrary number of crossings between $s$ and $s'$.
Let $\beta'$ denote the resulting braid\footnote{Formally, these are not braids as they do not have the same set of initial and terminal endpoints. 
To avoid confusion, let us call these objects {\it almost-braids}.} 
.
Let $\tau_\init$ denote a positive half-twist between $s^\init$ and $s'^\init$.
Let $\tau_\term$ denote a positive half-twist between $s^\term$ and $s'^\term$.
Then $\beta'\cdot\tau_\init=\tau_\term\cdot\beta'$. 
\end{obs}

\section{Breaking Braids Via First Closest Returns}\label{sect:gen-cabling}
\subsection{Cabling}
Our first construction extends that by Holmes~\cite{Holmes86} which generalised the cabling\footnote{Although Holmes only used the term cabling for iterated torus knots (not iterated horseshoe knots) we will use the term for both cases.
} construction for iterated torus knots to horseshoe braids.
Our description is in terms of braids rather than templates, but they are equivalent.
Note that Holmes did not consider braid equivalence of the pairs of braids created by this construction. 
He considered only the single unimodal braid that resulted.

We now give an informal description of the cabling procedure.
Given a cyclic braid $\beta\in B_p$, take the folding strand $s_0$.
`Break' the braid by disconnecting $s_0$ from $s_0^\term$ and `glue' it to a neighbouring point $r_0^\term$ which is not already the endpoint of some strand.
Let $r_1^\init$ and $s_1^\init$ denote the vertical projection of $r_0^\term$ and $s_0^\term$ respectively.
Let $s_1$ denote the strand emanating from $s_1^\init$.

Add a strand $r_1$ with initial point $r_1^\init$ and which unimodally follows $s_1$. 
Then $r_1$ has a terminal point $r_1^\term$ neighbouring $s_1^\term$.
Now repeat the process:
Let $r_2^\init$ and $s_2^\init$ denote the vertical projection of $r_1^\term$ and $s_1^\term$ respectively.
Add a strand $r_2$ unimodally following $s_2$, etc.

The final strand $r_{p-1}$ has initial point $r_{p-1}^\init$ and terminal point $r_0^\term$. 
However, it may not be possible for $r_{p-1}$ to follows $s_{p-1}$ unimodally, as $r_0^\term$ may lie on the wrong side of $s_0^\term$. 
Therefore $r_{p-1}$ unimodally follows $s_{p-1}$ until all other crossing in have been made, then makes a negative crossing between $s_{p-1}$ and $r_{p-1}$ before connecting to $r_0^\term$. 

Given a unimodal braid $\beta$ we can cable it in two distinct ways: to the left or right of the folding strand.
It can be shown that this pair of braids are conjugate or reverse-conjugate.
However, they are not both unimodal.

Observe that the cabling can be closed-up to form a braid whenever we land inside the interval containing the dynamical preimage of the folding point. 
More precisely, choose $q\in \NN$ so that $r_q^\term$ lies in an interval $(s_q^\term,s_{q'}^\term)$, not containing any terminal endpoints of $\beta$, and such that it maps over $s_0^\term$. 
Let $r_{q+1}^\init$ denote the vertical projection of $r_{q}^\term$.
At this last step, construct a strand $r_q$ emanating from $r_q^\term$ which, rather than completely following $s_q$, only follows $s_q$ until we can reconnect $r_q$ to $s_0^\term$ without creating any further crossings.
We call this a {\it generalised cabling at the dynamical preimage}.
Again there are two distinct cablings: to the left or right of the folding strand.
As in the cabling case, the pair of braids are conjugate or reverse-conjugate.
However, they are never both unimodal. 

In Construction~\ref{constr:gencabling2} below, we observe that the cabling can be closed-up at any point in time at which we lie in the interval containing the non-dynamical preimage of the folding point.
We call this {\it generalised cabling at the non-dynamical preimage}.
There are two preferred cablings, either side of the folding strand.
Theorem~\ref{thm:constr3} says this pair is conjugate or reverse-conjugate.
What is important is that, unlike the previous two constructions, the braids produced by this process are both unimodal.

\begin{landscape}
\begin{figure}[p]
\centering
\tiny
\psfrag{0in}{$\uline{0}^\init$}
\psfrag{1in}{$\uline{1}^\init$}
\psfrag{2in}{$\uline{2}^\init$}
\psfrag{3in}{$\uline{3}^\init$}
\psfrag{4in}{$\uline{4}^\init$}
\psfrag{0te}{$\uline{0}^\term$}
\psfrag{1te}{$\uline{1}^\term$}
\psfrag{2te}{$\uline{2}^\term$}
\psfrag{3te}{$\uline{3}^\term$}
\psfrag{4te}{$\uline{4}^\term$}
\begin{subfigure}[b]{0.4\textwidth}
\centering
\includegraphics[page=8,width=\textwidth]
{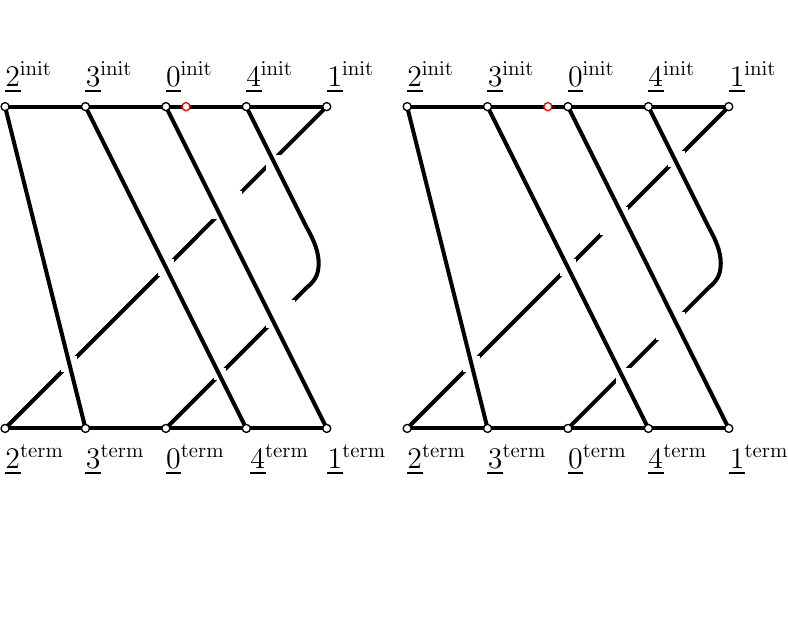}
%{2-3-0-4-1-braid6.eps}
\caption{Adding the first initial point $r_\pm^\init$ to $\beta$.}
\label{fig:2-3-0-4-1(a)}
\end{subfigure}
~\qquad
\begin{subfigure}[b]{0.4\textwidth}
\centering
\includegraphics[page=9,width=\textwidth]
{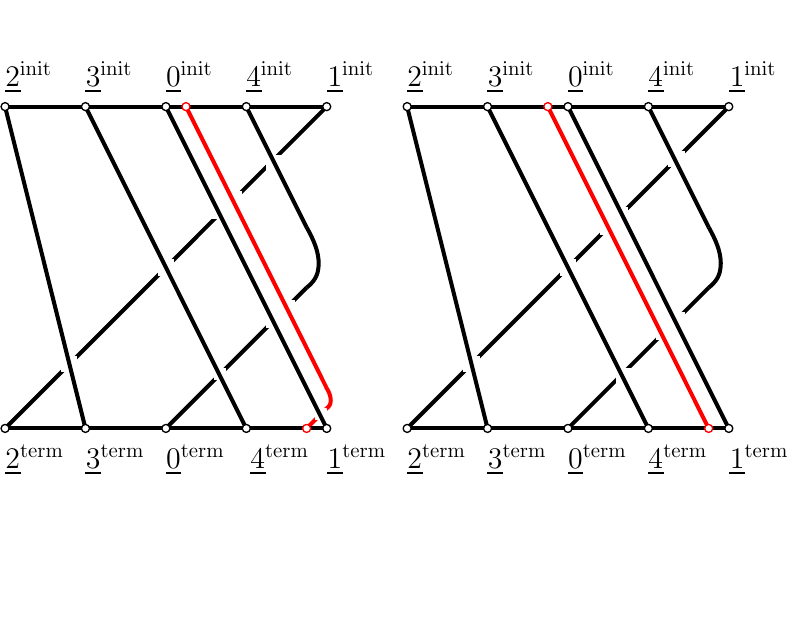}
%{2-3-0-4-1-braid5.eps}
\caption{Attaching the first strands $r_\pm(\uline{0})$.}
\label{fig:2-3-0-4-1(b)}
\end{subfigure}
~\qquad
\begin{subfigure}[b]{0.4\textwidth}
\centering
\includegraphics[page=10,width=\textwidth]
{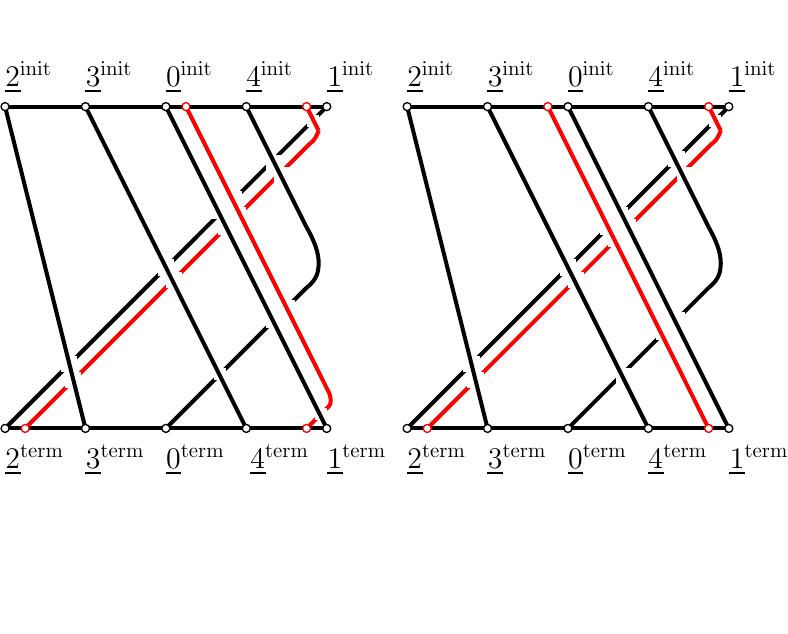}
%{2-3-0-4-1-braid4.eps}
\caption{Attaching the second strands $r_\pm(\uline{1})$.}
\label{fig:2-3-0-4-1(e)}
\end{subfigure}

\begin{subfigure}[b]{0.4\textwidth}
\centering
\includegraphics[page=11,width=\textwidth]
{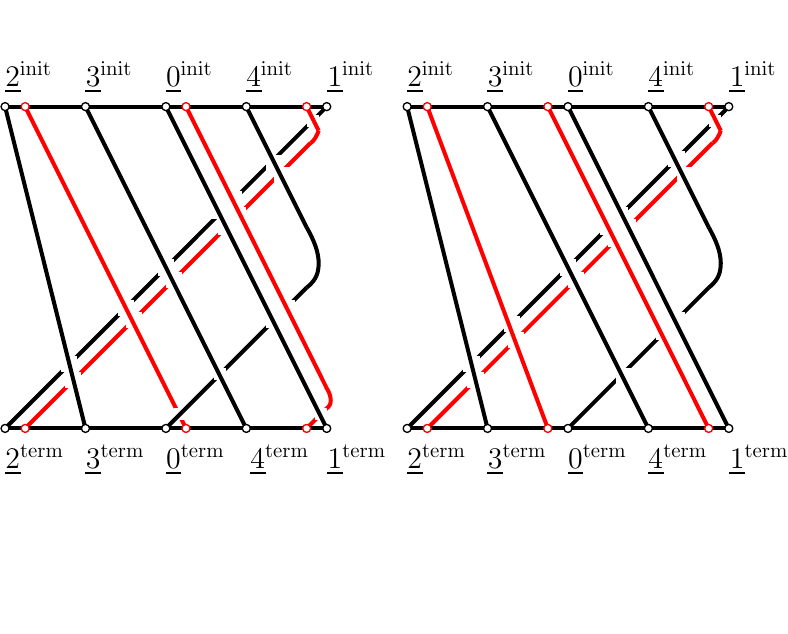}
%{2-3-0-4-1-braid3.eps}
\caption{Attaching the final strand $r_\pm(\uline{2})$.}
\label{fig:2-3-0-4-1(c)}
\end{subfigure}
~\qquad
\begin{subfigure}[b]{0.4\textwidth}
\centering
\includegraphics[page=12,width=\textwidth]
{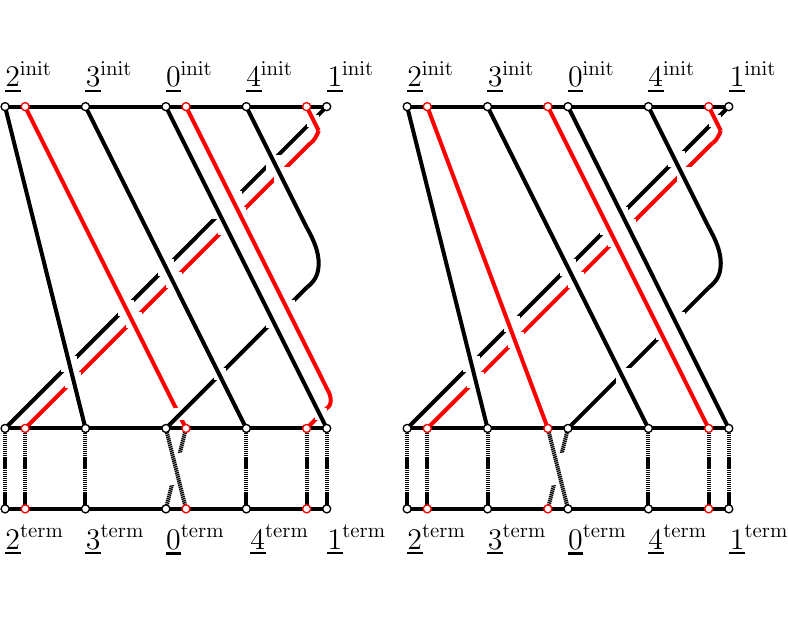}
%{2-3-0-4-1-braid2.eps}
\caption{Performing the positive half-twist $\tau$.}
\label{fig:2-3-0-4-1(d)}
\end{subfigure}
%~\qquad
%\begin{subfigure}[b]{0.4\textheight}
%\centering
%\includegraphics[width=\textwidth]{2-3-0-4-1-braid1.eps}
%\caption{....}
%\label{fig:2-3-0-4-1(f)}
%\end{subfigure}
\caption{Adding strands to two copies of $\beta$ from Example~\ref{eg:main-constr3} resulting in the new braids $\beta_-$ and $\beta_+$.}
\label{fig:eg-2-3-0-4-1}
\end{figure}
\end{landscape}

%HERE TWIST TAU IS POSITIVE I.E. WE USE NEGATIVE CROSSINGS
\begin{constr}[Generalised Cabling -- At the Non-Dynamical Preimage]\label{constr:gencabling2}
Let $\beta\in UB_p$ be cyclic and reconnectable at the non-dynamical preimage.
Let $\upsilon\in  U_p$ denote the corresponding unimodal permutation.
Let $E=[\uline{e}^-, \uline{e}^+]$ denote the interval containing the non-dynamical preimage of the folding point.
Then either $\rho(\uline{1},\uline{e}^-)$ is odd or $\rho(\uline{1},\uline{e}^+)$ is even, or both.
Let $\uline{q-1}$ denote either of the points $\uline{e}^-$ or $\uline{e}^+$ satisfying this parity condition.

Choose $\epsilon>0$ small.
Let $r_-\in (\uline{0}-\epsilon,\uline{0})$ and $r_+\in (\uline{0},\uline{0}+\epsilon)$.
%Call these respectively the {\it left} and {\it right new nearest neighbours} to the folding point.
Let $r_1\in (\uline{1}-\epsilon,\uline{1})$, $r_2\in (\uline{2},\uline{2}+\epsilon)$ and choose points $r_i\in (\uline{i}-\epsilon,\uline{i}+\epsilon)$, $i=3,\ldots,q-1$, which satisfy
\begin{equation}
r_i\in \left\{
\begin{array}{ll}
(\uline{i}-\epsilon,\uline{i}) & \mbox{if} \ \rho(\uline{1},\uline{i}) \ \mbox{is even} \\
(\uline{i},\uline{i}+\epsilon) & \mbox{if} \ \rho(\uline{1},\uline{i}) \ \mbox{is odd}
\end{array}
\right.
\end{equation}
For $i=\pm,1,2\ldots,q-1$, let $r_i^\init$ and $r_i^\term$ denote the projection of $r_i$ to $[-1,1]^\init$ and $[-1,1]^\term$ respectively. 
These will be endpoints for the strands constructed below.
 
First, let $r_\pm(\uline{0})$ denote the strand with initial endpoint $r_\pm^\init$ and terminal endpoint $r_1^\term$, which follows unimodally the strand $s(\uline{0})$.

Secondly, assume $r_\pm(\uline{i})$ have been constructed for $i=0,\ldots,j-1$, where $0<j\leq q-2$.
Let $r_+(\uline{j})=r_-(\uline{j})$ denote the strand with initial endpoint $r_i^\init$ and terminal endpoint $r_{i+1}^\term$ which follows unimodally the strand $s(\uline{i})$.

Thirdly, assume that the strands $r_\pm(\uline{i})$ have been constructed for all $i=0,1,\ldots,q-2$.
Construct the strand $r_\pm(\uline{q-1})$ as follows.
As $[\uline{e}^-,\uline{e}^+]$ maps over the folding point, one of the strands $s(\uline{e}^-)$ or $s(\uline{e}^+)$ has a terminal endpoint which lies on the same side of the folding point as its initial endpoint.
Call it $s$ and the other strand $s'$.
Let $r_\pm(\uline{q-1})$ denote the strand with initial endpoint $r^\init_{q-1}$ and terminal endpoint $r_\pm^\term$ which unimodally follows $s$ and $s'$.
We make modifications to $r_{\pm}(\uline{q-1})$ as follows:
\begin{enumerate}
\item[(I)] 
{\it [If $\uline{0}<\uline{p-1}$.]} %i.e. E<D
Since $r_+^\term>s(\uline{p-1})^\term$ and $r^\init_{q-1}<\uline{0}^\init<s(\uline{p-1})^\init$ we need an additional crossing between $r_+(\uline{q-1})$ and $s(\uline{p-1})$. 
Make a single {\it negative} crossing between $r_+(\uline{q-1})$ and $s(\uline{p-1})$ after all other crossings, resulting from $r_+(\uline{q-1})$ following $s$ and $s'$, have been made.

Since $r_-^\term<s(\uline{p-1})^\term$ and $r^\init_{q-1}<\uline{0}^\init<s(\uline{p-1})^\init$ no additional crossings between $r_-(\uline{q-1})$ and $s(\uline{p-1})$ are necessary.

\item[(II)] 
{\it [If $\uline{p-1}<\uline{0}$.]} %i.e. D<E
Since $r_-^\term<s(\uline{p-1})^\term$ and $s(\uline{p-1})^\init<\uline{0}^\init<r^\init_{q-1}$ we need an additional crossing between $r_+(\uline{q-1})$ and $s(\uline{p-1})$. 
Make a single {\it negative} crossing between $r_+(\uline{q-1})$ and $s(\uline{p-1})$ after all other crossings, resulting from $r_+(\uline{q-1})$ following $s$, have been made.

Since $r_+^\term>s(\uline{p-1})^\term$ and $s(\uline{p-1})^\init<\uline{0}^\init<r^\init_{q-1}$ no additional crossings between $r_+(\uline{q-1})$ and any other strand are necessary.
\end{enumerate}
Let $\alpha_\pm$ denote the braid consisting of the strands of $\beta$, together with the strands $r_\pm(\uline{j}),j=0,1,\ldots,q-1$ formed above. 
Note that $\alpha_\pm$ is not cyclic.
In fact, it is not necessarily unimodal. 
%\begin{enumerate}
%\item[(I)] $\alpha_-$ is unimodal. $\alpha_+$ is direct but not positive or unimodal.
%\item[(II)] $\alpha_-$ is direct, but not positive or unimodal. $\alpha_+$ is unimodal.
%\end{enumerate}
Let $\tau_\pm$ denote the positive half-twist between $r_\pm$ and $\uline{0}$.
(That is, $\tau_\pm$ is the Artin generator exchanging $\uline{0}$ and $r_\pm$ via a single positive crossing.)
Let $\beta_\pm=\tau_\pm\cdot\alpha_\pm$.
Note that $\beta_\pm$ is a cyclic braid.
Moreover, after cancelling appropriate positive and negative crossings between the strands $r_{\pm}(\uline{q-1})$ and $s(\uline{p-1})$ via Reidemeister moves and isotopy, we see that, in both cases (I) and (II), that $\beta_-$ and $\beta_+$ are unimodal.
\end{constr}
\begin{eg}\label{eg:main-constr3}
Consider the cyclic unimodal permutation $\upsilon$ from Example~\ref{eg:main}.
In cyclic notation it is given by $(\uline{2},\uline{3},\uline{0},\uline{4},\uline{1})$.
Then the interval containing the non-dynamical preimage is $E=[\uline{2},\uline{3}]$.
Since $\rho(\uline{1},\uline{2})=1$ is odd $\upsilon$ is reconnectable at the nondynamical preimage. 
Let $\beta$ denote the corresponding unimodal braid and $\beta_\pm$ the pair of braids from Construction~\ref{constr:gencabling2}.
If $\upsilon_\pm$ denotes the permutation corresponding to $\beta_\pm$ then
\begin{equation}
\upsilon_-=(\uline{2},\uline{7},\uline{3},\uline{5},\uline{0},\uline{4},\uline{6},\uline{1})
\qquad
\upsilon_+=(\uline{2},\uline{7},\uline{3},\uline{0},\uline{5},\uline{4},\uline{6},\uline{1})
\end{equation}
See Figure~\ref{fig:eg-2-3-0-4-1} for the construction of $\beta_\pm$ in this case.
\end{eg}

\begin{figure}
\centering
\tiny
\psfrag{sp-1in}{$s_{p-1}^\init$}
\psfrag{rq-1in}{$r_{q-1}^\init$}
\psfrag{sd-in}{$s_{\uline{d}^-}^\init$}
\psfrag{sd+in}{$s_{\uline{d}^+}^\init$}
\psfrag{s0in}{$s_{0}^\init$}
\psfrag{r-in}{$r_{-}^\init$}
\psfrag{r+in}{$r_{+}^\init$}
\psfrag{sigma}{$\sigma$}
\psfrag{alpha-}{$\alpha_-$}
\psfrag{tau-}{$\tau_-$}
\psfrag{alpha+}{$\alpha_+$}
\psfrag{tau+}{$\tau_+$}
\psfrag{sigma-1}{$\sigma^{-1}$}
\psfrag{r-te}{$r_{-}^\term$}
\psfrag{s0te}{$s_{0}^\term$}
\psfrag{r+te}{$r_{+}^\term$}
\psfrag{sd-+1te}{$s_{\uline{d^{-}+1}}^\term$}
\psfrag{sd++1te}{$s_{\uline{d^{+}+1}}^\term$}
\psfrag{r1te}{$r_{1}^\term$}
\psfrag{s1te}{$s_{1}^\term$}
\begin{subfigure}[b]{0.45\textwidth}
\centering
\includegraphics[page=4,width=\textwidth]
{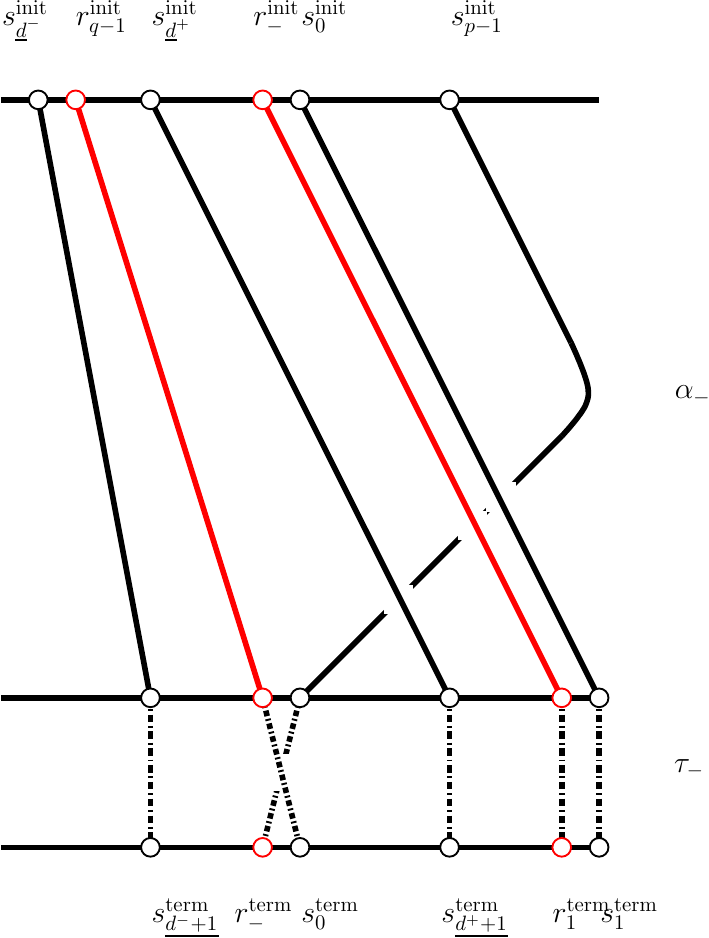}
%{braidsv2-cabl3I-a.eps}
\caption{$\beta_-$ in Construction~\ref{constr:gencabling2}(I).}
\label{fig:3I-a}
\end{subfigure}
~\qquad
\begin{subfigure}[b]{0.45\textwidth}
\centering
\includegraphics[page=5,width=\textwidth]
{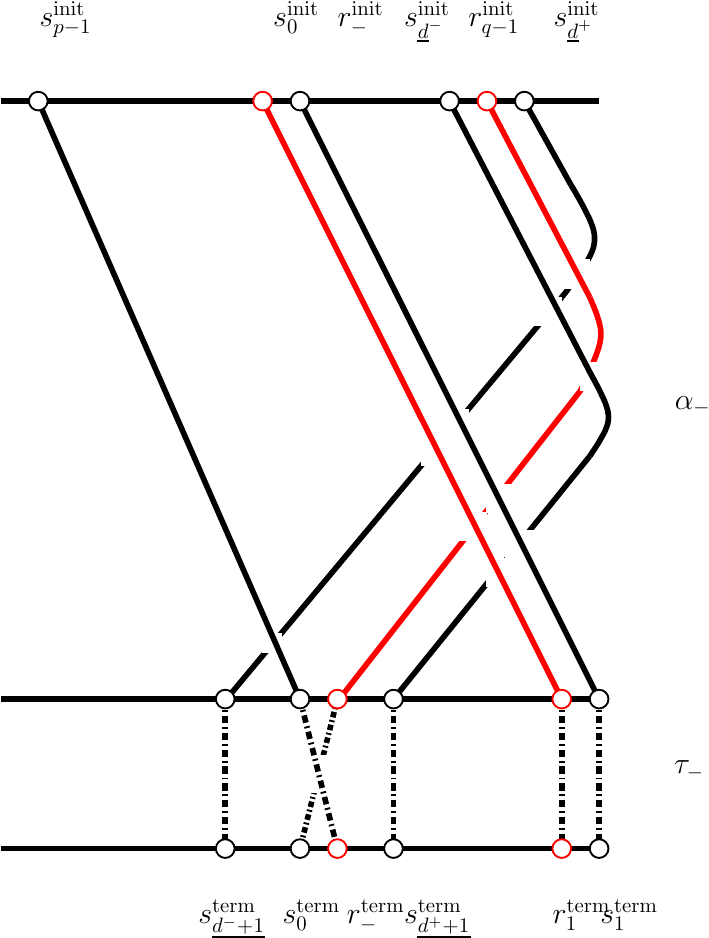}
%{braidsv2-cabl3II-a.eps}
\caption{$\beta_-$ in Construction~\ref{constr:gencabling2}(II).}
\label{fig:3II-a}
\end{subfigure}
~
\vspace{10pt}

\begin{subfigure}[b]{0.45\textwidth}
\centering
\includegraphics[page=6,width=\textwidth]
{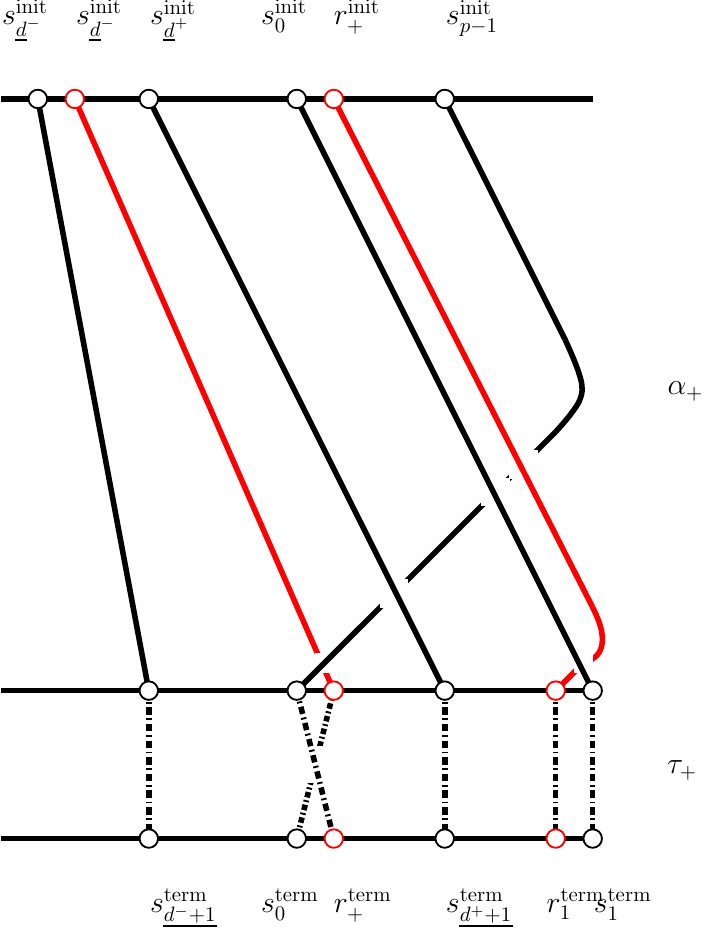}
%{braidsv2-cabl3I-b.eps}
\caption{$\beta_+$ in Construction~\ref{constr:gencabling2}(I).}
\label{fig:3I-b}
\end{subfigure}
~\qquad
\begin{subfigure}[b]{0.45\textwidth}
\centering
\includegraphics[page=7,width=\textwidth]
{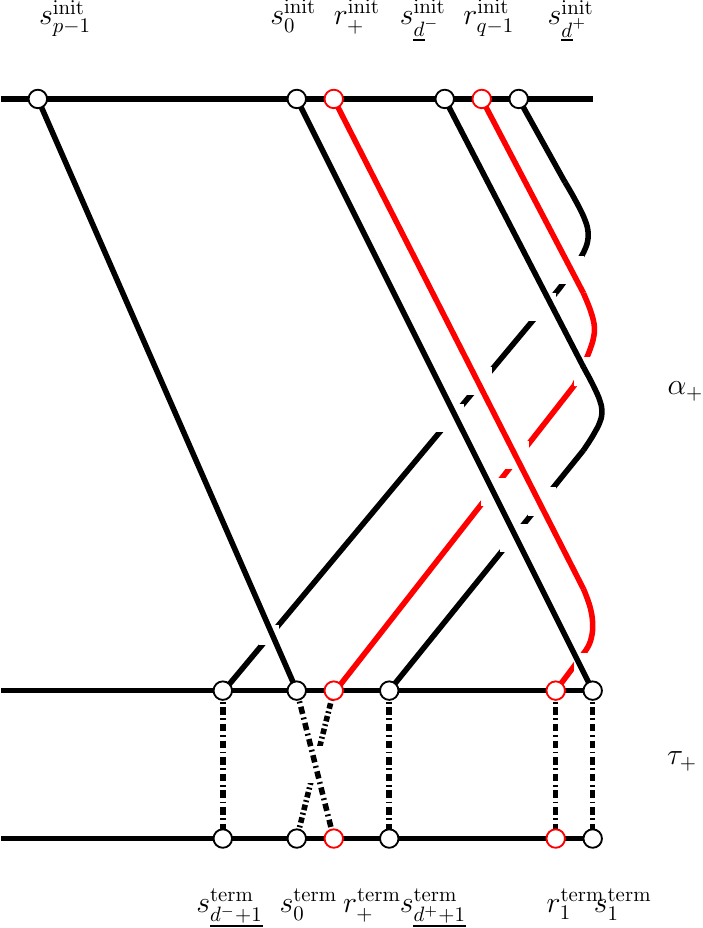}
%{braidsv2-cabl3II-b.eps}
\caption{$\beta_+$ in Construction~\ref{constr:gencabling2}(II).}
\label{fig:3II-b}
\end{subfigure}
\label{fig:3I}
\caption{The braids $\beta_-$ and $\beta_+$ from Construction~\ref{constr:gencabling2}. 
Only the necessary strands are shown.
Strands from the original braid $\beta$ are black, while new strands are red. 
(Only the first and last new strands are depicted.)}
\end{figure}

More generally, consider the case when the dynamical preimage of the folding point lies to the right of $m$.
Then $\beta_-$ and $\beta_+$ are depicted in Figures~\ref{fig:3I-a} and~\ref{fig:3I-b} respectively.
Identifying the sets of initial and terminal endpoints of $\beta_-$ and $\beta_+$ in an order-preserving manner, 
then precomposing $\beta_-$ with a positive half-twist between the strand $s(\uline{0})$ and its neighbour $r_{-}$, 
and post-composing by the inverse of this twist, yields the braid $\beta_+$.
A similar argument can be given when the dynamical preimage of the folding point lies to the left of $m$.
Then $\beta_-$ and $\beta_+$ are shown in Figures~\ref{fig:3II-a} and~\ref{fig:3II-b}.
However, in this case we precompose and post-compose by the same positive half-twist.
Hence we have the following Theorem.
%%%%%%%%%%%%%%%%%%%%%%%%%%%%%%%%%%%%%%%%%%%%%%%%%%%%%%%%%%%%%%%%%%%%%%%%%%%%%%%%%
%\begin{landscape}
%\end{landscape}
%%%%%%%%%%%%%%%%%%%%%%%%%%%%%%%%%%%%%%%%%%%%%%%%%%%%%%%%%%%%%%%%%%%%%%%%%%%
\begin{thm}\label{thm:constr3}
Let $p\in\NN$.
Let $\beta\in UB_p$ be cyclic.
Assume $\beta$ is reconnectable at the non-dynamical preimage.
Let $\beta_-$ and $\beta_+$ denote the braids from Construction~\ref{constr:gencabling2}
\begin{enumerate}
\item
%(conjugate)
If the non-dynamical preimage lies to the left of $m$ then $\beta_+\sim\beta_-$ 
\item 
%(reverse-conjugate)
If the non-dynamical preimage lies to the right of $m$ then $\beta_+\sim_r\beta_-$ 
\end{enumerate}
\end{thm}
Theorem~\ref{thm:constr3} and Remark~\ref{rmk:BE-equiv2} imply the following important corollary.
\begin{cor}
Let $p\in \NN$.
Let $\upsilon\in U_p$ be reconnectable at the non-dynamical preimage.
Let $\upsilon_-$ and $\upsilon_+$ denote the unimodal permutations from Construction~\ref{constr:gencabling2}.

If $f_-, f_+\in C([-1,1],[-1,1])$ are unimodal maps with periodic critical orbits $C_-$ and $C_+$ of type $\upsilon_-$ and $\upsilon_+$ respectively, then $(C_-,f_-)\sim_{BE}(C_+,f_+)$.
\end{cor}

\section{Breaking Braids Via Second Closest Returns}
\subsection{Generalised Cabling Process}
We now describe a second construction that, given a braid-equivalent pair of unimodal permutations $\upsilon_-$ and $\upsilon_+$ satisfying conditions given below, generates another pair of braid-equivalent unimodal permutations. 
Moreover, the new pair satisfy the same conditions and hence we can apply the construction once more.
The idea is the following.
Take $\upsilon_\pm$ as given by Construction~\ref{constr:gencabling2}. 
Break the connection at $\uline{0}^\term$ and glue the free strand to a point just outside the interval $[\uline{p}_-^\term,\uline{p}_+^\term]$, making a {\it second closest return}. 
Then consecutively add strands, starting from this point, that unimodally follow the strands from initial endpoints $\uline{p}_\pm^\init$, $\uline{p+1}_\pm^\init$, etc. 
Stop once a terminal endpoint lands inside $D$, the interval containing the dynamical preimage of the folding point. 
Finally, to close-up the braids, glue the last strand back to $\uline{0}^\term$.  
However, since we will define the process inductively we need to enlarge the space of such pairs (not just those coming from Construction~\ref{constr:gencabling2}).
Before giving the general construction, however, let us consider the following example.
\begin{eg}\label{eg:main-constr4}
Let $\beta_\pm$ denote the braids from Example~\ref{eg:main-constr3}.
For notational clarity, let us initially denote the underlying permutations by
\begin{equation}
\upsilon_-=(\uline{2}_-,\uline{7}_-,\uline{3}_-,\uline{5}_-,\uline{0}_-,\uline{4}_-,\uline{6}_-,\uline{1}_-)
\end{equation}
and 
\begin{equation}
\upsilon_+=(\uline{2}_+,\uline{7}_+,\uline{3}_+,\uline{0}_+,\uline{5}_+,\uline{4}_+,\uline{6}_+,\uline{1}_+)
\end{equation}
The braids $\beta_\pm$ are shown in Figure~\ref{fig:2-3-0-4-1-c2(a)}.
The initial and terminal endpoints of $\beta_\pm$ are denoted by $\uline{0}_\pm^\init,\uline{1}_\pm^\init,\ldots,\uline{7}_\pm^\init$ and $\uline{0}_\pm^\term,\uline{1}_\pm^\term,\ldots,\uline{7}_\pm^\term$ respectively.
Denote the corresponding strands by $s(\uline{0}_\pm),\ldots,s(\uline{7}_\pm)$.
We may assume that the endpoints of the strands in $\beta_\pm$ have been moved, in an order-preserving manner, so that for all $k\neq 5$, $\uline{k}_-^\init=\uline{k}_+^\init$ and $\uline{k}_-^\term=\uline{k}_+^\term$.
We may also assume that the strands have been deformed so that $s(\uline{k}_-)=s(\uline{k}_+)$ for all $k\neq 4$ or $5$ (i.e. any strand whose set of initial or terminal endpoints does not contain $\uline{5}_\pm^\init$ or $\uline{5}_\pm^\term$).
Therefore, for all $k\neq 5$, denote the points $\uline{k}_\pm^\init$ and $\uline{k}_\pm^\term$ by $\uline{k}^\init$ and $\uline{k}^\term$ respectively.
Similarly, for all $k\neq 4$ or $5$, denote the strands $s(\uline{k}_\pm)$ simply by $s(\uline{k})$. 
See Figure~\ref{fig:2-3-0-4-1-mixed}.
\begin{figure}[h]
\centering
\psfrag{0in}{$\uline{0}^{\init}$}
\psfrag{1in}{$\uline{1}^{\init}$}
\psfrag{2in}{$\uline{2}^{\init}$}
\psfrag{3in}{$\uline{3}^{\init}$}
\psfrag{4in}{$\uline{4}^{\init}$}
\psfrag{0te}{$\uline{0}^{\term}$}
\psfrag{1te}{$\uline{1}^{\term}$}
\psfrag{2te}{$\uline{2}^{\term}$}
\psfrag{3te}{$\uline{3}^{\term}$}
\psfrag{4te}{$\uline{4}^{\term}$}
\includegraphics[page=13,scale=0.45]
{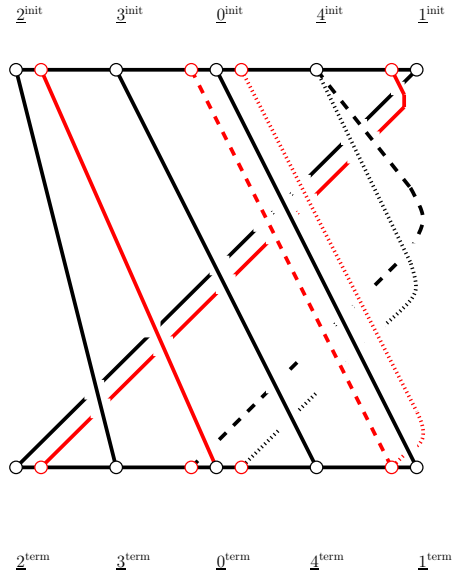}
%{2-3-0-4-1-braid-mixed-b.eps}
\caption{The braids $\beta_\pm$ with strands in the same diagram, from Example~\ref{eg:main-constr4}.
The strands with long dashes come from $\beta_-$, while the strands with short dashes come from $\beta_+$.
}
\label{fig:2-3-0-4-1-mixed}
\end{figure}
Let us do the following to the braids $\beta_-$ and $\beta_+$.
As shown in Figure~\ref{fig:2-3-0-4-1-c2(a)} add an initial endpoint, shown in green, to the left of $\uline{5}_-^\init$ in both diagrams.
Add a strand, also shown in green, from this new initial endpoint which unimodally follows the strand $s(\uline{5}_-)$.
See Figure~\ref{fig:2-3-0-4-1-c2(b)}.
Next, add a new initial endpoint directly above the terminal endpoint of this last strand, which neighbours $\uline{6}^\init$.
From this initial endpoint add a new strand which unimodally follows the next strand $s(\uline{6})$.
See Figure~\ref{fig:2-3-0-4-1-c2(c)}.
Observe that the new terminal endpoint lies inside the interval containing the dynamical preimage.
Add a final initial endpoint directly above this terminal endpoint, necessarily neighbouring $\uline{7}^\init$.
From this initial endpoint add a strand which unimodally follows the strand $s(\uline{7})$ until $s(\uline{7})$ has made all its crossings {\it except} possibly for a single crossing with $s(\uline{4})$.
Make the new strand form a negative crossing with the strand $s(\uline{7})$ before ending at a terminal endpoint directly below the very first initial endpoint we started with.
See Figure~\ref{fig:2-3-0-4-1-c2(d)}.
This gives a pair of unimodal braids, which are non-cyclic (there are exactly two cycles).
To form a pair of cyclic braids, perform a negative half-twist between this final terminal 
%%%%%%%%%%%%%%%%%%%%%%%%%%%%%%%%%%%%%%%%%%%%%%%%%
\begin{landscape}
\begin{figure}[htbp]
\centering
\tiny
\psfrag{0in}{$\uline{0}^\init$}
\psfrag{1in}{$\uline{1}^\init$}
\psfrag{2in}{$\uline{2}^\init$}
\psfrag{3in}{$\uline{3}^\init$}
\psfrag{4in}{$\uline{4}^\init$}
\psfrag{0te}{$\uline{0}^\term$}
\psfrag{1te}{$\uline{1}^\term$}
\psfrag{2te}{$\uline{2}^\term$}
\psfrag{3te}{$\uline{3}^\term$}
\psfrag{4te}{$\uline{4}^\term$}
\begin{subfigure}[b]{0.4\textwidth}
\centering
\includegraphics[page=1,width=\textwidth]
{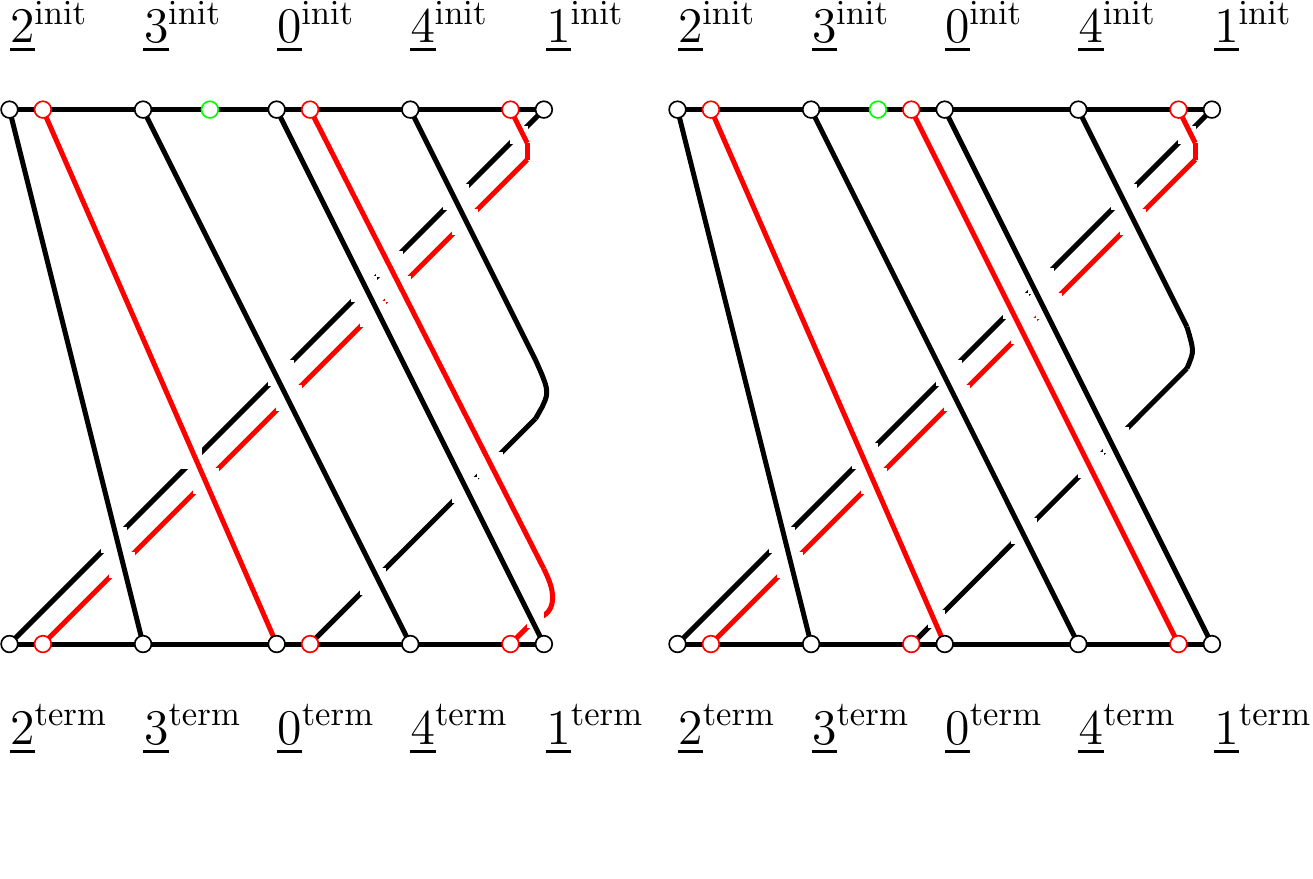}
%{2-3-0-4-1-c2-braid6.eps}
\caption{Adding the first initial point.}
\label{fig:2-3-0-4-1-c2(a)}
\end{subfigure}
~\qquad
\begin{subfigure}[b]{0.4\textwidth}
\centering
\includegraphics[page=2,width=\textwidth]
{BEiH1-figs_epsii-2015-03-30-crop.pdf}
%{2-3-0-4-1-c2-braid5.eps}
\caption{Attaching the first strands.}
\label{fig:2-3-0-4-1-c2(b)}
\end{subfigure}
~\qquad
\begin{subfigure}[b]{0.4\textwidth}
\centering
\includegraphics[page=3,width=\textwidth]
{BEiH1-figs_epsii-2015-03-30-crop.pdf}
%{2-3-0-4-1-c2-braid4.eps}
\caption{Attaching the second strands.}
\label{fig:2-3-0-4-1-c2(c)}
\end{subfigure}

\vspace{5 pt}

\begin{subfigure}[b]{0.4\textwidth}
\centering
\includegraphics[page=4,width=\textwidth]
{BEiH1-figs_epsii-2015-03-30-crop.pdf}
%{2-3-0-4-1-c2-braid3.eps}
\caption{Attaching the final strand.}
\label{fig:2-3-0-4-1-c2(d)}
\end{subfigure}
~\qquad
\begin{subfigure}[b]{0.4\textwidth}
\centering
\includegraphics[page=5,width=\textwidth]
{BEiH1-figs_epsii-2015-03-30-crop.pdf}
%{2-3-0-4-1-c2-braid2.eps}
\caption{Performing the half-twist $\tau$.}
\label{fig:2-3-0-4-1-c2(e)}
\end{subfigure}
%~\qquad
%\begin{subfigure}[b]{0.4\textheight}
%\centering
%\includegraphics[width=\textwidth]{2-3-0-4-1-braid1.eps}
%\caption{....}
%\label{fig:2-3-0-4-1(f)}
%\end{subfigure}
\caption{Adding strands to $\beta_\pm$ from Example~\ref{eg:main-constr3} resulting in the new braids $\beta_-^1$ and $\beta_+^1$.}
\label{fig:eg-2-3-0-4-1-c2}
\end{figure}
\end{landscape}
%%%%%%%%%%%%%%%%%%%%%%%%%%%%%%%%%%%%%%%%%%%%%%%%%%%%%
\noindent
endpoint and $s(\uline{0})^\term$.
There is a complication in that for one of these braids, there is another terminal endpoint, namely $\uline{5}_+^\term$, lying in between: we take the strand from this endpoint going {\it under} the half-twist.
See Figure~\ref{fig:2-3-0-4-1-c2(e)}.  
Denote the resulting braids by $\beta_-^1$ and $\beta_+^1$.
Conjugating $\beta_+^1$ by a positive half-twist between $\uline{0}$ and $\uline{5}_{+}$, the resulting braid can be deformed by a series of isotopies and Reidemeister moves into the braid $\beta_+^1$.
\end{eg}
%CHECK THAT CLOSEST RETURN TIMES HAVE BEEN DEFINED.
%CHECK THAT OPPOSITE POINTS HAVE BEEN DEFINED.
Next, we need the following notions to simplify exposition.
Given a unimodal permutation $\upsilon\in U_p$ and 
$\imath\in\{0,1,\ldots,p-1\}$, 
$\uline{\imath}\neq \uline{0}$, 
we add a point $\uline{\bar{\imath}}$ to the linearly ordered set $\{\uline{0},\uline{1},\ldots,\uline{p-1}\}$, which we call the {\it opposite} of $\uline{\imath}$ with respect to $\upsilon$, 
which satisfies
\begin{enumerate}
\item[(i)] $\uline{\bar{\imath}}<\uline{0}$ if and only if $\uline{\imath}>\uline{0}$ 
\item[(ii)] $\uline{\bar{\imath}}\in (\uline{\jmath},\uline{k})$ if and only if $\upsilon(\uline{\imath})\in (\upsilon(\uline{\jmath}),\upsilon(\uline{k}))$.
\end{enumerate}
\begin{rmk}
Observe that, if the set $\{1,\ldots,p-1\}$ is embedded in the line in an order-preserving manner and $\upsilon$ is realised by a unimodal endomorphism $f$, 
then the opposite of a point $\uline{\imath}$ just corresponds to the preimage of $f(\uline{\imath})$ under $f$ which is not $\uline{\imath}$. 
\end{rmk}
In what follows we only need to consider closest return times to the folding point, so we will simply say that $p$ is the closest return time. (See Section~\ref{subsect:unimodal} for the definition.)

For positive integers $p$ and $q$, let $U_{p,q}$ denote the set of cyclic unimodal permutations of length $p+q$ with closest return time $p$.
Given $\upsilon\in U_{p,q}$, 
let $C=[\uline{p},\uline{\bar{p}}]$, where $\uline{\bar{p}}$ denotes the opposite point of $\uline{p}$ with respect to $\upsilon$. 
Let $D=[\uline{d}^-,\uline{d}^+]$ denote the interval containing the dynamical preimage $\uline{d}=\uline{p+q-1}$ of $\uline{0}$. 
Let $D^{-}=[\uline{d}^{-},\uline{d}]$ and $D^{+}=[\uline{d},\uline{d}^{+}]$.
Let $E=[\uline{e}^-,\uline{e}^+]$ denote the interval containing the non-dynamical preimage.

\begin{eg}\label{eg:main-constr5}
Observe that the unimodal permutations 
constructed in Example~\ref{eg:main-constr3}, and also considered in Example~\ref{eg:main-constr4}, 
\begin{equation}
\upsilon_-=(\uline{2},\uline{7},\uline{3},\uline{5}_-,\uline{0},\uline{4},\uline{6},\uline{1}), 
\ \mbox{and} \
\upsilon_+=(\uline{2},\uline{7},\uline{3},\uline{0},\uline{5}_+,\uline{4},\uline{6},\uline{1})
\end{equation}
both lie in $U_{5,3}$, where we have kept the notation convention of Example~\ref{eg:main-constr4}.
%As suggested by Example~\ref{eg:main-constr4}, when considering the points $\uline{0}_\pm, \uline{1}_\pm,\ldots,\uline{p-1}_\pm$ embedded in the line we may assume $\uline{k}_-=\uline{k}_+$ for all $k\neq 5$, and denote them simply by $\uline{k}$.
%Moreover, observe that $\overline{\uline{5}}_-=\uline{5}_+$ and likewise $\overline{\uline{5}}_+=\uline{5}_-$.
%Also observe that if we remove $\uline{5}_\pm$, the set 
%$\{\uline{0}_-,\uline{1}_-,\ldots,\uline{p-1}_-\}$
%is in order-preserving bijection with
%$\{\uline{0}_+,\uline{1}_+,\ldots,\uline{p-1}_+\}$.
%Hence,we may apply the convention in Example~\ref{eg:main-constr4} of setting $\uline{k}_-=\uline{k}_+$ for all $k\neq 5$.
If we denote the intervals $C, D, D^\pm, E$ for the permutation $\upsilon_\pm$ by $C_\pm,D_\pm,D_\pm^\pm, E_\pm$, then we find that 
$
C_-=[\uline{5}_-,\uline{5}_+]=C_+
$,
$
D_-
%=[\uline{2}_-,\uline{3}_-]
=[\uline{2},\uline{3}]
%=[\uline{2}_+,\uline{3}_+]
=D_+
$,
$
D_-^-
%=[\uline{2}_-,\uline{7}_-]
=[\uline{2},\uline{7}]
%=[\uline{2}_+,\uline{7}_+]
=D_+^-
$,
$
D_-^+
%=[\uline{7}_-,\uline{3}_-]
=[\uline{7},\uline{3}]
%=[\uline{7}_+,\uline{3}_+]
=D_+^+
$,
$
E_-
%=[\uline{0}_-,\uline{4}_-]
=[\uline{0},\uline{4}]
$, and
$
E_+
%=[\uline{4}_+,\uline{6}_+]
=[\uline{4},\uline{6}] 
$.
\end{eg}

Now consider the general situation.
Let $\upsilon_\pm\in U_{p,q}$.
In cyclic notation denote $\upsilon_\pm$ by $(\uline{2}_\pm,\ldots,\uline{0}_\pm,\ldots,\uline{1}_\pm)$.
Assume the following:

\begin{enumerate}
\item[1.] 
(closest returns lie on opposite sides of the folding point)
$\uline{p}_-<\uline{0}_-$ and $\uline{0}_+<\uline{p}_+$.
\item[2.]  
(all remaining points are in order-preserving bijection)
for all $k,l\neq p$, $\uline{k}_-<\uline{l}_-$ if and only if $\uline{k}_+<\uline{l}_+$.
Hence, if $\uline{k}_\pm$ are embedded in the line in an order-preserving manner, as in the preceding example, 
we may assume that $\uline{k}_-=\uline{k}_+$ for all $k\neq p$.
Consequently, for each $k\neq p$, we may denote the point $\uline{k}_\pm$ by $\uline{k}$.
\item[3.]
(the closest return is not contained in the interval containing the dynamical preimage of the folding point)
$\uline{p}_\pm\notin\del D_\pm$.
%$C\cap D=\emptyset$
\item[4.]
(the dynamical preimage of the folding point and the dynamical preimage of the closest return lie on opposite sides of the folding point)
either $\uline{p+q-1}_\pm<\uline{0}_\pm<\uline{p-1}_\pm$,
or $\uline{p-1}_\pm<\uline{0}_\pm<\uline{p+q-1}_\pm$.
\end{enumerate}
Let $C_\pm, D_\pm, D_\pm^\pm, E_\pm$ denote the corresponding intervals for $\upsilon_\pm$.
Then by (1) we may assume that $\uline{p}_-=\overline{\uline{p}}_+$ and hence $C_-=C_+$. 
Denote this interval by $C$.
Properties (2) and (3) then imply that $D_-=D_+$. 
Denote this interval by $D$.
Similarly $D_-^\pm=D_+^\pm$ and we may denote this interval by $D^\pm$.
Note that, as we saw in the previous Example~\ref{eg:main-constr5}, $E_-$ and $E_+$ do not necessarily coincide.

By the discussion in Section~\ref{sect:gen-cabling}, Condition (2) implies the following property, that will be used in the proof of Theorem~\ref{thm:2ndmain} below.
\begin{enumerate}
\item[$2'$.]
Let $\beta_\pm$ denote the canonical unimodal braid of $\upsilon_\pm$. 
Then there exists a braid $\gamma$ with $p+q$ strands, 
with the property that outside of $C=[\uline{p}_{-},\uline{p}_{+}]$, $\gamma$ is trivial, and such that one of the following holds:
\begin{enumerate}
\item[$2_{+}'$.]
$\beta_-\cdot \gamma=\gamma\cdot\beta_+$
\item[$2_{-}'$.] 
$\beta_-\cdot \gamma=\gamma^{-1}\cdot\beta_+$
\end{enumerate}
\end{enumerate}
Finally we will need the following trivial observation
\begin{enumerate}
\item[5.]
(existence of the transit time from the closest return to the interval containing the dynamical preimage of the folding point)
The orbit segment $\uline{p+1},\uline{p+2},\ldots,\uline{p+q-1}$, of the image of the closest return point to the interval $D$ does not intersect the interval $C$.
\end{enumerate}
Call $t=q-2$ the {\it transit time}. 
We say it is a {\it left transit time} if the parity, given by $\rho(\uline{p+1},\uline{p+t+1})$, is even and a {\it right transit time} if this parity is odd.
\begin{eg}
Let us continue consider $\upsilon_\pm$ from Example~\ref{eg:main-constr5}.
Then the transit time is $t=1$, since $q=3$. It is a right transit time. If we imagine a point neighbouring $\uline{p+1}$ to the left, then the image of that point lies in the interval containing the dynamical preimage, i.e. we land inside the interval containing the dynamical preimage after a single iterate. Moreover, this image lies to the right of the dynamical preimage of the folding point, therefore it is a right transit time. 
\end{eg}

\begin{rmk}
Any pair $(\upsilon_-,\upsilon_+)$ of equivalent unimodal permutations from Construction~\ref{constr:gencabling2} satisfies properties (1)--(5). 
In particular, the unimodal permutations from Example~\ref{eg:main-constr5} satisfy these conditions.
\end{rmk}

\begin{constr}[Generalised Cabling Process.]\label{constr:gencablingprocess1}
Let $(\upsilon_-,\upsilon_+)$ denote a pair of cyclic unimodal permutations satisfying properties (1)--(4).
Recall that $C_\pm, D_\pm$, etc. denote the corresponding intervals for $\upsilon_\pm$.
%We may assume $\uline{p}_-=\overline{\uline{p}}_+$. 
%Hence we can suppose that $C_-=C_+$ and we denote this simply by $C$.
%We also assume that $D_-=D_+$ and denote it by $D$.
%This can be assumed provided that $p_\pm$ does not lie on the boundary of $D_\pm$.
Recall that we may assume $C_-=C_+$ and $D_-=D_+$.
Denote these intervals by $C$ and $D$ respectively.
Observe that $\upsilon_-(D)=\upsilon_+(D)$ and this interval contains $\uline{0}$ and $\uline{p}_\pm$ in its interior.

Let $\beta_\pm$ denote the canonical braid corresponding to $\upsilon_{\pm}$.
Denote the set of strands for $\beta_\pm$ by $S_\pm$. 
For $k=0,1,\ldots,p+q-1$, denote the $k$-th strand by $s_\pm(\uline{k})$.
The strand $s_\pm(\uline{k})$ lies between initial endpoint $s_{\pm}(\uline{k})^\init=\uline{k}_\pm^\init$ and terminal endpoint $s_{\pm}(\uline{k})^\term=\uline{k+1}_\pm^\term$, where addition is taken modulo $p+q$. 
By property (4) above we may assume, after applying an isotopy if necessary, that $s_-(\uline{k})=s_+(\uline{k})$ for all $k\neq p-1$ or $p$.
%\begin{rmk}
%By the above it follows that $\uline{k}_{-}=\uline{k}_{+}$ for $k\neq p$. 
%We therefore drop $\pm$ from this notation and denote these points by $\uline{k}$.
%\end{rmk}
\begin{rmk}
We will also use the notation $s_k$ for the projection to $[-1,1]$ of the initial endpoint of $s(\uline{k})$ for $k=0,1,\ldots,p_-,p_+,\ldots,p+q-1$.
%The reason for this is to make the distinction between the endpoints of the strands from the original braids $\beta_\pm$ and those of new strands added during the construction.
\end{rmk}

Let $t$ be the transit time. Construct a pair $(\upsilon_-^1,\upsilon_+^1)$ of unimodal permutations as follows:

%def of additional points
To begin, we define points $r_i$, $r_i^\init$ and $r_i^\term$ for $i=0, 1,\ldots, t+1$.
Choose $\epsilon>0$ small.
Let $r_1\in (\uline{p+1}-\epsilon,\uline{p+1})$,  
and choose points $r_i\in(\uline{p+i}-1,\uline{p+i}+1)$, $i=2,\ldots, t+1$, which satisfy
\begin{equation}
r_i\in\left\{
\begin{array}{ll}
(\uline{p+i}-\epsilon,\uline{p+i}) & \mbox{if} \ \rho(\uline{p+1},\uline{p+i}) \ \mbox{even} \\
(\uline{p+i},\uline{p+i}+\epsilon) & \mbox{if} \ \rho(\uline{p+1},\uline{p+i}) \ \mbox{odd}
\end{array}
\right.
\end{equation}
Finally, choose $r_0\in (\uline{p}_{-}-1,\uline{p}_{+}+1)$ which satisfies
\begin{equation}
r_0\in\left\{
\begin{array}{ll}
(\uline{p}_{-}-\epsilon,\uline{p}_{-}) & 
\mbox{if} \ D<\uline{0} \ \mbox{and} \ t \ \mbox{a right time, or} \
\uline{0}<D \ \mbox{and} \ t \ \mbox{a left time.}\\
(\uline{p}_{+},\uline{p}_{+}+\epsilon) & 
\mbox{if} \  D<\uline{0} \ \mbox{and} \ t \ \mbox{a left time, or} \
\uline{0}<D \ \mbox{and} \ t \ \mbox{a right time.}
\end{array}
\right.
\end{equation}
For $i=0,1,2,\ldots, t+1$ let $r_i^\init$ and $r_i^\term$ denote the projections of $r_i$ onto $[-1,1]^\init$ and $[-1,1]^\term$ respectively. 
As before, these will be endpoints for the strands constructed inductively below.

%def of additional strands
First, let $r_\pm(\uline{1})$ denote the strand with initial endpoint $r_1^\init$ and terminal endpoint $r_2^\term$ which follows unimodally the strand $s_\pm(\uline{p+1})$.

Secondly, assume that the strands $r_\pm(\uline{i})$ have been constructed for $i=1,\ldots,j-1$, for some $0<j<t+1$.
Let $r_\pm(\uline{j})$ denote the strand with initial endpoint $r_j^\init$ and terminal endpoint $r_{j+1}^\term$ which follows unimodally the strand $s_\pm(\uline{p+j})$. 

Thirdly, assume that the strands $r_\pm(\uline{i})$ have been constructed for all $i=1,2,\ldots,t$. 
Then $r_{\pm}(\uline{t+1})$ denotes the strand  with initial endpoint $r_{t+1}^\init$ and terminal endpoint $r_{0}^\term$ which, in each of the four cases below, satisfies the following: 
\begin{enumerate}
\item[(Ii)]
{\it [If $\uline{p+q-1}<\uline{0}$ and $t$ is a left transit time.]}%i.e. D<E and...

Let $r_{-}(\uline{t+1})$ unimodally follow the strand $s_{-}(\uline{p+t+1})$ rightwardly until $s_{-}(\uline{p+t+1})$ has only one crossing to make, which is necessarily with $s_{-}(\uline{p-1})$, before reconnecting to its terminal endpoint.
Then $r_{-}(\uline{t+1})$ forms a {\it negative} crossing with $s_{-}(\uline{p+t+1})$, after which it makes a {\it positive} crossing with $s_{-}(\uline{p-1})$ before connecting to its terminal endpoint.
See Figure~\ref{fig:4I-ia}.

Let $r_{+}(\uline{t+1})$ unimodally follow the strand $s_{+}(\uline{p+t+1})$ rightwardly until $s_{+}(\uline{p+t+1})$ has no more crossings to make before reconnecting to its terminal endpoint. 
Then $r_{+}(\uline{t+1})$ forms a {\it negative} crossing with $s_{+}(\uline{p+t+1})$, after which it makes a {\it positive} crossing with $s_{+}(\uline{p-1})$ before connecting to its terminal endpoint. 
See Figure~\ref{fig:4I-ib}.
\item[(Iii)]
{\it [If $\uline{p+q-1}<\uline{0}$ and $t$ a right transit time.]}%i.e. D<E and...

Let $r_{-}(\uline{t+1})$ unimodally follow the strand $s_{-}(\uline{p+t+1})$ rightwardly until $s_{-}(\uline{p+t+1})$ has only one crossing to make, which is necessarily with $s_{-}(\uline{p-1})$, before reconnecting to its terminal endpoint.
Then $r_{-}(\uline{t+1})$ forms a {\it negative} crossing with $s_{-}(\uline{p+t+1})$ before connecting to its terminal endpoint.
See Figure~\ref{fig:4I-iia}.

Let $r_{+}(\uline{t+1})$ unimodally follow the strand $s_{+}(\uline{p+t+1})$ rightwardly until $s_{+}(\uline{p+t+1})$ has no further crossings to make before reconnecting to its terminal endpoint.
Then $r_{+}(\uline{t+1})$ forms a {\it negative} crossing with $s_{+}(\uline{p+t+1})$ before connecting to its terminal endpoint.
See Figure~\ref{fig:4I-iib}.
\item[(IIi)]
{\it [If $\uline{0}<\uline{p+q-1}$ and $t$ a left transit time.]}%i.e. E<D and...

Let $r_{-}(\uline{t+1})$ unimodally follow $s_{-}(\uline{p+t+1})$ until $s_{-}(\uline{p+t+1})$ has no more crossings to make before reconnecting to its terminal endpoint. 
Then $r_{-}(\uline{t+1})$ forms a {\it negative} crossing with $s_{-}(\uline{p+t+1})$, after which it makes a {\it positive} crossing with $s_{-}(\uline{p-1})$ before connecting to its terminal endpoint. 
See Figure~\ref{fig:4II-ia}.

Let $r_{+}(\uline{t+1})$ unimodally follow $s_{+}(\uline{p+t+1})$ until $s_{+}(\uline{p+t+1})$ has only one crossing to make, which is necessarily with $s_{+}(\uline{p-1})$, before reconnecting to its terminal endpoint. 
Then $r_{+}(\uline{t+1})$ forms a {\it positive} crossing with $s_{+}(\uline{p-1})$, after which it makes a {\it negative} crossing with $s_{+}(\uline{p+t+1})$ before connecting to its terminal endpoint.
See Figure~\ref{fig:4II-ib}.
\item[(IIii)]
{\it [If $\uline{0}<\uline{p+q-1}$ and $t+1$ a right transit time.]}%i.e. E<D and...

Let $r_{-}(\uline{t+1})$ unimodally follow $s_{-}(\uline{p+t+1})$ until $s_{-}(\uline{p+t+1})$ has no further crossings to make before connecting to its terminal endpoint.
Then $r_{-}(\uline{t+1})$ forms a single {\it negative} crossing with $s_{-}(\uline{p+t+1})$ before connecting to its terminal endpoint.
See Figure~\ref{fig:4II-iia}.

Let $r_{+}(\uline{t+1})$ unimodally follow $s_{+}(\uline{p+t+1})$ until $s_{+}(\uline{p+t+1})$ has only one crossing to make, which is necessarily with $s_{+}(\uline{p-1})$, before connecting to its terminal endpoint.
Then $r_{+}(\uline{t+1})$ forms a single {\it negative} crossing with $s_{+}(\uline{p+t+1})$ before connecting to its terminal endpoint.
See Figure~\ref{fig:4II-iib}.
\end{enumerate}
Finally, let $r_{\pm}(\uline{0})$ denote the strand with initial endpoint $r_0^\init$ and terminal endpoint $r_1^\term$ which, in cases (Ii) and (IIii), follows unimodally the strand  $s(\uline{p}_{+})$ and, in cases (Iii) and (IIi), follows unimodally the strand $s(\uline{p}_{-})$.  

Let $\alpha_\pm^1\in B_{p+q+t+2}$ denote the braid consisting of the strands of $\beta_\pm$ together with the strands $r_\pm(\uline{i})$, $i=0,1,\ldots,t+1$, formed above.
(As before, after relabelling the endpoints this gives a well-defined braid diagram normalised as in Section~\ref{sect:unimodalbraids}. However, for the moment we will keep the labelling as above.)
%%%%%%%%%%%%%%%%%%%%%%%%%%%%%%%%%%%%%%%%%%%%%%%%%%%%%%%%%%%%%%%%%%%%
%\begin{landscape}
\begin{figure}
\tiny
\psfrag{sp-1in}{$s_{p-1}^\init$}
\psfrag{sp+q-1in}{$s_{p+t+1}^\init$}
\psfrag{rt1in}{$r_{t+1}^\init$}
\psfrag{s0in}{$s_{0}^\init$}
\psfrag{sp-in}{$s_{p_{-}}^\init$}
\psfrag{sp+in}{$s_{p_{+}}^\init$}
\psfrag{r0in}{$r_{0}^\init$}
\psfrag{rt+1in}{$r_{t+1}^\init$}
\psfrag{sp+t+1in}{$s_{p+t+1}^\init$}
\psfrag{sigma}{$\sigma$}
\psfrag{alpha-1}{$\alpha_-^1$}
\psfrag{tau-1}{$\tau_-^1$}
\psfrag{alpha+1}{$\alpha_+^1$}
\psfrag{tau+1}{$\tau_+^1$}
\psfrag{sigmainv}{$\sigma^{-1}$}
\psfrag{delta-1}{$\delta_-^{1}$}
\psfrag{epsilon-1}{$\varepsilon_-^{1}$}
\psfrag{delta+1}{$\delta_+^{1}$}
\psfrag{epsilon+1}{$\varepsilon_+^{1}$}
\psfrag{sp-te}{$s_{p_{-}}^\term$}
\psfrag{s0te}{$s_{0}^\term$}
\psfrag{sp+te}{$s_{p_{+}}^\term$}
\psfrag{sp+1te}{$s_{p+1}^\term$}
\psfrag{r0te}{$r_{0}^\term$}
\psfrag{r1te}{$r_{1}^\term$}
\psfrag{s1te}{$s_{1}^\term$}
\begin{subfigure}[b]{0.45\textwidth}
\centering
\includegraphics[page=6,width=\textwidth]
{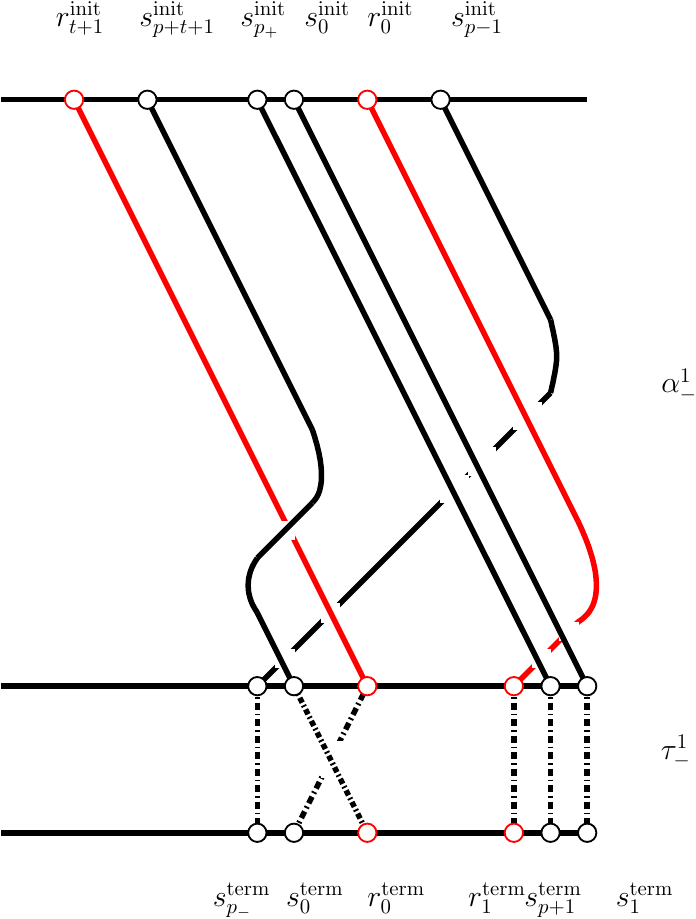}
%{braidsv2-cabl4I-ia.eps}
\caption{$\beta_-$ in Construction~\ref{constr:gencablingprocess1}(Ii).}
\label{fig:4I-ia}
\end{subfigure}
~\qquad
\begin{subfigure}[b]{0.45\textwidth}
\centering
\includegraphics[page=7,width=\textwidth]
{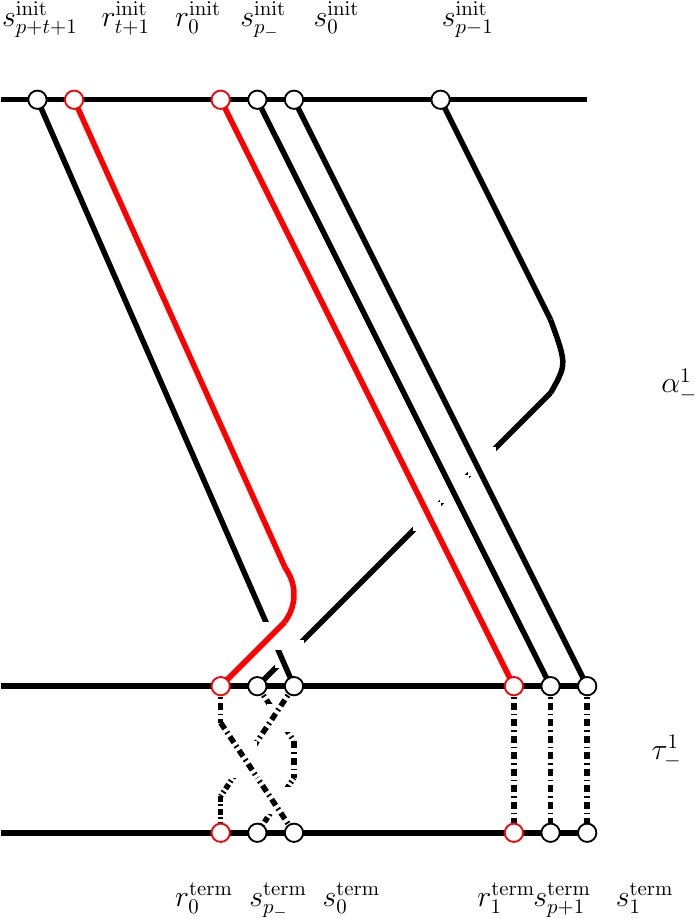}
%{braidsv2-cabl4I-iia.eps}
\caption{$\beta_-$ in Construction~\ref{constr:gencablingprocess1}(Iii).}
\label{fig:4I-iia}
\end{subfigure}
~
\vspace{10pt}

\begin{subfigure}[b]{0.45\textwidth}
\centering
\includegraphics[page=8,width=\textwidth]
{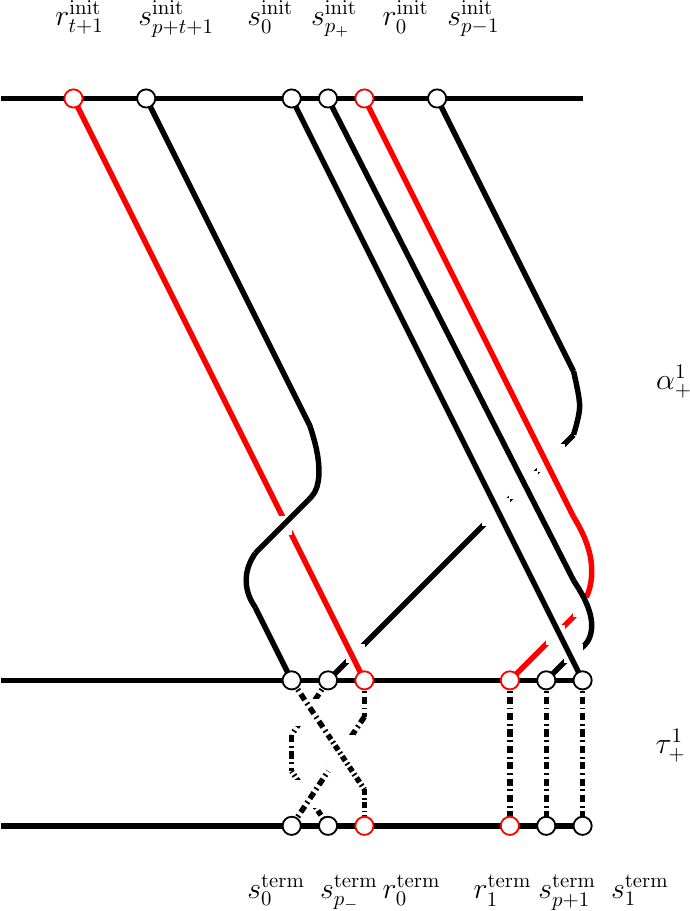}
%{braidsv2-cabl4I-ib.eps}
\caption{$\beta_+$ in Construction~\ref{constr:gencablingprocess1}(Ii).}
\label{fig:4I-ib}
\end{subfigure}
~\qquad
\begin{subfigure}[b]{0.45\textwidth}
\centering
\includegraphics[page=9,width=\textwidth]
{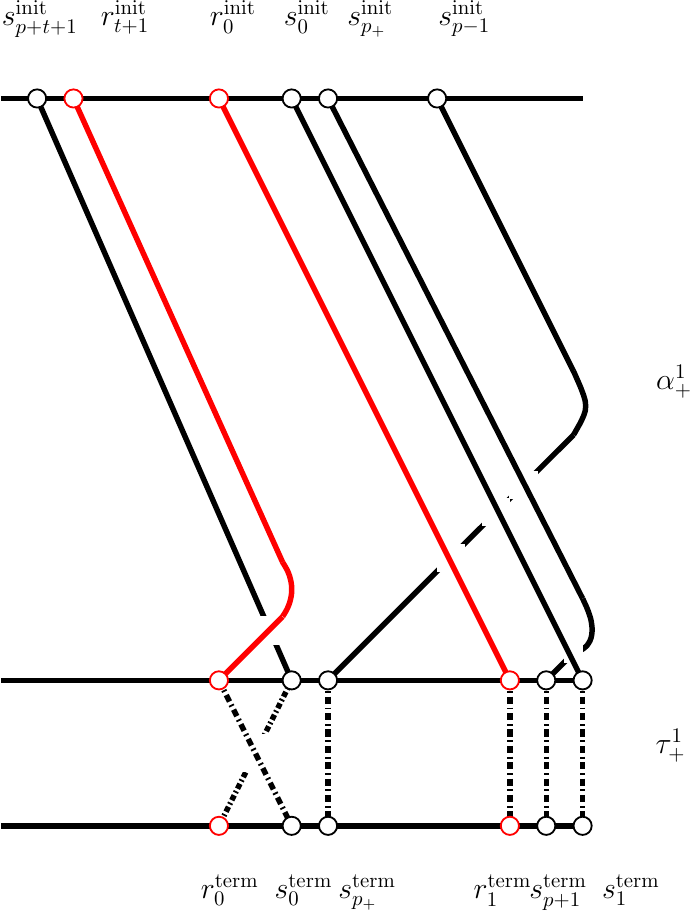}
%{braidsv2-cabl4I-iib.eps}
\caption{$\beta_+$ in Construction~\ref{constr:gencablingprocess1}(Iii).}
\label{fig:4I-iib}
\end{subfigure}
\caption{$\beta_-$ and $\beta_+$ in Construction~\ref{constr:gencablingprocess1}(I).}
\label{fig:4I}
\end{figure}
%\end{landscape}
%%%%%%%%%%%%%%%%%%%%%%%%%%%%%%%%%%%%%%%%%%%%%%%%%%%%%%%%%%%%%%%%%%
%\begin{landscape}
\begin{figure}
\tiny
\psfrag{sp-1in}{$s_{p-1}^\init$}
\psfrag{sp+q-1in}{$s_{p+t+1}^\init$}
\psfrag{rt1in}{$r_{t+1}^\init$}
\psfrag{s0in}{$s_{0}^\init$}
\psfrag{sp-in}{$s_{p_{-}}^\init$}
\psfrag{sp+in}{$s_{p_{+}}^\init$}
\psfrag{rt+1in}{$r_{t+1}^\init$}
\psfrag{sp+t+1in}{$s_{p+t+1}^\init$}
\psfrag{r0in}{$r_{0}^\init$}
\psfrag{sigma}{$\sigma$}
\psfrag{alpha-1}{$\alpha_-^1$}
\psfrag{tau-1}{$\tau_-^1$}
\psfrag{alpha+1}{$\alpha_+^1$}
\psfrag{tau+1}{$\tau_+^1$}
\psfrag{sigmainv}{$\sigma^{-1}$}
\psfrag{sp-te}{$s_{p_{-}}^\term$}
\psfrag{s0te}{$s_{0}^\term$}
\psfrag{sp+te}{$s_{p_{+}}^\term$}
\psfrag{sp+1te}{$s_{p+1}^\term$}
\psfrag{r0te}{$r_{0}^\term$}
\psfrag{r1te}{$r_{1}^\term$}
\psfrag{s1te}{$s_{1}^\term$}
\begin{subfigure}[b]{0.45\textwidth}
\centering
\includegraphics[page=10,width=\textwidth]
{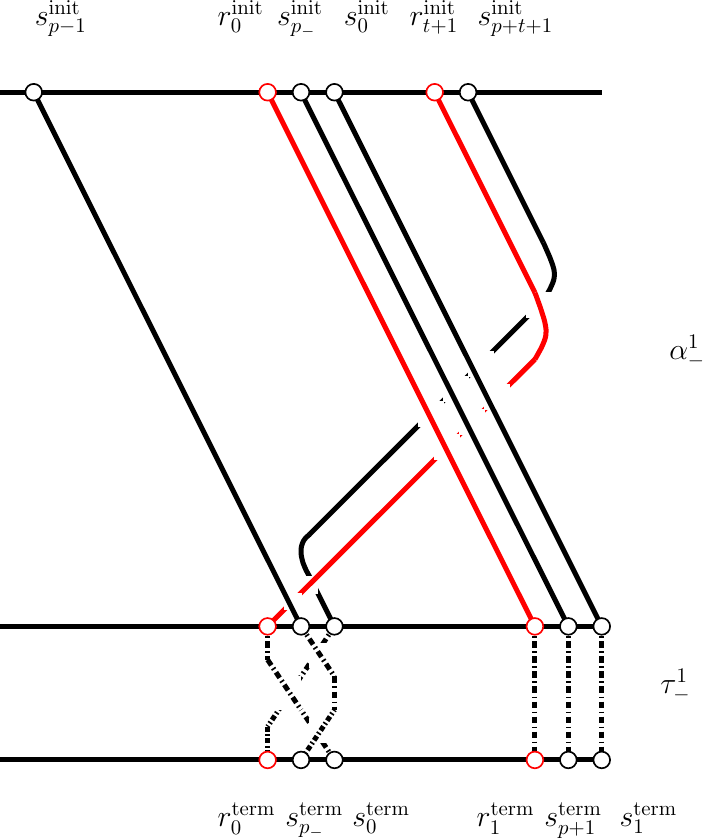}
%{braidsv2-cabl4II-ia.eps}
\caption{$\beta_-$ in Construction~\ref{constr:gencablingprocess1}(IIi).}
\label{fig:4II-ia}
\end{subfigure}
~\qquad
\begin{subfigure}[b]{0.45\textwidth}
\centering
\includegraphics[page=11,width=\textwidth]
{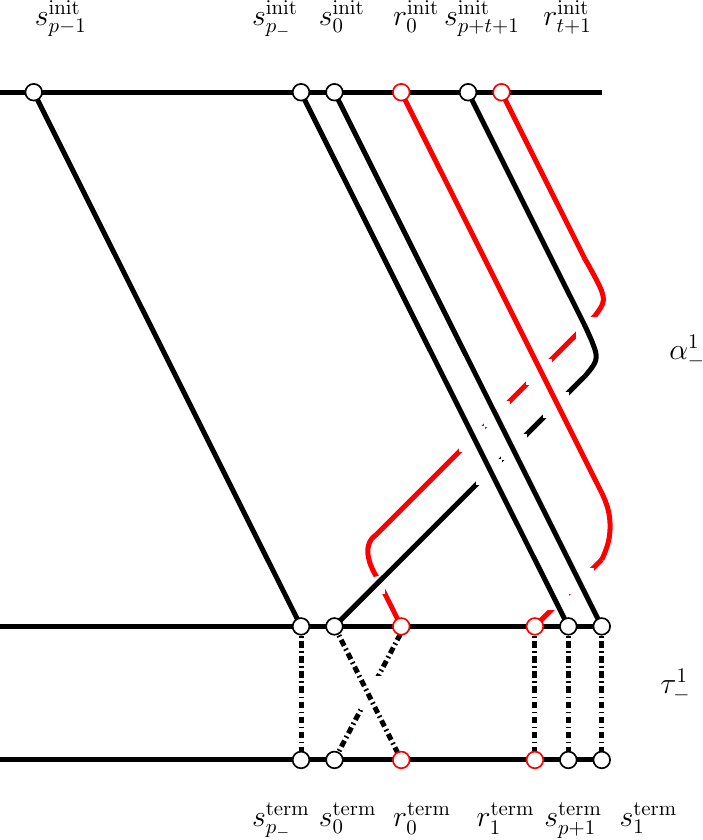}
%{braidsv2-cabl4II-iia.eps}
\caption{$\beta_-$ in Construction~\ref{constr:gencablingprocess1}(IIii).}
\label{fig:4II-iia}
\end{subfigure}
~
\vspace{10pt}

\begin{subfigure}[b]{0.45\textwidth}
\centering
\includegraphics[page=12,width=\textwidth]
{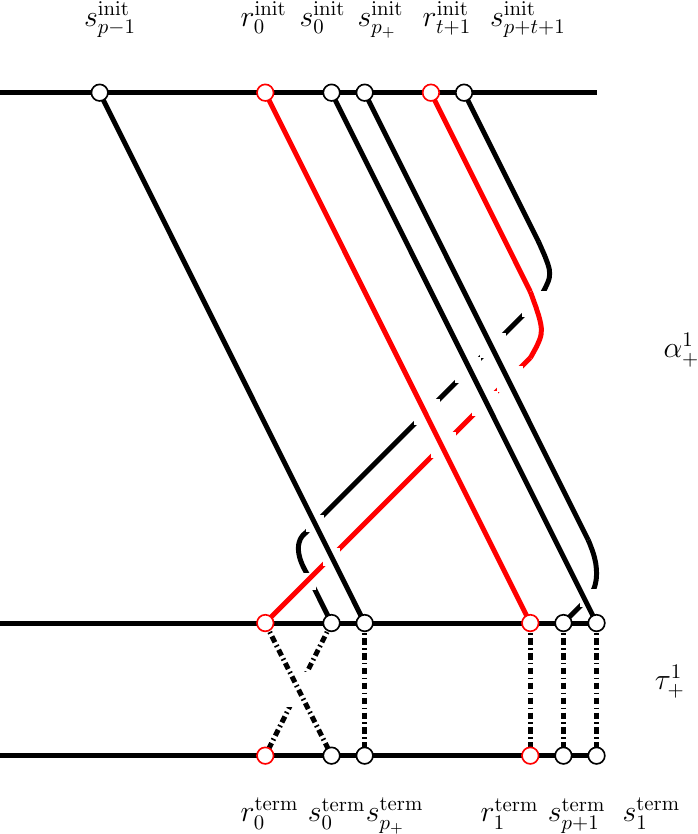}
%{braidsv2-cabl4II-ib.eps}
\caption{$\beta_+$ in Construction~\ref{constr:gencablingprocess1}(IIi).}
\label{fig:4II-ib}
\end{subfigure}
~\qquad
\begin{subfigure}[b]{0.45\textwidth}
\centering
\includegraphics[page=13,width=\textwidth]
{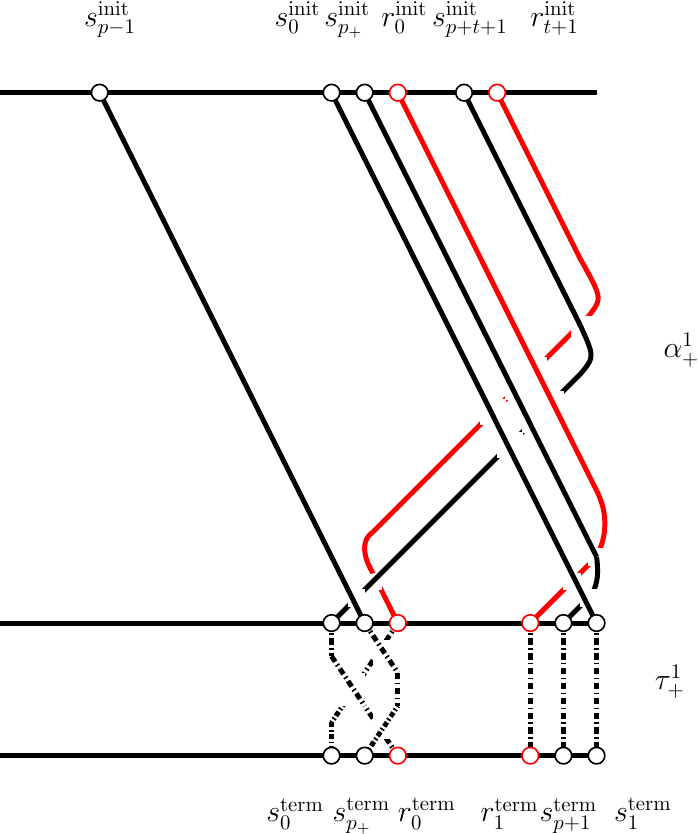}
%{braidsv2-cabl4II-iib.eps}
\caption{$\beta_+$ in Construction~\ref{constr:gencablingprocess1}(IIii).}
\label{fig:4II-iib}
\end{subfigure}
\caption{$\beta_-$ and $\beta_+$ in construction~\ref{constr:gencablingprocess1}(II).}
\label{fig:4II}
\end{figure}
%\end{landscape}
%%%%%%%%%%%%%%%%%%%%%%%%%%%%%%%%%%%%%%%%%%%%%%%%%%%%%%%%%%%%%%%%%%%%%
Note that $\alpha_\pm^1$ is not a cyclic braid.
Note also that $\alpha_\pm^1$ is not necessarily a unimodal braid.
In fact, in all cases $\alpha_\pm^1$ is direct but neither positive nor unimodal.
See Figures~\ref{fig:4I-ia}--\ref{fig:4II-iib}.

As in the previous construction, we now compose $\alpha_\pm^1$ with a positive half-twist $\tau_\pm^1$. 
However, in terms of the diagram there is an extra complication due to either $s_{p_-}$ or $s_{p_+}$ lying between $s_0$ and $r_0$. 
(Recall that, for each $k$,  $s_k$ denotes the projection to $[-1,1]$ of the initial endpoint of $s(\uline{k})$.)
Therefore let $\tau(i,i+1)$ denote the Artin generator positively interchanging $i$ and $i+1$. 
Then
\begin{enumerate}
\item[(Ii)]
$\tau_-^1$ denotes a single positive half-twist  $\tau(s_0,r_0)$ between $s_0$ and $r_0$.

$\tau_+^1$ denotes $\tau(s_0,s_{p_+})^{-1}\cdot \tau(s_{p_+},r_0)\cdot \tau(s_0,s_{p_+})$.
\item[(Iii)]
$\tau_-^1$ denotes $\tau(s_{p_-},s_0)\cdot\tau(r_0,s_{p_-})\cdot \tau(s_{p_-},s_0)^{-1}$.

$\tau_+^1$ denotes a single positive half-twist $\tau(r_0,s_0)$ between $r_0$ and $s_0$. 
\item[(IIi)]
$\tau_-^1$ denotes $\tau(s_{p_-},s_0)^{-1}\cdot\tau(r_0,s_{p_-})\cdot\tau(s_{p_-},s_0)$.

$\tau_+^1$ denotes a single positive half-twist $\tau(r_0,s_0)$ between $r_0$ and $s_0$.
\item[(IIii)]
$\tau_-^1$ denotes a single positive half-twist $\tau(s_0,r_0)$ between $s_0$ and $r_0$.

$\tau_+^1$ denotes $\tau(s_{p_+},r_0)^{-1}\cdot\tau(s_0,s_{p_+})\cdot\tau(s_{p_+},r_0)$
\end{enumerate} 
Let $\beta_\pm^1=\tau_\pm^1\cdot\alpha_\pm^1$. 
Note that $\beta_\pm^1$ is a cyclic braid.
Moreover, after cancelling appropriate crossings between strands $r_\pm(\uline{t+1})$, $s_\pm(\uline{p+t+1})$ and $s_\pm(\uline{p-1})$ via Reidemeister moves and isotopy, we find that $\beta_-^1$ and $\beta_+^1$ are both unimodal.
Again, see Figures~\ref{fig:4I-ia}--\ref{fig:4II-iib}.
\end{constr}

Let use investigate what is going on in more detail.
Consider the case (Ii).
Then $\beta_-$ and $\beta_+$ are depicted in Figures~\ref{fig:4I-ia} and~\ref{fig:4I-ib} respectively.
dentifying the sets of initial and terminal endpoints of $\beta_-$ and $\beta_+$ in an order-preserving manner, then precomposing $\beta_-$ with a positive half-twist between the strand $s(\uline{0})$ and its neighbour $s(\uline{p}_-)$, and post-composing by the inverse of this twist, yields the braid $\beta_+$.
A similar argument can be given in cases (Iii)--(IIii).
See Figures~\ref{fig:4I-iia}--~\ref{fig:4II-iib}.
However, in cases (IIi) and (IIii) we pre-compose and post-compose by the same positive half-twist.
Hence we have the following.
\begin{thm}\label{thm:2ndmain}
Given $(\upsilon_-,\upsilon_+)$, satisfying properties (1)--(4). 
Let $\beta_-$ and $\beta_+$ denote the corresponding unimodal braids.
Let $(\beta_-^1,\beta_+^1)$ denote the pair of braids produced from Construction~\ref{constr:gencablingprocess1}.
\begin{enumerate}
\item 
If $\beta_-\sim\beta_+$ then $\beta_-^1\sim \beta_+^1$.
\item
If $\beta_-\sim_r\beta_+$ then $\beta_-^1\sim_r\beta_+^1$.
\end{enumerate}
Moreover, the corresponding unimodal permutations $(\upsilon_-^1,\upsilon_+^1)$ also satisfy the properties (1)--(4). 
\end{thm}

\section{Applications to the H\'enon Family}
We have constructed two mechanisms for constructing braid equivalences.
However, we would like to restrict ourselves to equivalences realised in the H\'enon family.

%%%%%%%%NEW%%%%%%%%%%%%%
Before proceeding, let us give a brief description of the parameter space of the quadratic family and the H\'enon family.
Recall from the introduction that $f_a(x)=a-x^2$ denotes the quadratic family and $F_{a,b}(x,y)=(a-x^2-by,x)$ denotes the H\'enon family.
In the $(a,b)$-plane, for $b$ positive the map $F_{a,b}$ is an orientation-preserving diffeomorphism.
For $b\in (0,1)$ the map $F_{a,b}$ is area-contracting. 
In fact, $\mathrm{Jac}(F_{a,b})(x,y)=b$ for all $(a,b)\in\RR^2$ and $(x,y)\in\RR^2$.

The parabola
\begin{equation}
(1+b)^2+4a=0
\end{equation}
is the saddle-node bifurcation locus. 
For all parameters $(a,b)$ on this curve $F_{a,b}$ possesses a unique fixed point.
The parabola passes through the $a$-axis at $a=-1/4$.
All parameters to the left of this curve possess no fixed point and the iterates of all points escape to infinity.
All parameters to the right of the curve possess two fixed points, one saddle and one sink, and hence the non-wandering set is non-trivial.

The curve
\begin{equation}\label{eq:Devaney-Nitecki}
a=(5+2\sqrt{5})(1+b)^2
\end{equation}
lies to the right of the saddle-node bifurcation curve above. 
It was shown by Devaney-Nitecki~\cite{DevaneyNitecki79} that for all parameters $(a,b)$ lying to the right of this curve the map $F_{a,b}$ possesses a full horseshoe, and all points either escape to infinity or converge to this invariant set under iteration.  
%Actually, the homoclinic bifurcation parameter $a=2$ extends in the $(a,b)$-plane to a pinched laminar set of bifurcation curves (see Holmes and Whitley~\cite{HolmesWhitley84}).

For a fixed $b$ sufficiently small, increasing $a$ from the saddle-node bifurcation locus to the horseshoe locus the map $F_{a,b}$ undergoes a period-doubling cascade. 
Each period-doubling bifurcation curve is algebraic and they accumulate upon an analytic curve which intersects the $a$-axis at the Feigenbaum-Collet-Tresser parameter $a=1.401...$. 
For parameters to the left of the accumulation of period-doubling, the map $F_{a,b}$ has simple dynamics: there are finitely many periodic orbits each of period $2^n$ for some non-negative integer $n$.   

After this accumulation of period-doubling less is known.
However, restricting to the $a$-axis we know more.
Here between the accumulation of period-doubling and the horseshoe locus, uncountably many bifurcations occur. 
For each periodic kneading sequence there corresponds a hyperbolic component, 
i.e. an interval such that for each parameter in this interval $F_{a,b}$ has a periodic attractor whose itinerary is determined by the given periodic kneading sequence.
For example, there is an interval around the point $a=1.7549..$, for which every parameter has an attractive cycle of period three. 
(This parameter is actually the critically-periodic parameter, or centre, of the unique period-three hyperbolic component.)

These hyperbolic intervals extend to open subsets of the $(a,b)$-plane, where the periodic attractor persists.
The loci of all such parameters, for fixed periods, were first considered by El-Hamouly and Mira~\cite{ElHamoulyMira82}. 
Many components of this locus have the following structure: there exists a main `body' out of which four `limbs' emanate.
Then limbs do not intersect; the union of the limbs and body is simply connected; two of the limbs intersect $\{b=1\}$; the two remaining limbs intersect $\{b=0\}$.
(Similar configurations have been observed for the one-dimensional cubic family. 
These configurations have been called {\it swallow configurations} by Milnor~\cite{Milnor92}.)
Numerical investigations into the braid-equivalences exhibited in the $(a,b)$-plane were carried out by Holmes~\cite{Holmes86}.
(See also Sannami~\cite{Sannami88}.)
However, currently very little is understood about these configurations. 
Their apparent prevalence, in the chaotic parameter region for the H\'enon family as well as in other families, also requires explanation.
\begin{rmk}
The braid equivalences we constructed were based on the initial unimodal permutation $\upsilon$ possessing a non-dynamical preimage. 
This can only be satisfied if the corresponding kneading sequence satisfies $\kappa\succ 10\mathrm{C}$ or, equivalently, has hyperbolic parameter interval lying to the right of the period-three hyperbolic parameter interval. 
Consequently, all numerical example given below intersect the $a$-axis in the interval $[1.7549...,2]$.
However, in~\cite{deCarvalhoHallHazard15} we will describe a generalisation of the construction of braid equivalences given here which do not have this restriction. 
\end{rmk}
%%%%%%%%%%%%%%%%%%%%%%%%%%%%
%%%%%%%%%%%%%%%%%%%%%%%%%%%%

Following the numerical evidence given below, we ask the following questions.
Let $\upsilon_-$ and $\upsilon_+$ be an arbitrary pair of combinatorial types from Construction~\ref{constr:gencabling2}.
\begin{quote}
{\it Question A.}
Let $a_-, a_+\in [-1/4,2]$ be such that $f_{a_-}$ and $f_{a_+}$ have critical orbits $c_-$ and $c_+$ of types $\upsilon_-$ and $\upsilon_+$ respectively. 
Let $C_-$ and $C_+$ denote the corresponding periodic orbits for $F_{a_-,0}$ and $F_{a_+,0}$ respectively. 
Does there exist a braid equivalence in the family $F_{a,b}$ connecting $(C_-,F_{a_-,0})$ and $(C_+,F_{a_+,0})$?
\end{quote}
The above question only deal with braid equivalences coming from Construction~\ref{constr:gencabling2}.
Now, given an initial braid equivalent pair $(\upsilon_-^0,\upsilon_+^0)$, let us consider a sequence of equivalent pairs $(\upsilon_-^i,\upsilon_+^i)$ coming from Construction~\ref{constr:gencablingprocess1}.
\begin{quote}
{\it Question B.}
For each positive integer $i$, let $a_-^i, a_+^i\in [-1/4,2]$ be parameters such that $f_{a_-^i}$ and $f_{a_+^i}$ have critical orbits $c_-^i$ and $c_+^i$ of types $\upsilon_-^i$ and $\upsilon_+^i$ respectively.
Let $C_-^i$ and $C_+^i$ denote the corresponding periodic orbits for $F_{a_-^i,0}$ and $F_{a_+^i,0}$ respectively.
Does there exist, for each $i$,  a braid equivalence in the family $F_{a,b}$ connecting $(C_-^i,F_{a_-^i,0})$ and $(C_+^i,F_{a_+^i,0})$?
Are the paths $\gamma^i$ in the $(a,b)$-plane which realise these braid equivalences pairwise disjoint?
\end{quote}

%\section{Numerical Supporting Evidence}
\newcommand{\A}{\ensuremath{a}}
\newcommand{\B}{\ensuremath{b}}
\newcommand{\X}{\ensuremath{x}}
\newcommand{\Y}{\ensuremath{y}}
\newcommand{\AZEROINIT}{\ensuremath{a_0}}
\newcommand{\BZEROINIT}{\ensuremath{b_0}}
\newcommand{\AONEINIT}{\ensuremath{a_1}}
\newcommand{\BONEINIT}{\ensuremath{b_1}}
\newcommand{\XINIT}{\ensuremath{x_{\mathrm{in}} }}
\newcommand{\YINIT}{\ensuremath{y_{\mathrm{in}} }}
\newcommand{\XZERO}{\ensuremath{x_{0} }}
\newcommand{\YZERO}{\ensuremath{y_{0} }}
\newcommand{\XRIGHT}{\ensuremath{x_{\mathrm{r}} }}
\newcommand{\YRIGHT}{\ensuremath{y_{\mathrm{r}} }}
\newcommand{\XYSEP}{\ensuremath{\delta_{\mathrm{sep}} }}
\newcommand{\ABSEP}{\ensuremath{\epsilon_{\mathrm{sep}} }}
\newcommand{\PMIN}{\ensuremath{p_{\mathrm{min}} }}
\newcommand{\PMAX}{\ensuremath{p_{\mathrm{max}} }}
\renewcommand{\P}{\ensuremath{p}}
\newcommand{\Q}{\ensuremath{q}}
\newcommand{\MAXIT}{\ensuremath{N_\mathrm{transit}}}

\begin{figure}[htp]
\centering
\includegraphics
%[bb=0 0 1080 1080]
[scale=0.3]
{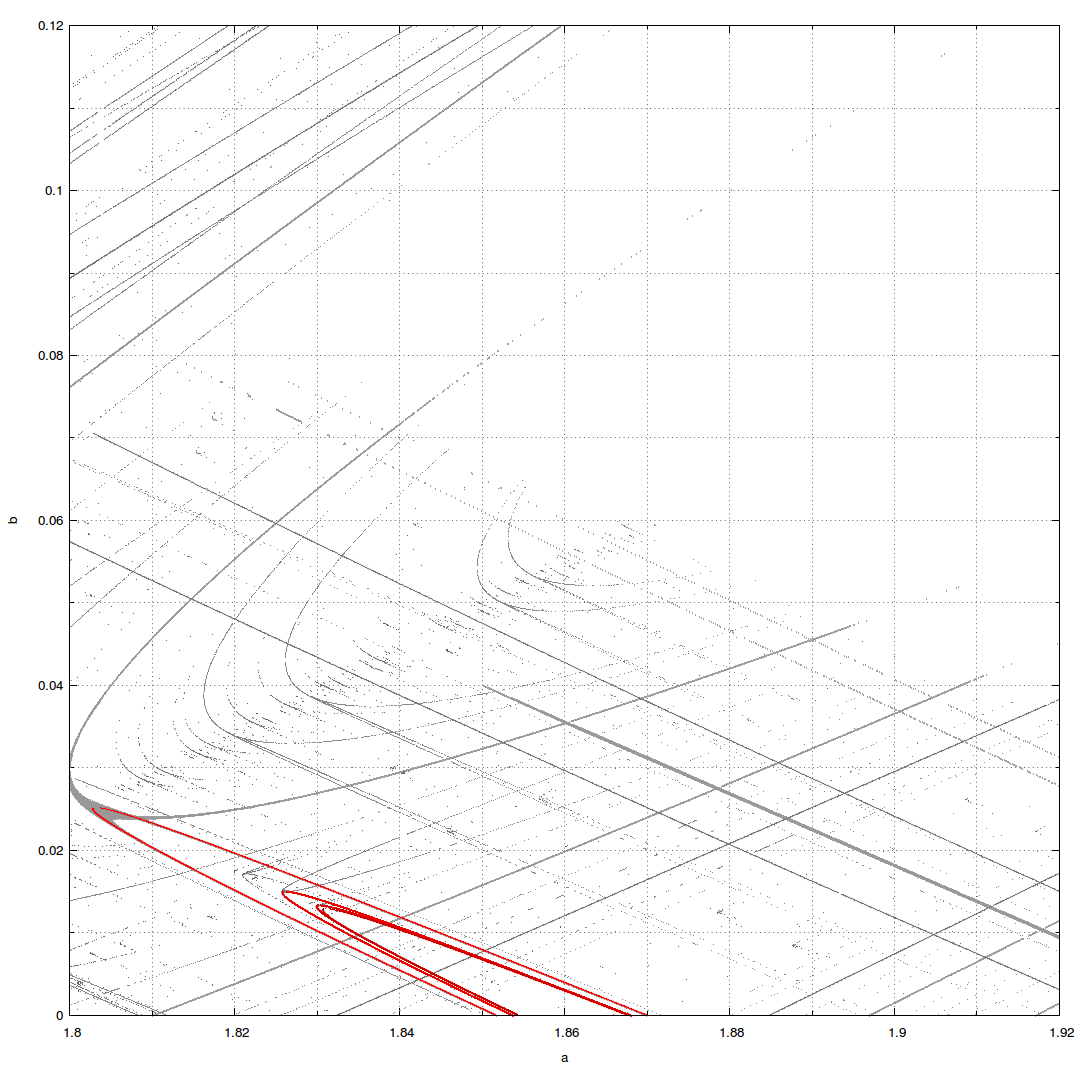}
\caption{The isotracal curves of periods $8,11,14$ and $17$ generated from head $1001C$ are shown in red. 
The grey region show the scatterplots for periods $8,11,14$ and $17$. 
The darker the colour, the higher the period.}
\label{fig:path+plot1001C}
\end{figure}

We now give numerical evidence suggesting A and B are true. 
For simplicity we represent unimodal combinatorial types by the itinerary of the critical point (with respect to the standard partition $I_0=[0,m)$, $I_C={m}$, $I_1=(m,p-1]$.)
Figures~\ref{fig:per8}--~\ref{fig:per11} showed plots of parameters for which $F_{a,b}$ possessed an attracting periodic orbit of a fixed period $p$. 
These regions connected distinct degenerate H\'enon parameters. 
Numerically we construct a curve between such parameters as follows.
First, take a unimodal permutation $\upsilon$ satisfying the hypotheses of Construction~\ref{constr:gencabling2}, 
i.e. it is  reconnectable at the non-dynamical preimage, and the non-dynamical preimage lies to the left of the folding point.
Call $\upsilon$ the {\it head}.
Apply Construction~\ref{constr:gencabling2}, giving a braid-equivalent pair $\upsilon_-^0$ and $\upsilon_+^0$.
Then apply Construction~\ref{constr:gencablingprocess1} inductively, giving braid-equivalent pairs $\upsilon_-^i$ and $\upsilon_+^i$ for $i=1,2,3$.
This was repeated for various $\upsilon$. 
See Table~\ref{table:BEbig2} for orbits listed by itinerary and grouped by head, with the associated permutation given in cyclic notation.
They are listed in pairs obtained as just described: the first pair is obtained from the head by Construction~\ref{constr:gencabling2} and subsequent ones from Construction~\ref{constr:gencablingprocess1}, inductively.

For each pair $\upsilon_-^i$ and $\upsilon_+^i$ we computed the superattracting parameters $a_{-}$ and $a_{+}$ in the quadratic family with critical orbit of type $\upsilon_{-}$ and $\upsilon_{+}$ respectively.
Then the parameter locus of
\begin{equation}\label{eq:Newton}
F_{\A,\B}^\P(\XINIT,\YINIT)=(\XINIT,\YINIT), \qquad \tr DF^\P_{\A,\B}(\XINIT,\YINIT)=0.
\end{equation}
passing through the parameters $(\A_{-},0)$ and $(\A_{+},0)$ was computed. 
We call a parameter curve lying in the locus~\eqref{eq:Newton} a {\it zero isotracal path}.
(More generally an {\it isotracal path} satisfies~\eqref{eq:Newton}, but with the trace set to some fixed constant instead of zero.)
We compute an isotracal curve iteratively by starting from the initial data given by 
\begin{equation}\label{eq:Newton-init}
\A=\A_\pm, \quad  
\B=0, \quad 
\XINIT=\A_\pm, \quad
\YINIT=0
\end{equation}
On each slice $\{\B=\B_0\}$, Newton's method was used to find $\A=\A(\B_0), \XINIT=\XINIT(\B_0)$ and $\YINIT=\YINIT(\B_0)$ satisfying equation~\eqref{eq:Newton}. 
%(Note that with $\B$ fixed the system~\eqref{eq:Newton} has three equations in three unknowns so the standard Newton map may be considered.)  
The value of $\B$ was then incremented and the values of $\A, \XINIT$ and $\YINIT$ from the previous step were used as initial data. 
The algorithm terminated once $\A_-(\B)$ and $\A_+(\B)$ were sufficiently close.

\begin{figure}[ht]
\begin{subfigure}[b]{0.5\linewidth}
\centering
\includegraphics
%[bb=0 0 1080 1080]
[scale=0.11]
{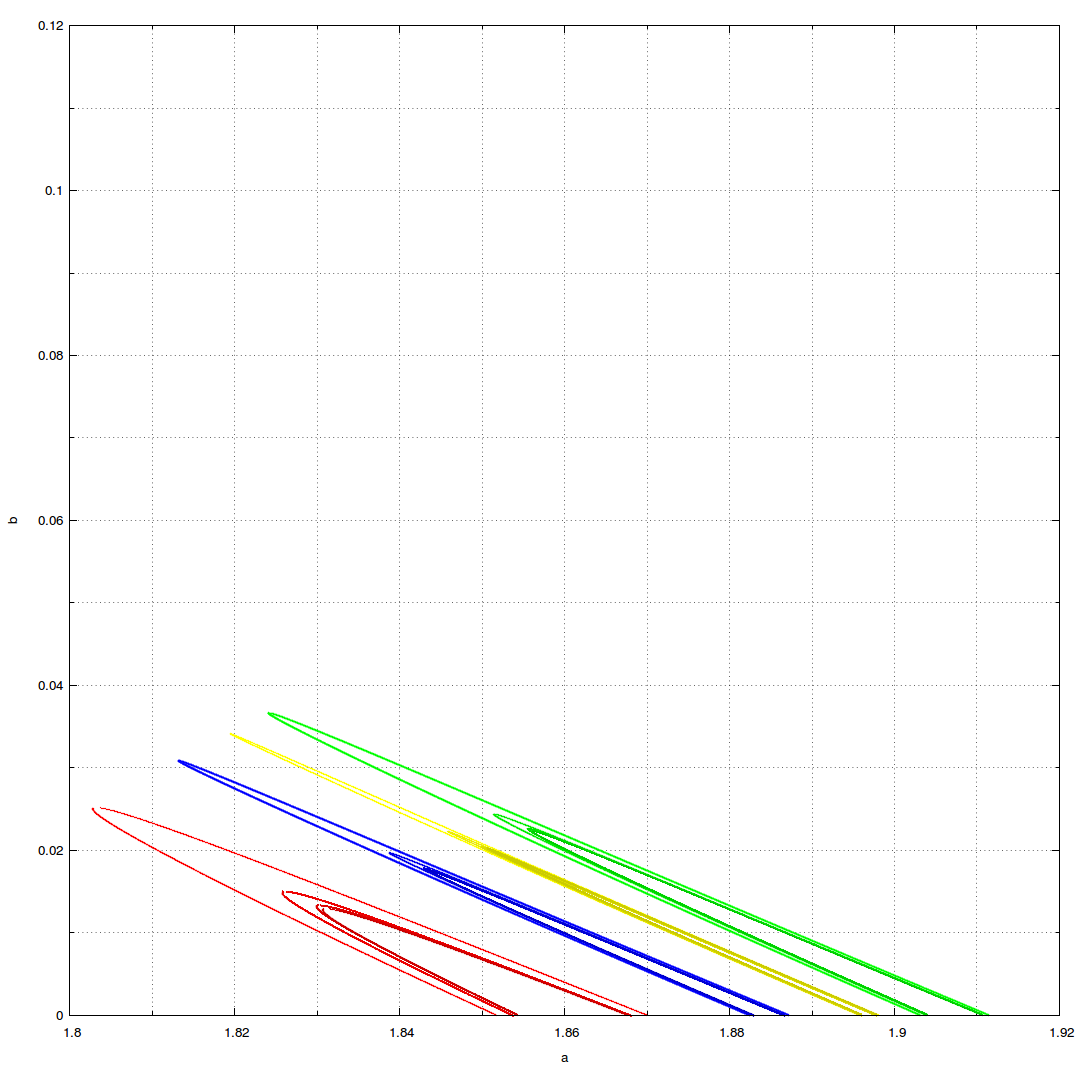}
\caption{Zero isotracal paths}
\label{fig:paths}
\end{subfigure}
~
\begin{subfigure}[b]{0.5\linewidth}
\centering
\includegraphics
%[bb=0 0 1080 1080]
[scale=0.11]
{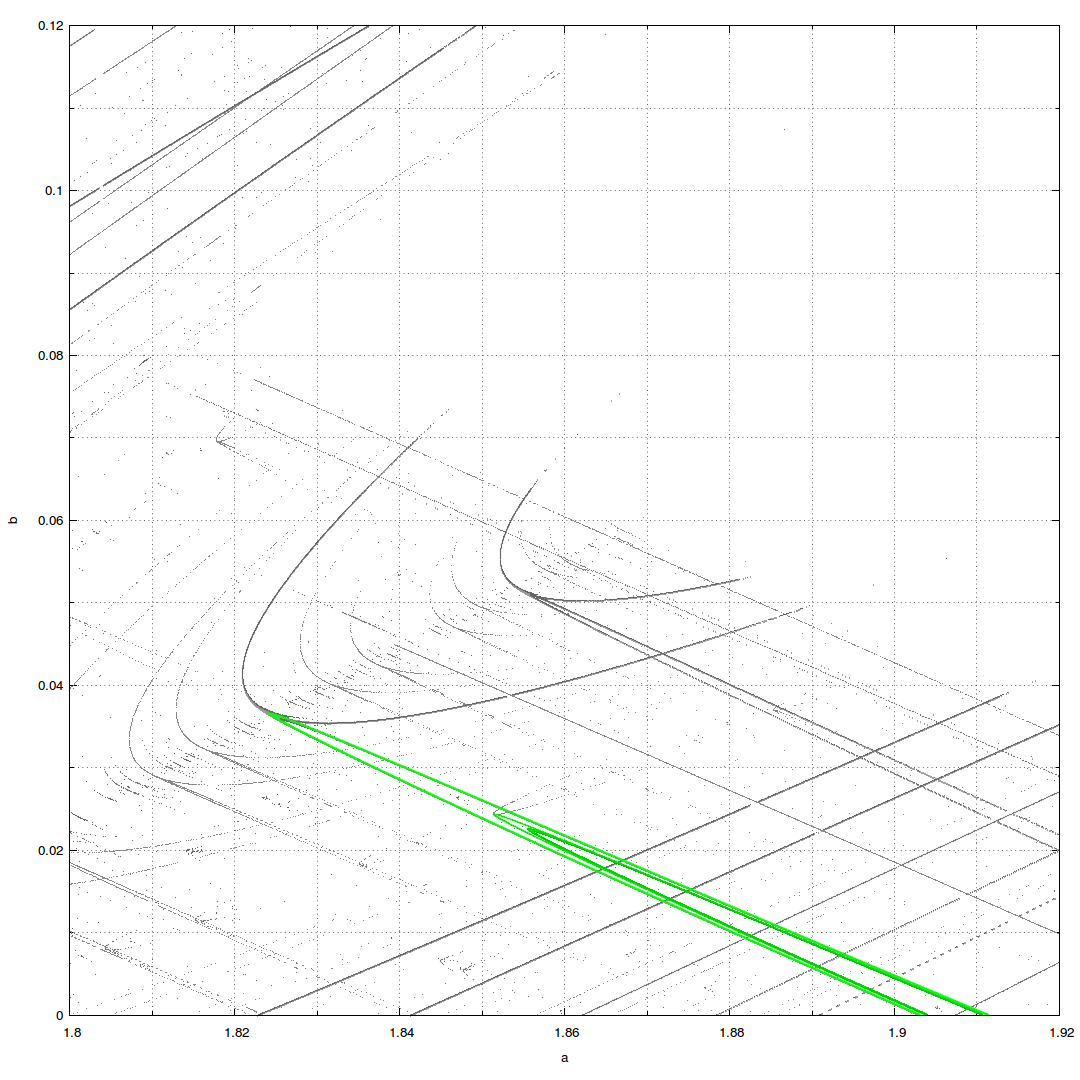}
\caption{Head $10011C$}
\label{fig:path+plot10011C}
\end{subfigure}

\vspace{0.2cm}
\begin{subfigure}[b]{0.5\linewidth}
\centering
\includegraphics
%[bb=0 0 1080 1080]
[scale=0.11]
{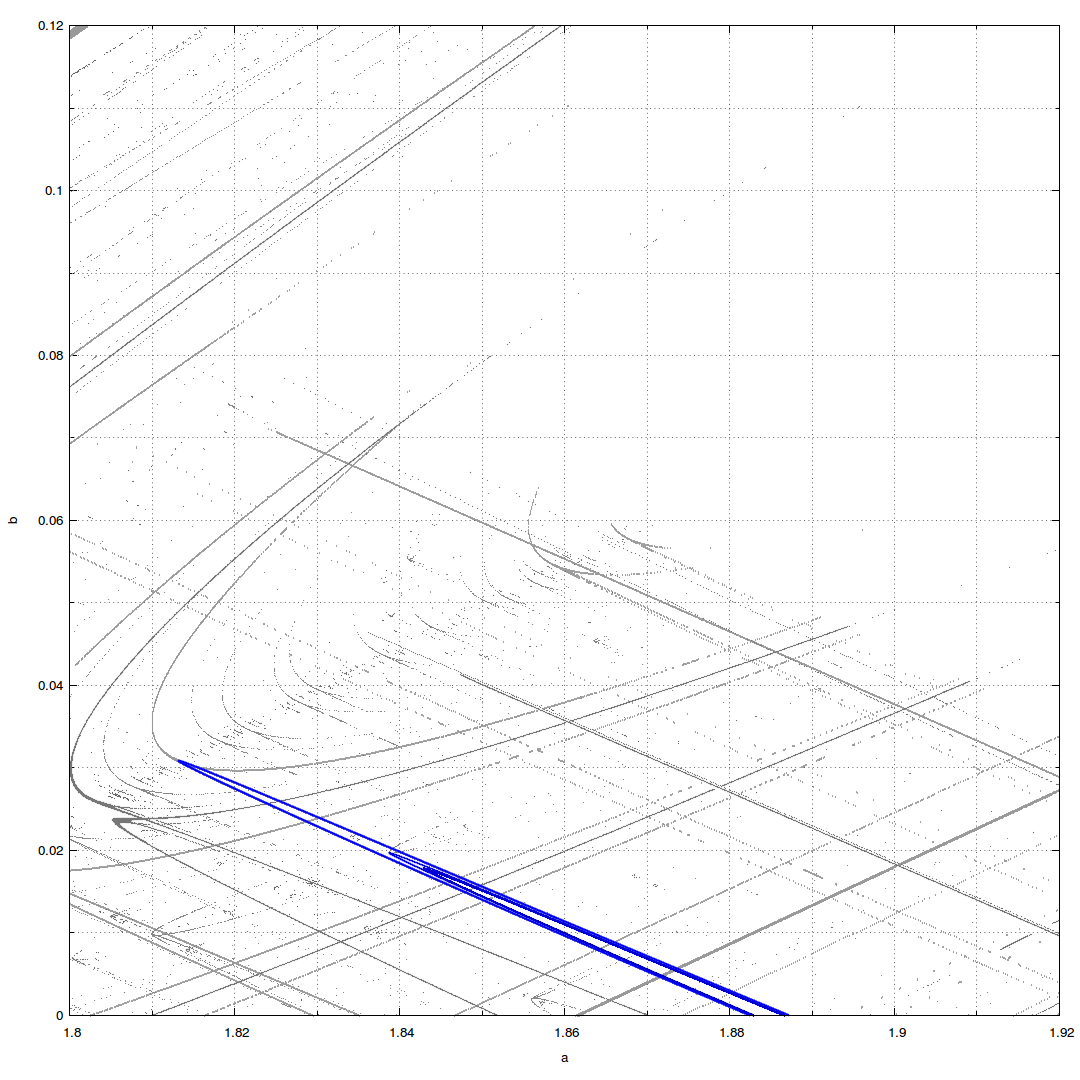}
\caption{Head $100111C$}
\label{fig:path+plot100111C}
\end{subfigure}
~
\begin{subfigure}[b]{0.5\linewidth}
\centering
\includegraphics
%[bb=0 0 1080 1080]
[scale=0.11]
{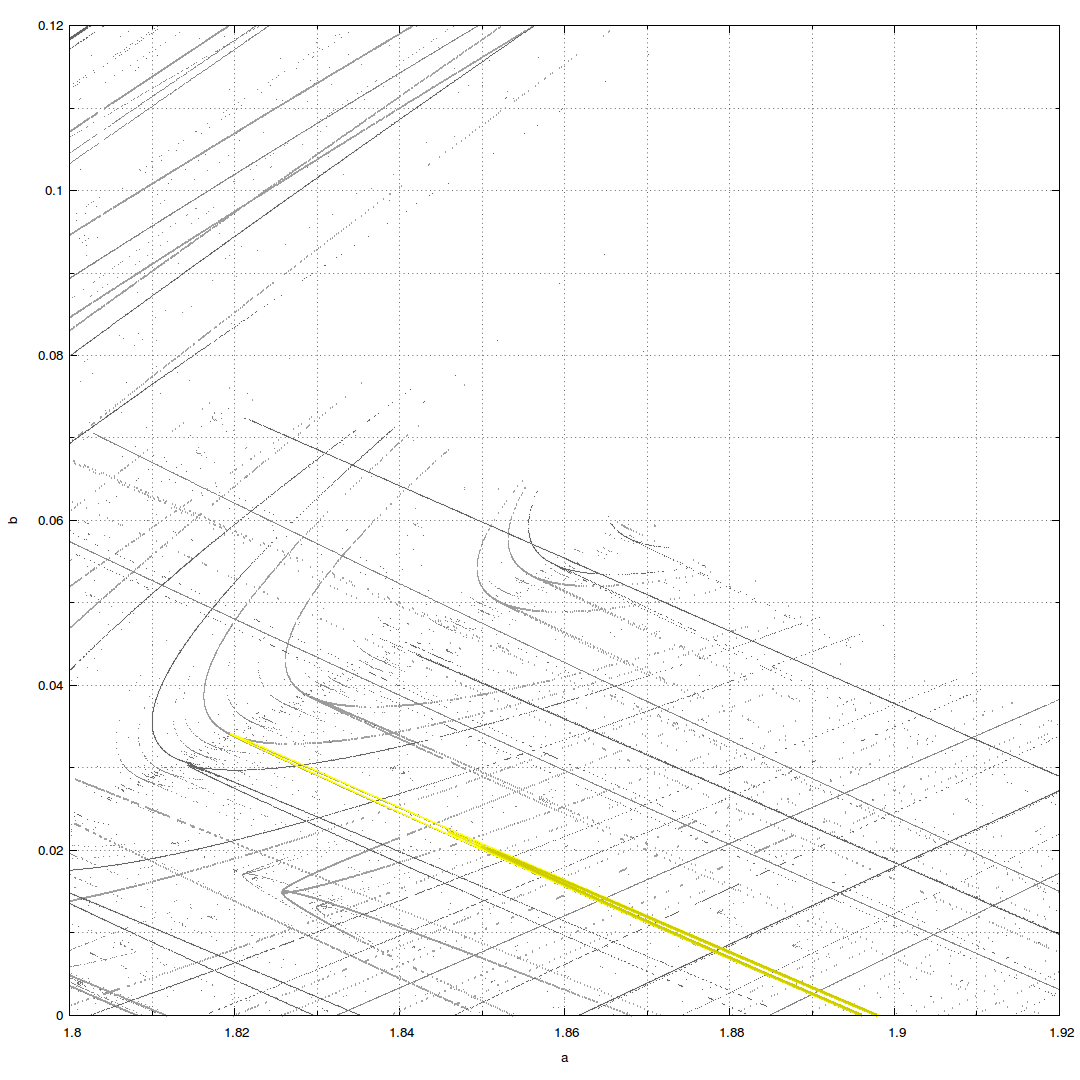}
\caption{Head $100111C$}
\label{fig:path+plot1001111C}
\end{subfigure}
\caption{For the H\'enon family $F_{a,b}$, in Figure~\ref{fig:paths} zero isotracal paths are plotted for heads $1001C$, $10011C$, $100111C$ and $1001111C$. 
In Figures~\ref{fig:path+plot1001C},~\ref{fig:path+plot10011C} and~\ref{fig:path+plot100111C} 
the same paths are plotted, for a fixed head, together with scatterplots of parameters with attracting points of the same period.}
\label{fig:isotracal}
\end{figure}

Table~\ref{table:BEbig} shows some of the braid equivalences which were realised in the H\'enon family. 
%For the readers convenience Table~\ref{table:BEbig2} is included showing the combinatorial type, both in terms of the itinerary of the critical orbit and in cyclic notation.
The paths in the parameter region of $(a,b)\in[1.8,1.9]\times [0.0,0.1]$ are shown in Figure~\ref{fig:paths}.
The red curves show the period $8$,$11$,$14$ and $17$ curves associated with head $1001C$ given in Table~\ref{table:BEbig}. They are shaded so that the darker the curve is, the higher the period.
Similarly, the green curves show the period $9$,$12$,$15$ and $18$ curves associated with the head $10011C$ from Table~\ref{table:BEbig}.
The blue curves show the period $10$,$13$,$16$ and $19$ curves associated with the head $100111C$ from Table~\ref{table:BEbig}.
The yellow curves show the period $11$,$14$,$17$ and $20$ curves associated with the head $1001111C$ from Table~\ref{table:BEbig}.
In each of these cases, the darker the curve is, the higher the period.

These plots were then superimposed with the scatterplots in Figures~\ref{fig:path+plot1001C}-\ref{fig:path+plot1001111C}. 
The grey region denoting the data from the first algorithm.
So, for example, Figure~\ref{fig:path+plot1001C} shows the curves of period $8$, $11$, $14$ and $17$ in Figure~\ref{fig:paths} together with the scatterplot data from the introduction for periods $8$, $11$, $14$, and $17$, where the darker the grey is the lower the period. Figures~\ref{fig:path+plot10011C} and~\ref{fig:path+plot100111C} are similar.

Table~\ref{table:BEbig} also gives the prefix and decoration to compare the current mechanism with that given by the first and second authors (see~\cite{dCH3} for more details on prefixes and decorations). 
For example consider the kneading sequences $10011110C$ and $10011010C$ (see Table~\ref{table:BEbig}). 
These have the same prefix $1001$ but different decorations, $110$ and $010$ respectively, showing that braid equivalent sequences constructed by the mechanism in this article may have different decorations.

%HENON-LIKE MAPS
We can weaken Questions A and B by asking if these braid equivalences are realisable in a more general class of maps containing the H\'enon family.
The typical generalisation of the H\'enon family is that of {\it H\'enon-like maps}. 
These are maps of the form
\begin{equation}
F(x,y)=(f(x)-\varepsilon(x,y),x)
\end{equation} 
where $f$ is unimodal on some interval $J$ and $\epsilon\colon J\times J\to \RR$ satisfies $\del_y\epsilon>0$.
These are diffeomorphisms onto their images which appear, after a suitable change of variables, when considering maps in the neighbourhood of a homoclinic bifurcation~\cite{PandT}.
\begin{quote}
{\it Question A'.}
Given unimodal maps $f_-$ and $f_+$ of type $\upsilon_-$ and $\upsilon_+$ respectively, does there exist a family $F_t$, $t\in [-1,1]$, of H\'enon-like diffeomorphisms such that $F_{-1}=\iota(f_-)$ and $F_{+1}=\iota(f_+)$, where $\iota$ is some embedding of the set of unimodal maps into the boundary of the space of H\'enon-like diffeomorphisms?
\end{quote}
and
\begin{quote}
{\it Question B'.}
For each positive integer $i$, let $f_{-,i}$ and $f_{+,i}$ have critical orbits $c_-^i$ and $c_+^i$ of types $\upsilon_-^i$ and $\upsilon_+^i$ respectively.
Let $C_-^i$ and $C_+^i$ denote the corresponding periodic orbits for the degenerate H\'enon-like maps $F_{-,i,0}$ and $F_{+,i,0}$ respectively.
Does there exist, for each $i$, a one-parameter family of H\'enon-like maps  $F_{i,t}$ realising a braid equivalence  connecting $(C_-^i,F_{-,i,0})$ and $(C_+^i,F_{+,i,0})$?
\end{quote}

\small
%\addcontentsline{toc}{section}{Bibliography}
\bibliographystyle{amsplain}
\bibliography{BEiH1-2015-06-07}

\providecommand{\bysame}{\leavevmode\hbox to3em{\hrulefill}\thinspace}
\providecommand{\MR}{\relax\ifhmode\unskip\space\fi MR }
% \MRhref is called by the amsart/book/proc definition of \MR.
\providecommand{\MRhref}[2]{%
  \href{http://www.ams.org/mathscinet-getitem?mr=#1}{#2}
}
\providecommand{\href}[2]{#2}
\begin{thebibliography}{10}

\bibitem{BirmanBook}
Joan~S. Birman, \emph{Braids, links, and mapping class groups}, Princeton
  University Press, Princeton, N.J., 1974, Annals of Mathematics Studies, No.
  82. \MR{0375281 (51 \#11477)}

\bibitem{Boyland84}
Philip Boyland, \emph{Braid types and a topological method of proving positive
  entropy}, Boston University, 1984.

\bibitem{Boyland94}
Philip Boyland, \emph{Topological methods in surface dynamics}, Topology Appl.
  \textbf{58} (1994), no.~3, 223--298. \MR{1288300 (95h:57016)}

\bibitem{dCH3}
Andr{\'e} de~Carvalho and Toby Hall, \emph{The forcing relation for horseshoe
  braid types}, Experimental Mathematics \textbf{11} (2002), no.~2, 271--288.

\bibitem{deCarvalhoHallHazard15}
Andr{\'e} de~Carvalho, Toby Hall, and Peter Hazard, \emph{Braid equivalence in
  the {H}{\'e}non family {II}}, in preparation, 2015.

\bibitem{DevaneyNitecki79}
R.~Devaney and Z.~Nitecki, \emph{Shift automorphisms in the {H}\'enon mapping},
  Comm. Math. Phys. \textbf{67} (1979), no.~2, 137--146. \MR{539548
  (80f:58035)}

\bibitem{ElHamoulyMira82}
Hassan El~Hamouly and Christian Mira, \emph{Singularit\'es dues au feuilletage
  du plan des bifurcations d'un diff\'eomorphisme bi-dimensionnel}, C. R. Acad.
  Sci. Paris S\'er. I Math. \textbf{294} (1982), no.~12, 387--390. \MR{659728
  (83m:58056)}

\bibitem{Hall94}
Toby Hall, \emph{The creation of horseshoes}, Nonlinearity \textbf{7} (1994),
  861--924.

\bibitem{VLHansen1}
Vagn~Lundsgaard Hansen, \emph{Braids and coverings: selected topics}, London
  Mathematical Society Student Texts, vol.~18, Cambridge University Press,
  Cambridge, 1989, With appendices by Lars G{\ae}de and Hugh R. Morton.
  \MR{1247697 (94g:57004)}

\bibitem{Holmes86}
Philip Holmes, \emph{Knotted periodic orbits in suspensions of {S}male's
  horseshoe: period multiplying and cabled knots}, Phys. D \textbf{21} (1986),
  no.~1, 7--41. \MR{860006 (88b:58112)}

\bibitem{KKY92}
Ittai Kan, H{\"u}seyin Ko{\c{c}}ak, and James~A. Yorke, \emph{Antimonotonicity:
  concurrent creation and annihilation of periodic orbits}, Ann. of Math. (2)
  \textbf{136} (1992), no.~2, 219--252. \MR{1185119 (94c:58135)}

\bibitem{Lyubich00}
Mikhail Lyubich, \emph{The quadratic family as a qualitatively solvable model
  of chaos}, Notices Amer. Math. Soc. \textbf{47} (2000), no.~9, 1042--1052.
  \MR{1777885 (2001g:37063)}

\bibitem{Milnor92}
John Milnor, \emph{Remarks on iterated cubic maps}, Experiment. Math.
  \textbf{1} (1992), no.~1, 5--24. \MR{1181083 (94c:58096)}

\bibitem{PandT}
Jacob Palis, Jr. and Floris Taken, \emph{Hyperbolicity \& sensitive chaotic
  dynamics at homoclinic bifurcations}, Cambridge Studies in Advanced
  Mathematics, vol.~35, Cambridge University Press, 1993.

\bibitem{Sannami88}
Atsuro Sannami, \emph{A topological classification of the periodic orbits of
  the {H}enon family}, Tech. report, Hokkaido University, 1988.

\end{thebibliography}

\begin{landscape}
\scriptsize

\begin{center}
\begin{longtable}{ll}

\caption[Braid Equivalences.]{Braid Equivalences in Cyclic Notation.}\label{table:BEbig2} \\

\hline \hline \\[-2ex] 
\multicolumn{1}{l}{critical itinerary} & 
\multicolumn{1}{l}{cyclic notation}  \\[0.5ex] \hline
\\[-1.8ex]
\endfirsthead

\multicolumn{2}{c}{\tablename}{ \thetable{} -- Continued} \\
\hline \hline \\[-2ex]
\multicolumn{1}{l}{critical itinerary} & 
\multicolumn{1}{l}{cyclic notation}  \\[0.5ex] \hline
\\[-1.8ex]
\endhead

\hline
\multicolumn{2}{r}{{Continued on Next Page\ldots}} \\
\endfoot
\\ \hline \hline
\endlastfoot

1001010C & 2,7,3,5,0,4,6,1\\
1001110C & 2,7,3,0,5,4,6,1\\
1001010010C & 2,7,10,3,8,5,0,4,9,6,1\\
1001110010C & 2,7,10,3,8,0,5,4,9,6,1\\
1001010010110C & 2,7,13,10,3,8,5,0,11,4,9,12,6,1\\
1001110010110C & 2,7,13,10,3,8,0,5,11,4,9,12,6,1\\
1001010010110110C & 2,7,16,13,10,3,8,5,0,14,11,4,9,12,15,6,1\\
1001110010110110C & 2,7,16,13,10,3,8,0,5,14,11,4,9,12,15,6,1\\ 
\hline
10011110C & 2,8,3,0,6,4,5,7,1\\
10011010C & 2,8,3,6,0,4,5,7,1\\
10011110010C & 2,8,11,3,9,0,6,4,5,10,7,1\\
10011010010C & 2,8,11,3,9,6,0,4,5,10,7,1\\
10011110010110C & 2,8,14,11,3,9,0,6,12,4,5,10,13,7,1\\
10011010010110C & 2,8,14,11,3,9,6,0,12,4,5,10,13,7,1\\
10011010010110110C & 2,8,17,14,11,3,9,6,0,15,12,4,5,10,13,16,7,1\\
10011110010110110C & 2,8,17,14,11,3,9,0,6,15,12,4,5,10,13,16,7,1\\ 
\hline
100111010C & 2,9,3,7,0,5,4,6,8,1\\
100111110C & 2,9,3,0,7,5,4,6,8,1\\
100111010010C & 2,9,12,3,10,7,0,5,4,6,11,8,1\\
100111110010C & 2,9,12,3,10,0,7,5,4,6,11,8,1\\
100111010010110C & 2,9,15,12,3,10,7,0,13,5,4,6,11,14,8,1\\
100111110010110C & 2,9,15,12,3,10,0,7,13,5,4,6,11,14,8,1\\
100111010010110110C & 2,9,18,15,12,3,10,7,0,16,13,5,4,6,11,14,17,8,1\\
100111110010110110C & 2,9,18,15,12,3,10,0,7,16,13,5,4,6,11,14,17,8,1\\ 
\hline
1001111110C & 2,10,3,8,0,6,4,5,7,9,1\\
1001111010C & 2,10,3,8,6,0,4,5,7,9,1\\
1001111110010C & 2,10,13,3,11,0,8,6,4,5,7,12,9,1\\
1001111010010C & 2,10,13,3,11,8,0,6,4,5,7,12,9,1\\
1001111110010110C & 2,10,16,13,3,11,0,8,14,6,4,5,7,12,15,9,1\\
1001111010010110C & 2,10,16,13,3,11,8,0,14,6,4,5,7,12,15,9,1\\
1001111110010110110C & 2,10,19,16,13,3,11,0,8,17,14,6,4,5,7,12,15,18,9,1\\
1001111010010110110C & 2,10,19,16,13,3,11,8,0,17,14,6,4,5,7,12,15,18,9,1\\ 
\hline
100110010100110C & 2,11,6,15,3,12,7,9,0,4,13,8,14,5,10,1\\
100110011100110C & 2,11,6,15,3,12,7,0,9,4,13,8,14,5,10,1\\
1001100101001100100110C & 2,11,18,6,15,22,3,12,19,7,16,9,0,4,13,20,8,21,14,5,17,10,1\\
1001100111001100100110C & 2,11,18,6,15,22,3,12,19,7,16,0,9,4,13,20,8,21,14,5,17,10,1\\
10011001010011001001101100110C & 2,11,25,18,6,15,29,22,3,12,26,19,7,16,9,0,23,4,13,27,20,8,21,28,14,5,17,24,10,1\\
10011001110011001001101100110C & 2,11,25,18,6,15,29,22,3,12,26,19,7,16,0,9,23,4,13,27,20,8,21,28,14,5,17,24,10,1\\
100110010100110010011011001101100110C & 2,11,32,25,18,6,15,36,29,22,3,12,33,26,19,7,16,9,0,30,23,4,13,34,27,20,8,21,28,35,14,5,17,24,31,10,1\\
100110011100110010011011001101100110C & 2,11,32,25,18,6,15,36,29,22,3,12,33,26,19,7,16,0,9,30,23,4,13,34,27,20,8,21,28,35,14,5,17,24,31,10,1
\end{longtable}

\begin{longtable}{ccclll}

\caption[Braid Equivalences.]{Braid Equivalences and Associated Data.}\label{table:BEbig} \\

\hline \hline \\[-2ex]
\multicolumn{1}{l}{head} &
\multicolumn{1}{c}{period} & 
\multicolumn{1}{c}{unimodal $a$} & 
\multicolumn{1}{l}{critical itinerary} & 
\multicolumn{1}{l}{prefix} & 
\multicolumn{1}{l}{decoration} \\[0.5ex] \hline
\\[-1.8ex]
\endfirsthead

\multicolumn{3}{c}{\tablename}{ \thetable{} -- Continued} \\
\hline \hline \\[-2ex]
\multicolumn{1}{l}{head} & 
\multicolumn{1}{c}{period} & 
\multicolumn{1}{c}{unimodal $a$} & 
\multicolumn{1}{l}{critical itinerary} & 
\multicolumn{1}{l}{prefix} & 
\multicolumn{1}{l}{decoration} \\[0.5ex] \hline
\\[-1.8ex]
\endhead

\hline
\multicolumn{3}{r}{{Continued on Next Page\ldots}} \\
\endfoot
\\ \hline \hline
\endlastfoot
1001C 
&8 & 1.85173004941 & 1001010C & 1001 & 10 \\
%&&& 2,7,3,5,0,\fix,4,6,1 &&\\
&8 & 1.87000388083 & 1001110C & 1001 & 10 \\
%&&& 2,7,3,0,5,\fix,4,6,1 &&\\
&11 & 1.85376146047 & 1001010010C & 1001 & 10010 \\
%&&& 2,7,10,3,8,5,0,\fix,4,9,6,1 &&\\
&11 & 1.86808014899 & 1001110010C & 1001 & 10010 \\
%&&& 2,7,10,3,8,0,5,\fix,4,9,6,1 &&\\
&14 & 1.85420779956 & 1001010010110C & 1001 & 10010110 \\
%&&& 2,7,13,10,3,8,5,0,11,\fix,4,9,12,6,1 &&\\
&14 & 1.86768154112 & 1001110010110C & 1001 & 10010110 \\
%&&& 2,7,13,10,3,8,0,5,11,\fix,4,9,12,6,1 &&\\
&17 & 1.85429374549 & 1001010010110110C & 1001 & 10010110110 \\
%&&& 2,7,16,13,10,3,8,5,0,14,11,\fix,4,9,12,15,6,1 &&\\
&17 & 1.86760846133 & 1001110010110110C & 1001 & 10010110110 \\
%&&& 2,7,16,13,10,3,8,0,5,14,11,\fix,4,9,12,15,6,1 &&\\ 
\hline
10011C 
&9 & 1.90311677305 & 10011110C & 1001 & 110 \\
%&&& 2,8,3,0,6,4,\fix,5,7,1 &&\\
&9 & 1.91144463147 & 10011010C & 1001 & 010 \\
%&&& 2,8,3,6,0,4,\fix,5,7,1 &&\\
&12 & 1.90383983136 & 10011110010C & 1001 & 110010 \\
%&&& 2,8,11,3,9,0,6,4,\fix,5,10,7,1 &&\\
&12 & 1.91074213574 & 10011010010C & 1001 & 010010 \\
%&&& 2,8,11,3,9,6,0,4,\fix,5,10,7,1 &&\\
&15 & 1.90396328656 & 10011110010110C & 1001 & 110010110 \\
%&&& 2,8,14,11,3,9,0,6,12,4,\fix,5,10,13,7,1 &&\\
&15 & 1.91062464072 & 10011010010110C & 1001 & 010010110 \\
%&&& 2,8,14,11,3,9,6,0,12,4,\fix,5,10,13,7,1 &&\\
&18 & 1.90398296779 & 10011010010110110C & 1001 & 010010110110 \\
%&&& 2,8,17,14,11,3,9,6,0,15,12,4,\fix,5,10,13,16,7,1 &&\\
&18 & 1.91060625301 & 10011110010110110C & 1001 & 110010110110 \\
%&&& 2,8,17,14,11,3,9,0,6,15,12,4,\fix,5,10,13,16,7,1 &&\\ 
\hline
100111C
&10 & 1.88240793375 & 100111010C & 1001 & 1010 \\
%&&& 2,9,3,7,0,5,\fix,4,6,8,1 &&\\
&10 & 1.88717220640 & 100111110C & 1001 & 1110 \\
%&&& 2,9,3,0,7,5,\fix,4,6,8,1 &&\\
&13 & 1.88286957759 & 100111010010C & 1001 & 1010010 \\
%&&& 2,9,12,3,10,7,0,5,\fix,4,6,11,8,1 &&\\
&13 & 1.88672860609 & 100111110010C & 1001 & 1110010 \\
%&&& 2,9,12,3,10,0,7,5,\fix,4,6,11,8,1 &&\\
&16 & 1.88295631360 & 100111010010110C & 1001 & 1010010110 \\
%&&& 2,9,15,12,3,10,7,0,13,5,\fix,4,6,11,14,8,1 &&\\
&16 & 1.88664609851 & 100111110010110C & 1001 & 1110010110 \\
%&&& 2,9,15,12,3,10,0,7,13,5,\fix,4,6,11,14,8,1 &&\\
&19 & 1.88297117918 & 100111010010110110C & 1001 & 1010010110110 \\
%&&& 2,9,18,15,12,3,10,7,0,16,13,5,\fix,4,6,11,14,17,8,1 &&\\
&19 & 1.88663211208 & 100111110010110110C & 1001 & 1110010110110 \\
%&&& 2,9,18,15,12,3,10,0,7,16,13,5,\fix,4,6,11,14,17,8,1 &&\\ 
\hline
1001111C
&11 & 1.89575422001 & 1001111110C & 1001 & 11110 \\
%&&& 2,10,3,8,0,6,4,\fix,5,7,9,1 &&\\
&11 & 1.89808902792 & 1001111010C & 1001 & 11010 \\
%&&& 2,10,3,8,6,0,4,\fix,5,7,9,1 &&\\
&14 & 1.89596323393 & 1001111110010C & 1001 & 11110010 \\
%&&& 2,10,13,3,11,0,8,6,4,\fix,5,7,12,9,1 &&\\
&14 & 1.89787945132 & 1001111010010C & 1001 & 11010010 \\
%&&& 2,10,13,3,11,8,0,6,4,\fix,5,7,12,9,1 &&\\
&17 & 1.89600029216 & 1001111110010110C & 1001 & 11110010110 \\
%&&& 2,10,16,13,3,11,0,8,14,6,4,\fix,5,7,12,15,9,1 &&\\
&17 & 1.89784256570 & 1001111010010110C & 1001 & 11010010110 \\
%&&& 2,10,16,13,3,11,8,0,14,6,4,\fix,5,7,12,15,9,1 &&\\
&20 & 1.89600636282 & 1001111110010110110C & 1001 & 11110010110110 \\
%&&& 2,10,19,16,13,3,11,0,8,17,14,6,4,\fix,5,7,12,15,18,9,1 &&\\
&20 & 1.89783655715 & 1001111010010110110C & 1001 & 11010010110110 \\
%&&& 2,10,19,16,13,3,11,8,0,17,14,6,4,\fix,5,7,12,15,18,9,1 &&\\ 
\hline
10011001C
&16 & 1.93213488504 & 100110010100110C & 10011001 & 100110 \\
%&&& 2,11,6,15,3,12,7,9,0,4,13,\fix,8,14,5,10,1 &&\\
&16 & 1.93235574679 & 100110011100110C & 10011001 & 100110 \\
%&&& 2,11,6,15,3,12,7,0,9,4,13,\fix,8,14,5,10,1 &&\\
&23 & 1.93213896428 & 1001100101001100100110C & 10011001 & 1001100100110 \\
%&&& 2,11,18,6,15,22,3,12,19,7,16,9,0,4,13,20,\fix,8,21,14,5,17,10,1 &&\\
&23 & 1.93235156026 & 1001100111001100100110C & 10011001 & 1001100100110 \\
%&&& 2,11,18,6,15,22,3,12,19,7,16,0,9,4,13,20,\fix,8,21,14,5,17,10,1 &&\\
&30 & 1.93213911573 & 10011001010011001001101100110C & 10011001 & 10011001001101100110 \\
%&&& 2,11,25,18,6,15,29,22,3,12,26,19,7,16,9,0,23,4,13,27,20,\fix,8,21,28,14,5,17,24,10,1 &&\\
&30 & 1.93235140816 & 10011001110011001001101100110C & 10011001 & 10011001001101100110 \\
%&&& 2,11,25,18,6,15,29,22,3,12,26,19,7,16,0,9,23,4,13,27,20,\fix,8,21,28,14,5,17,24,10,1 &&\\
&37 & 1.93213912112 & 100110010100110010011011001101100110C & 10011001 & 100110010011011001101100110 \\
%&&& 2,11,32,25,18,6,15,36,29,22,3,12,33,26,19,7,16,9,0,30,23,4,13,34,27,20,\fix,8,21,28,35,14,5,17,24,31,10,1 &&\\
&37 & 1.93235140286 & 100110011100110010011011001101100110C & 10011001 & 100110010011011001101100110 \\
%&&& 2,11,32,25,18,6,15,36,29,22,3,12,33,26,19,7,16,0,9,30,23,4,13,34,27,20,\fix,8,21,28,35,14,5,17,24,31,10,1 &&
\end{longtable}

\end{center}
\end{landscape}

\end{document}